\documentclass[11pt]{amsart}
\usepackage{amsthm}
\usepackage{amsmath}
\usepackage{amstext}
\usepackage{amsfonts}
\usepackage{amssymb}
\usepackage{latexsym}
\usepackage{amscd}
\usepackage{xypic}
\theoremstyle{plain}
\newtheorem{theorem}{Theorem}[section]

\newtheorem{proposition}[theorem]{Proposition}
\newtheorem{corollary}[theorem]{Corollary}

\theoremstyle{definition}
\newtheorem{definition}[theorem]{Definition}

\newtheorem{problem}[theorem]{Problem}

\newtheorem{example}[theorem]{Example}
\newtheorem{remark}[theorem]{Remark}
\newtheorem{remarks}[theorem]{Remarks}

\numberwithin{equation}{theorem}

\newtheorem{exercise}[theorem]{Exercise}
\newtheorem*{acknowledgement}{Acknowledgement}
\newcommand{\quadand}{\quad\text{and}\quad}
\newcommand{\lieg}{\operatorname{\mathfrak{g}}}
\begin{document}
\title[Group schemes and surfaces in positive characteristic]{Topics in group schemes and surfaces in positive characteristic}

\author[]{Nikolaos Tziolas}
\address{Department of Mathematics and Statistics, University of Cyprus, Nicosia 1678, Cyprus}
%\curraddr{}
\email{tziolas@ucy.ac.cy}


\thanks{From February 2023 to August 2023 the author was supported by the Leverhulme Visiting Professorship VP1-2022-004.}
\thanks{The paper is for the final Edge volume in the ADUF journal.}

%\dedicatory{\mydate}

\begin{abstract}
This is a survey paper on algebraic surfaces in positive characteristic based on a series of lectures that I gave at the University of Edinburgh in March 2023. It is focused on certain positive characteristic phenomena like infinitesimal group schemes and their actions on algebraic surfaces as well as the failure in positive characteristic of certain fundamental characteristic zero results like the Kodaira vanishing theorem. Many explicit examples are presented.
\end{abstract}

\maketitle
\tableofcontents

These notes are based on a series of lectures that I gave in March 2023 at the University of Edinburgh. The main purpose of writing them is in order to present certain differences between the geometry of algebraic varieties in characteristic zero and in characteristic $p>0$ and to familiarise the readers with certain positive characteristic techniques. 

Geometry in positive characteristic is not only of interest for its own sake but it is also strongly connected to geometry in characteristic zero. In fact many fundamental results in characteristic zero, like Mori's theorem of the rational connectedness of a smooth Fano variety in characteristic zero, are proved via characteristic $p>0$ methods.

There are many differences between geometry in characteristic zero and in positive characteristic. In particular, many fundamental vanishing theorems that hold in characteristic zero (like the Kodaira vanishing and the Akizuki-Nakano-Kodaira vanishing theorem), as well as fundamental inequalities like the Bogomolov-Miyaoka-Yau inequality, fail to hold in general in positive characteristic. The automorphism group scheme and the Picard scheme may not be reduced and bounds for their size that hold in characteristic zero do not hold in positive characteristic (see Section~\ref{automorphisms}). Many of the positive characteristic "pathologies" were investigated by D. Mumford in a series of papers~\cite{Mu61},~\cite{Mu62},~\cite{Mu67}.

I did not intend to present a wide range of topics, an excellent reference for which is the paper of C. Liedtke~\cite{Liedtke13}, but rather to emphasise on  some particular topics that I considered of special interest, some of them directly related to my own personal research. 

I tried to  present in some detail theory that exists exclusively in positive characteristic, give references for further reading and to exhibit through explicit examples certain fundamental differences between the characteristic zero and the positive characteristic case. In addition, I provided proofs, or the main parts of them,  of selected results  with the aim to familiarise the reader with certain techniques of the subject that might be useful to someone who intends to work in these areas.

The main topics presented here are infinitesimal group schemes and their actions on algebraic varieties, the structure of the automorphism group scheme of an algebraic variety and the failure of the Kodaira Vanishing Theorem in positive characteristic. 

The material presented in these notes is organised as follows.

\textbf{In section 2}  some basic results about purely inseparable field extensions and morphisms of schemes are presented. The notion of a restricted Lie algebra is defined and  the Jacobson correspondence between purely inseparable field extensions of height $1$ and restricted Lie algebras of derivations is presented. The section concludes with generalizations to schemes.

\textbf{In section 3} some fundamental results about group schemes in positive characteristic are presented and discussed. Some topics are. The Lie algebra of a group scheme. Non smooth group schemes and other unusual properties of group schemes in positive characteristic. Infinitesimal group schemes. The Demazure-Gabriel correspondence between height $1$ group schemes and restricted Lie algebras. Actions of infinitesimal group schemes, and in particular of $\alpha_p$ and $\mu_p$,  on algebraic varieties. The quotient of a variety by an infinitesimal group scheme. The singularities of the quotient of a smooth surface by a nontrivial $\alpha_p$ or $\mu_p$ action. Many examples are also given.

\textbf{In section 4} the structure of the automorphism group scheme of a smooth variety of general type in positive characteristic is discussed. Existing results in characteristic zero about its size and structure are reviewed and it is shown by explicit examples that almost all fail in positive characteristic. Finally some structure results in positive characteristic are presented.

\textbf{In section 5} the failure of the Kodaira Vanishing theorem in positive characteristic is discussed. Conditions under which it holds are presented and an explicit example of a smooth projective surface of general type in any characteristic  is given where both the Kodaira Vanishing Theorem and the Bogomolov-Miyaoka-Yau inequality fail.

\begin{acknowledgement}
From February 2023 until August 2023 I was a Leverhulme Visiting Professor (grant reference VP1-2022-004)  at the School of Mathematics of the University of Edinburgh. The material presented in these notes is based on a series of lectures that I gave at the University of Edinburgh in March 2023. I would like to thank the Leverhulme Trust for its valuable support and for giving me the opportunity to visit the University of Edinburgh.
\end{acknowledgement}

\section{Frobenius and  purely inseparable maps.}\label{sec-1}
\subsection{The Frobenius map.}

\begin{definition}
Let $\pi \colon X \rightarrow \mathrm{Spec} k$ be a scheme defined over a field $k$ of characteristic $p>0$. The Frobenius morphism $F \colon X \rightarrow X$ is the map of ringed spaces which is the identity on the topological space $X$ and for any $U \subset X$ open, $F(U) \colon \mathcal{O}_X(U)\rightarrow \mathcal{O}_X(U)$ is given by $F(U)(a)=a^p$, for any $a\in \mathcal{O}_X(U)$.
\end{definition}
The Frobenius map $F\colon X \rightarrow X$, and more generally its $n$-iterated product $F^n$,  is a map of ringed spaces but it is not a map of schemes over $k$. To remedy this the geometric Frobenius map $F^{(p^n)}$ is defined from the following diagram.
\[
\xymatrix{
X.\ar[drr]^{F^n}  \ar[ddr]^{\pi} \ar[dr]^{F^{(p^n)}} &.  &\\
&    X^{(p^n)}\ar[d] \ar[r]    &         X \ar[d]^{\pi} & \\
& \mathrm{Spec} k \ar[r]^{F^n} & \mathrm{Spec} k                                                 
}
\]
Where $F^n\colon \mathrm{Spec}k \rightarrow \mathrm{Spec} k$ is the $n$-iterated Frobenius map, $X^{(p^n)}$ is the fiber product and $F^{(p^n)}\colon X \rightarrow X^{(p^n)}$ is the map induced from the universal property of the fiber product. This is a $k$-map. 

Suppose that $k$ is a perfect field. Then the  Frobenius map $F\colon k \rightarrow k$ is an isomorphism and therefore $X^{(p^n)}$ is isomorphic to $X$ as abstract schemes but not necessarily as $k$-schemes as shown by Example~\ref{example Frobenius}. Locally,  if $X=\mathrm{Spec}A$, where $A$ is a $k$-algebra,   $X^{(p^n)}=\mathrm{Spec}A$ as a ringed space and its $k$-algebra structure is given by $\lambda \cdot a=\lambda^{1/p^n}a$, for any $a\in A$ and $\lambda \in k$. 

In particular, a straight forward chase through the definitions shows that $(\mathbb{A}^n_k)^{(p^n)}$ is isomorphic to $\mathbb{A}^n_k$ as $k$-schemes and the geometric Frobenius map 
\[
F^{(p^n)}\colon \mathbb{A}^n_k\rightarrow \mathbb{A}^n_k
\]
is induced from the the homomorphism of $k$-algebras $\Phi\colon k[x_1,\ldots, x_n]\rightarrow k[x_1,\ldots, x_n]$ given by $\Phi(f(x_1,\ldots,x_n))=f(x_1^{p^n},\ldots, x_n^{p^n})$.

\begin{example}\label{example Frobenius}
Let $k$ be a perfect field and $E$ be the elliptic curve over $k$ given by $y^2=x^3+ax+b$, $ab\not= 0$. Then $E^{(p)}$ is the elliptic curve given by the equation $y^2=x^3+a^px+b^p$. The $j$-invariant of $E$ is $6912a^3/(4a^3+27b^2)$ and the $j$-invariant of $E^{(p)}$ is $6912a^{3p}/(4a^{3p}+27b^{2p})$. Therefore, $E$ and $E^{(p)}$ are not isomorphic as $k$-schemes, even though they are isomorphic as abstract schemes.
\end{example}

If $X$ is a variety of dimension $n$ over a perfect field $k$ of Characteristic $p$ then the absolute and the geometric Frobenius are finite morphisms of degree $p^n$.

Finally, many properties of the singularities of $X$ are determined by properties of the Frobenius endomorphism of $X$. In particular, $X$ is smooth if and only the Frobenius map $F\colon X \rightarrow X$ is flat~\cite{Kunz69}. A very informative reference for the relation between the Frobenius and the singularities of a variety $X$ is the book by Karen Smith and Karl Schwede~\cite{SwSm24}.

\subsection{Purely inseparable field extensions.}

\begin{definition}
Let $E$ be a field of characteristic $p>0$. A field extension $E\subset F$ is called purely inseparable if and only if for any $a\in F$, there exists $n>0$ such that $a^{p^n}\in E$. The least such $n>0$ is called the height of $a$ and is denoted by $h(a)$. 

Suppose that $E\subset F$ is finite. Then there exists $r>0$ such that $F^{p^r} \subset E$. The least such $r>0$ is called the height of $F$ over $E$. 
\end{definition}

\begin{example}
Let $E$ be a field of characteristic $p>0$, $a\in E-E^p$ and $n>0$. The polynomial $x^{p^n}-a$ is then irreducible over $E$~\cite{Lang02} and therefore $F=E[x]/(x^{p^n}-a)$ is a field. The extension $E\subset F$ is purely inseparable of height $n$ and degree $p^n$.

More generally, let $a_1,\ldots, a_n \in E$ and $e_1, \ldots, e_n \in \mathbb{N}$ be positive integers such that $a_i \not\in k^p [a_1^{p^{e_1-1}}, \ldots, a_{i-1}^{p^{e_{i-1}-1}}]$, for all $i=1,\ldots, n$. Then the tensor product
\[
F=\frac{k[x_1]}{(x_1^{p^{e_1}}-a_1)}\otimes_k\cdots \otimes_k \frac{k[x_1]}{(x_n^{p^{e_n}}-a_n)}
\] is a field and $E\subset F$ is a purely inseparable extension. An extension of this form is called modular~\cite{Sweedler68}. All extensions of height 1 are modular~\cite{Sweedler68} but not every purely inseparable extension is as shown by the following example of Sweendler~\cite{Sweedler68}.
\end{example}
\begin{example}
Let $k=\mathbb{F}_p(a,b,c)$ be the field of rational functions in three variables of $\mathbb{F}_p$. Let $E=k[z,xz-y]$, where $z^{p^2}=c$, $x^p=a$ and $y^p=b$. This is an extension of height 2 which is not modular.
\end{example}

\subsection{Restricted Lie algebras.}
Let $k$ be a field. A Lie algebra over $k$ is a vector space $\mathfrak{g}$ over $k$ together with a $k$-bilinear map $\mathfrak{g} \times \mathfrak{g}\rightarrow \mathfrak{g}$ denoted by $(x,y)\mapsto [x,y]$ such that
\begin{enumerate}
\item $[x,x]=0$, for all $ x\in \mathfrak{g}$.
\item $[x,[y,z]]+[y,[z.x]]+[z,[x,y]]=0.$
\end{enumerate}
Let $x\in \mathfrak{g}$. The $k$-linear map $\mathrm{ad}_x \colon \mathfrak{g}\rightarrow \mathfrak{g}$ defined by $\mathrm{ad}_x(y)=[x,y]$ is called the adjoint endomorphism of $\mathfrak{g}$. 

Suppose that $k$ has characteristic $p>0$. A \emph{restricted Lie algebra} over $k$ is a Lie algebra $\mathfrak{g}$ over $k$, together with a map $\mathfrak{g}\rightarrow \mathfrak{g}$, $x\mapsto x^{[p]}$,
called the \emph{$p$-map}, subject to the following three axioms: 
\begin{enumerate}
\item We have $\mathrm{ad}_{x^{[p]}}= (\mathrm{ad}_x)^p$ for all vectors $x\in \mathfrak{g}$.
%\smallskip
\item Moreover $(\lambda\cdot x)^{[p]} = \lambda^p\cdot x^{[p]}$ for all  vectors $x\in \mathfrak{g}$ and scalars $\lambda\in k$.
%\smallskip
\item The formula $(x+y)^{[p]}= x^{[p]} + y^{[p]} + \sum_{r=1}^{p-1}s_r(x,y)$ holds for all $x,y\in\mathfrak{g}$.
\end{enumerate}
Here the summands $s_r(x,y)$ are   universal expressions defined by 
$$
s_r(x,y)=-\frac{1}{r}\sum_u  (\mathrm{ad}_{t_{u(1)}}\circ \mathrm{ad}_{t_{u(2)}}\circ \ldots \circ\mathrm{ad}_{t_{u(p-1)}} )(y),
$$
where the the index runs over all maps $u:\{1,\ldots,p-1\}\rightarrow \{x,y\}$ taking the value $x$ exactly $r$ times.
For $p=2$ the expression simplifies to $s_1(x,y)=[x,y]$, whereas $p=3$ gives $s_1(x,y)=[y,[x,t]]$ and $s_2(x,y)=[x,[x,y]]$.
Restricted $p$-Lie algebras were  introduced and studied by Jacobson \cite{Jac37},~\cite{Jac41},~\cite{Jac67}   and   also  go  under the name   \emph{$p$-Lie algebras}. More details can be found  to the monographs of Demazure and Gabriel \cite{DG70}, in particular Chapter II, Strade and Farnsteiner~\cite{FaSt88} or~\cite{Milne17}.

The prototype example of a restricted Lie algebra and the motivating example for its definition is the Lie algebra of $k$-derivations $\mathrm{Der}_k(R)$ of a $k$-algebra $R$, where $k$ is a field of characteristic $p>0$. In this case the Lie algebra structure is given by defining $[D_1,D_2]=D_1D_2-D_2D_1$, for any $k$-derivations $D_1,D_2$ of $R$ and the $p$-map is defined by $D^{[p]}=D^p$. Note that unlike in characteristic zero, in characteristic $p>0$ $D^p$ is a derivation. Indeed, let $x,y \in R$. Then by Leibniz's formula,
\[
D^p(xy)=xD^py+yD^px+\sum_{j=1}^{p-1}\binom{p}{j}D^jxD^{p-j}y.
\]
Now in characteristic $p>0$, $p$ divides $\binom{p}{j}$, for $1\leq j \leq p-1$ and therefore $D^P(xy)=xD^py+yD^px$ and hence $D^p$ is a $k$-derivation of $R$. Hence $\mathrm{Der}_k(R)$ is a restricted Lie algebra. The conditions on the $p$-map in the definition of a restricted $p$-Lie algebra are simply the properties that the $p$-map satisfies in the case of $\mathrm{Der}_k(R)$. 

In particular, when $R_n=k[x_1,\ldots,x_n]/(x_1^p,\ldots,x_n^p)$,  $\mathrm{Der}_k(R_n)$ is called the Jacobson-Witt algebra. This is a restricted simple Lie algebra. Restricted simple Lie algebras have been classified by R. Block and R. Wilson~\cite{BlWi88}.

Let $\mathfrak{g}$ be a restricted Lie algebra. A subalgebra is a subspace $\mathfrak{h}$ of $\mathfrak{g}$ such that  for any $x,y \in \mathfrak{h}$, $[x,y] \in \mathfrak{h}$ and for any $x\in \mathfrak{h}$, $x^{[p]}\in \mathfrak{h}$. An ideal of $\mathfrak{g}$ is a subspace $\mathfrak{I}$ of $\mathfrak{g}$ such that for any $a\in \mathfrak{g}$, $x\in\mathfrak{I}$, $[a,x]\in\mathfrak{I}$ and for any $x\in \mathfrak{I}$, $x^{[p]}\in \mathfrak{I}$. 

\subsection{Jacobson correspondence between purely inseparable extensions of height 1 and Lie algebras of derivations.}

Let $E$ be a field of characteristic $p>0$. Then $E^p=\{a^p,\; a\in E\}$ is a subfield of $E$. Then, as explained in the previous section, the $E$-vector space  $\mathfrak{g}=\mathrm{Der}(E)$ of derivations of $E$ is a Lie algebra over $E$ and a restricted $p$-Lie algebra over $E^p$. 

Let $F\subset E$ be a height 1 field extension. This means that for any $a\in E$, $a^p\in F$ and therefore $E^p\subset F\subset E$. Let $\mathfrak{g}(F)=\mathrm{Der}_F(E)$ be the set of $F$-derivations of $E$. This is a Lie subalgebra of $\mathfrak{g}$ which is also restricted over $F$.

Let $\mathcal{C}$ be the set of subfields $F$ of $E$ such that $E^p\subset F\subset E$. Let also $\mathcal{I}$ be the set of finite dimensional $E$-vector subspaces $\mathfrak{h}$ of $\mathfrak{g}=\mathrm{Der}(E)$ such that for any $D_1,D_2\in \mathfrak{h}$, $[D_1,D_2]\in \mathfrak{h}$ and for any $D\in \mathfrak{h}$, $D^p\in \mathfrak{h}$. In other words, $\mathfrak{h}$ is a restricted Lie subalgebra of $\mathfrak{g}$ over $E^p$ which is also an $E$-vector space (but not necessarily a restricted Lie algebra over $E$). For any $\mathfrak{h}\in \mathcal{I}$, let $I(\mathfrak{h})\subset E$ be the set of $x\in E$ such that $Dx=0$, for all $D\in \mathfrak{h}$. This is a subfield of $E$ such that $E^p\subset I(\mathfrak{h})\subset E$ and moreover, the extension $I(\mathfrak{h})\subset E$ is of height 1.

The next theorem gives a one to one correspondence between subfields $F$ of $E$ such that the extension $F \subset E$ is of height 1 and restricted Lie subalgebras of derivations of $E$.

\begin{theorem}[Jacobson Correspondence,~\cite{Jac41},~\cite{Jac44},~\cite{Jac67},~\cite{Bou81}]\label{Jacobson Correspondence}
Let $E$ be a field of characteristic $p>0$. The mappings $F\mapsto \mathfrak{g}(F)$ and $\mathfrak{h}\mapsto I(\mathfrak{h})$ are bijections between the set $\mathcal{C}$ of subfields $F\subset E$ such that $E^p\subset F\subset E$ and the set $\mathcal{I}$ of restricted Lie subalgebras of 
$\mathfrak{g}=\mathrm{Der}(E)$ which are also $E$-vector spaces. These mappings are mutually inverse.

Moreover, let $F\in \mathcal{C}$ and $\mathfrak{h}\in \mathcal{I}$ such that $F=I(\mathfrak{h})$. Then $[E:F]=p^d$, where $d=\dim_E \mathfrak{h}$.

\end{theorem}

The notions of purely inseparable extensions and the classification of height 1 extensions given by Jacobson in Theorem~\ref{Jacobson Correspondence} generalize naturally to algebraic schemes.

\subsection{Purely inseparable morphisms of schemes}
\begin{definition}
Let $\pi \colon X\rightarrow Y$ be a surjective finite morphism where $X$ and $Y$ are integral schemes over a perfect field $k$ of characteristic $p>0$. The map $\pi$ is called purely inseparable if and only if for any $x\in x$ and $y=\pi(x)\in Y$, the map $\mathcal{O}_{Y,y}\rightarrow \mathcal{O}_{X,x}$ is purely inseparable, i.e., for any $a \in  \mathcal{O}_{X,x}$, there exists $n>0$ such that $a^{p^n}\in \mathcal{O}_{Y,y}$. In particular the map between the function fields $k(Y)\rightarrow k(X)$ is purely inseparable. 

The height of $\pi$ is the smallest integer $n>0$ such that for any $x\in X$ and $a\in \mathcal{O}_{X,x}$, $a^{p^n}\in \mathcal{O}_{Y,y}$, where $y=\pi(x)$.

\end{definition}
Note that in the generality that the notion of a purely inseparable map $\pi \colon X\rightarrow Y$ is stated in the previous definition, it is possible that the height of $\pi$ is different than the height of the corresponding extension of function fields $k(Y)\subset k(X)$. However, if $Y$ is normal, then the two heights are the same.

\begin{proposition}
Let $\pi \colon X \rightarrow Y$ be a purely inseparable map of schemes such that $Y$ is normal. Then the height of $\pi$ is the same as the height of the extension $k(Y)\subset k(X)$.
\end{proposition}
\begin{proof}
This is a local result. So one may assume that $X=\mathrm{Spec}A$, $Y=\mathrm{Spec}B$ and $\pi$ is induced by an injective ring homomorphism $f \colon B \rightarrow A$ which induces a field extension $k(B)\subset k(A)$. Therefore there exists a commutative diagram
\[
\xymatrix{
B\ar[r]^f  \ar[d] &  A \ar[d] \\
k(B) \ar[r] & k(A)
}
\]
where $f$ is finite and $B$ integrally closed in $k(B)$. Let $n$ be the height of the extension $k(B)\subset k(A)$. Let $a\in A$. Then $a^{p^n}\in k(B)$. To prove the proposition it suffices to show that this element is actually an element of $B$. Since $a$ is integral over $B$, $a^{p^n}\in k(B)$ is also integral over $B$ and since $B$ is integrally closed it follows that $a^{p^n}\in B$. This concludes the proof.

\end{proof}

Purely inseparable maps appear naturally and are very important in the classification of surfaces, and higher dimension, in positive characteristic~\cite{BM77},~\cite{BM76}. Some of the simplest examples are the $p$-cyclic covers which are defined in the next two examples.

\begin{example}
Let $f\in k[x_1,x_2,\ldots,x_n]$ be a homogeneous polynomial of degree $mp$, where $k$ is a field of characteristic $p>0$ and $m>0$. Let $\mathbb{P}(m,1,\ldots, 1)$ be the $n$-dimensional weighted projective space with coordinates $y,x_i$, $i=1,\ldots,n$ and weights $m,1,\ldots, 1$. Let  $Z\subset \mathbb{P}(m,1,\ldots,1)$ be the hypersurface of $\mathbb{P}(m,1,\ldots,1)$ defined by the equation $y^p-f(x_0,\ldots,x_n)=0$. Then for general choice of $f$, $Z$ is normal and the projection $Z\rightarrow \mathbb{P}^{n-1}$ given by $[y,x_0,\ldots,x_n]\mapsto [x_0,\ldots,x_n]$ is a purely inseparable morphism of height 1.

\end{example}

\begin{example}
Let $X$ be a normal projective variety over a field $k$ of characteristic $p>0$.  Let $\mathcal{L}$ be an invertible sheaf on $X$ and $s \in H^0(X,\mathcal{L}^p)$ be a nonzero section. Then the section $s$ defines a map $\mathcal{L}^{-p}\rightarrow \mathcal{O}_X$ which gives the sheaf of $\mathcal{O}_X$-modules $R(\mathcal{L})=\oplus_{i=1}^{p-1}\mathcal{L}^{-i}$ a sheaf of $\mathcal{O}_X$-algebras structure. If in addition, $s$ is not a $p$-power, which means $s\not=t^p$, for some $t\in H^0(X,\mathcal{L})$, then $R(\mathcal{L})$ is a sheaf of integral domains and the map of schemes $\pi \colon Z =\mathrm{Spec}(R(\mathcal{L})) \rightarrow X$, is a purely inseparable map of degree $p$. Locally if $X=\mathrm{Spec} A$ and $s\in A$ corresponds to the section $s$, the map $\pi$ is simply the inclusion $A\subset A[x]/(x^p-s)$.

\end{example}

The Jacobson correspondence between purely inseparable height 1 extensions of fields and restricted Lie algebras of derivations presented in Theorem~\ref{Jacobson Correspondence} can be generalized in the case of purely inseparable morphisms of height 1 between normal varieties as follows.

First notice that if $\pi \colon X \rightarrow Y$ is purely inseparable of height 1, then there exists a factorization
\begin{equation}\label{diag1}
\xymatrix{ 
X \ar[rr]^F\ar[dr]_{\pi}  & & X^{(p)} \\
  &   Y \ar[ur]_{\phi} & \\
}
\end{equation}

Let now $X$ be a normal variety defined over a perfect field $k$ of characteristic $p>0$. Let $\mathcal{F} \subset \mathcal{T}_X$ be a coherent subsheaf of the tangent sheaf $\mathcal{T}_X=\mathcal{H}om_X(\Omega_{X/k},\mathcal{O}_X)$ of $X$.  $\mathcal{F}$ is called saturated if $\mathcal{T}_X/\mathcal{F}$ is torsion free and a subbundle if both $\mathcal{T}_X$ and $\mathcal{T}_X/\mathcal{F}$ are locally free. A saturated subsheaf $\mathcal{F}$ of $\mathcal{T}_X$ is always a subbundle on the codimension $\geq 2$ open subset $U$ of $X$ where $X$ is smooth.

A saturated coherent subsheaf $\mathcal{F}$ of $ \mathcal{T}_X$  is called a foliation on $X$ if and only if $\mathcal{F}$ is closed under Lie bracket and it is also  $p$-closed~\cite{Ekedahl87},~\cite{MiPe91}. These conditions mean that if $U \subset X$ is any open subset of $X$ and any section $D \in \mathcal{F}(U)$, then $D^p \in \mathcal{F}(U)$ and $[D,D^{\prime}]\in \mathcal{F}(U)$ for any $D^{\prime}\in \mathcal{F}(U)$, $D, D^{\prime}$ viewed as a $k$-derivations of $\mathcal{O}_U$.            
                                                                                                                                                             
Let $\mathcal{F} \subset T_X$ be a foliation on $X$. We define $Ann(\mathcal{F})\subset \mathcal{O}_X$ to be the sheaf of subrings of $\mathcal{O}_X$ on $X$ defined by 
\[
Ann(\mathcal{F}) (U)=\{a\in \mathcal{O}_X(U)|\; D(a)=0, \forall D\in \mathcal{F}(U)\},
\]
where $U \subset X$ is any open subset, and the sections of $\mathcal{F}$ over $U$ are viewed as $k$-derivations of $\mathcal{O}_U$. Let $S(\mathcal{F})=\mathrm{Spec} \left(Ann(\mathcal{F})\right)$. Then $\mathcal{O}_X^p \subset Ann(\mathcal{F})$ and hence there exists a factorization 
\begin{equation}\label{sec0-diag-1}
\xymatrix{ 
X \ar[rr]^F\ar[dr]_{\pi}  & & X^{(p)} \\
  &   S(\mathcal{F}) \ar[ur]_{\phi} & \\
}
\end{equation}
where $F \colon X \rightarrow X^{(p)}$ is the geometric Frobenius and $\pi$ is purely inseparable of height $1$.

Conversely, let $\pi \colon X \rightarrow Y$ be a purely inseparable map of height $1$. Then there exists a factorization as in diagram (\ref{diag1}).  Let $\mathcal{F}=Ann(\mathcal{O}_Y)\subset \mathcal{T}_X$ be defined by $\mathcal{F}(U)=\{D \in \mathcal{T}_X(U) |\; D(a)=0, \forall a\in \mathcal{O}_Y(U)\}$, for any open $U \subset X$. Then $\mathcal{F}$ is a foliation. The next theorem says that these two operations establish a one-to-one correspondence between foliations and normal varieties between $X$ and $X^{(p)}$.
  
\begin{theorem}[~\cite{Ekedahl87},~\cite{MiPe91}]\label{Jacobson schemes}  
Let $X$ be a normal variety defined over a perfect field $k$ of characteristic $p>0$. The mappings  $\mathcal{F}\mapsto S(\mathcal{F})$ and $Y\mapsto Ann(\mathcal{O}_Y)$ are bijections between the sets of  foliations of the tangent sheaf $T_{X}$ and purely inseparable maps of height 1 $\pi \colon X \rightarrow Y$. These mappings are mutually inverse.

Moreover, let $\mathcal{F} \subset \mathcal{T}_X$ be a foliation on $X$ of rank $r$ and let $Y =\mathrm{Spec}\left(Ann(\mathcal{F})\right)$. Then $[K(X):K(Y)]=p^r$
  \end{theorem}
 The next proposition gives some basic properties of the map $\pi \colon X \rightarrow S(\mathcal{F})$ and an adjunction formula for purely inseparable maps.
\begin{proposition}[\cite{Ekedahl87}~\cite{MiPe91}]\label{general-structure-theory}
Let $X$ be a normal variety defined over a perfect field of characteristic $p>0$. Let $\mathcal{F} \subset T_X$ be a foliation, $Y=S(\mathcal{F})$ and $\pi \colon X\rightarrow Y$ the corresponding map. Then. 
\begin{enumerate}
\item $Y$ is a normal variety.
\item Suppose that $X$ is smooth. Then $Y$ is smooth if and only if $\mathcal{F}$ is a subbundle.
\item Suppose that $X$ is smooth and $\mathcal{F}$ be a subbundle of $\mathcal{T}_X$. Then there exists a natural exact sequence
\begin{equation}\label{eq 1}
0 \rightarrow F^{\ast}(\mathcal{F}^{\ast}) \rightarrow \pi^{\ast}\Omega_{Y/k} \rightarrow \Omega_{X/k} \rightarrow \mathcal{F}^{\ast} \rightarrow 0,
\end{equation}
where $F \colon X \rightarrow X$ is the absolute Frobenius map. 
\end{enumerate}
\end{proposition}

\begin{corollary}\label{adjunction 2}
With assumptions as in the previous proposition,
\[
\omega_X=\left(\pi^{\ast}\omega_Y \otimes \left(\wedge^r\mathcal{F}\right)^{p-1}\right)^{\ast\ast},
\]
where $r$ is the rank of $\mathcal{F}$.
\end{corollary}    
\begin{proof}
Since $X$ is normal, to prove the formula it suffices to prove it on any open subset of $X$ of codimension $\geq 2$. Now because $X$ is normal, there exists an open subset $U$ of $X$ of codimension $\geq 2$ such that $U$ is smooth and the restriction of $\mathcal{F}$ on $U$ is a subbundle. Then Therefore it can be assumed that $X$, and hence $Y$ is smooth. and $\mathcal{F}$ locally free of rank $r$. Then the formula follows from the exact sequence (\ref{eq 1}).

\end{proof}              

More information about purely inseparable maps can be found in~\cite{DG70}, \cite{MiPe91}, \cite{Ekedahl87}, \cite{SP}, \cite{Tz17}.

\section{Group schemes in positive characteristic.}
Groups schemes  is a fundamental class of  algebraic schemes which, due to its importance,  has been heavily studied. Standard references in the subject are~\cite{DG70},~\cite{Milne17} and~\cite{Wat79}. The purpose of this section is to highlight some important differences in the structure of group schemes in positive characteristic in comparison to the characteristic zero case. In particular to present some results concerning the structure of  infinitesimal group schemes and their actions on algebraic varieties.

\subsection{Definition of group schemes and Hopf algebras}
\begin{definition}
Let $G$ be an algebraic scheme over a field $k$ and let $m\colon G\times G \rightarrow G$, $e\colon \mathrm{Spec} k \rightarrow G$, $i\colon G \rightarrow G$ be morphisms. $G$ is called a group scheme over $k$ if and only if the following diagrams commute.
\begin{enumerate}
\item \[
\xymatrix{
G\times G \times G \ar[rr]^{m\times 1_G} \ar[d]^{1_G \times m} & &G\times G \ar[d]^{m} \\
G\times G \ar[rr]^m && G
}
\]
\item 
\[
\xymatrix{
\mathrm{Spec} k \times G \ar[d]_{\cong} \ar[rr]^{e\times 1_G} & & G\times G \ar[d]^m & &G \times \mathrm{Spec} k \ar[ll]_{1_G \times e}\ar[d]^{\cong}\\
G \ar[rr]^{1_G} & &G && G\ar[ll]_{1_G}
}
\]
\item 
\[
\xymatrix{
           &.   G\times G \ar[dr]^m &           \\
G\ar[dr]_{(i,1_G)} \ar[ur]^{(1_G,i)} \ar[r] & \mathrm{Spec} k \ar[r]^e & G\\
           &     G\times G\ar[ur]_m
}
\]

\end{enumerate}
\end{definition}

An important and very useful equivalent definition is by using the functor of points of a scheme. 

\begin{definition}\label{definition group scheme}
An algebraic scheme $G$ over a field $k$ is a group scheme over  $k$ if and only if  for any scheme $S$ over $k$, the set of $k$-morphisms $\underline{G}(S)=\mathrm{Hom}(S, G)$ has a group structure and for any map of $k$-schemes $S\rightarrow T$, the induced map $\underline{G}(T) \rightarrow \underline{G}(S)$ is a group homomorphism. 
\end{definition}

Suppose that $G=\mathrm{Spec} A$ is an affine group scheme. Then the maps $m$, $e$ and $i$ in  Definition~\ref{definition group scheme} correspond to maps $m^{\ast}\colon A\rightarrow A\otimes_k A$, $e^{\ast} \colon A  \rightarrow k$ and $i^{\ast} \colon A\rightarrow A$. The conditions $(1)$, $(2)$, $(3)$ in Definition~\ref{definition group scheme} correspond to certain conditions on the maps $m^{\ast}$, $i^{\ast}$ and $e^{\ast}$. (the reader may deduce them directly from Definition~\ref{definition group scheme} or refer to~\cite[Page 7]{Pink05} or~\cite[Page 65]{Milne17}). This leads to the following notion of a Hopf algebra.

\begin{definition}[~\cite{Sweendler1967},~\cite{Milne17},~\cite{Pink05}~\cite{Abe77}]
A Hopf algebra over a field $k$ is a commutative  $k$-algebra $A$ together with $k$-algebra homomorphisms $m^{\ast} \colon A \rightarrow A\otimes A$, $e^{\ast} \colon A \rightarrow k$ and $i^{\ast} \colon A \rightarrow A$ such that the following diagrams are commutative.
\begin{enumerate}
\item 
\[
\xymatrix{
A\otimes_k A \ar[rr]^{m^{\ast}\otimes1_A} & &A\otimes_kA\otimes_k A \\
A\ar[u]^{m^{\ast}} \ar[rr]^{m^{\ast}} && A\otimes_k A\ar[u]_{1_A\otimes m^{\ast}}
}
\]
\item 
\[
\xymatrix{
k \otimes_k A && A\otimes_k A \ar[ll]_{e^{\ast}\otimes 1_A}  \ar[rr]^{1_A\otimes e^{\ast}} && A\otimes_k k\\
                   && A\ar[urr]_{\cong}. \ar[u]^{m^{\ast}} \ar[ull]^{\cong} && 
}
\]
\item
\[
\xymatrix{
A && A\otimes_k A \ar[ll]_{(i^{\ast},1_A)}\ar[rr]^{(1_A,i^{\ast})} && A \\
k\ar[u] && A\ar[ll]_{e^{\ast}} \ar[rr]^{e{\ast}} \ar[u]^{m^{\ast}} && A\ar[u]
}
\]
\end{enumerate}

\end{definition}

Comparing the definition of a group scheme to this of a Hopf algebra it immediately follows that a $k$-algebra $A$ is a Hopf algebra with maps $m^{\ast}$, $e^{\ast}$ and $i^{\ast}$ as above, if and only if $G=\mathrm{Spec} A$ is a group scheme over $k$ with multiplication $m$, identity element $e$ and inverse map $i$ induced by $m^{\ast}$, $e^{\ast}$ and $i^{\ast}$.

\subsection{Cartier duality}
Suppose that $G=\mathrm{Spec} A$ is a commutative group scheme. Let $\mu \colon A\otimes_k A\rightarrow A$  be the $k$-algebra homomorphism defined by $\mu(a_1\otimes a_2)=a_1a_2$ and $\varepsilon \colon k \rightarrow A$ be the $k$-algebra structure of $A$. Let also $m^{\ast}\colon A \rightarrow A\otimes A$, $i^{\ast} \colon A\rightarrow A$ and $e^{\ast} \colon A\rightarrow k$ be the $k$-algebra maps defining  the group scheme structure of $G$ as in the previous section. Let $A^{\ast}=\mathrm{Hom}_k(A,k)$ be the dual of $A$ as a $k$-vector space. Dualizing the maps $m^{\ast}$, $i^{\ast}$, $e^{\ast}$, $\mu$ and $\varepsilon$ one gets the following diagram.
\[
\xymatrix{
k \ar@/^/[rr]^{e^{\ast\ast}} && A^{\ast} \ar@/^/[ll]^{\varepsilon^{\ast}}  \ar@/^/[rr]^{\mu^{\ast}}  \ar@(dl,dr)[]_{i^{\ast\ast}}&& A^{\ast}\otimes_k A^{\ast}\ar@/^/[ll]^{m^{\ast\ast}}\\
                    &&                                                                         &&
}
\]
Then the $k$-vector space $A^{\ast}$ becomes a commutative (because $G$ is commutative) ring with multiplication defined by $m^{\ast\ast}$. Then one can easily check that the $k$-algebra, as defined previously, $A^{\ast}$ together with the maps $\mu^{\ast}$, $i^{\ast\ast}$ and $\varepsilon^{\ast}$ is a Hopf algebra. Therefore the  affine scheme $G^{\ast}=\mathrm{Spec} A^{\ast}$ is a group scheme over $k$. 

\begin{definition}The group scheme $G^{\ast}$ is called the Cartier dual of $G$. 
\end{definition}
For basic properties of Cartier duality the reader may refer to~\cite{Milne17} or~\cite{Pink05}.

\begin{example}\label{Cartier dual of mp}
Let $G=\mathbb{Z}/p\mathbb{Z}$ be the cyclic group of order $p>0$, where $p$ is any prime number. Let $k$ be a field of characteristic $p$ and $(\mathbb{Z}/p\mathbb{Z})_k$ the constant group scheme over $k$ defined in Example~\ref{ex constant gs}. Then, as shown in Examples~\ref{ex constant gs},~\ref{ex DG},
\[
(\mathbb{Z}/p\mathbb{Z})^{\ast}_k\cong D(\mathbb{Z}/p\mathbb{Z})\cong \mu_p.
\]
Similarly one may show~\cite{Pink05},~\cite{Milne17} thar $\alpha_p^{\ast}\cong \alpha_p$ and $\mu_p^{\ast}\cong (\mathbb{Z}/p\mathbb{Z})_k$.

\end{example}

\subsection{Group schemes are smooth in characteristic zero}
Let $G$ be a group scheme over a field $k$ and let $e\in G$ be the closed point corresponding to the identity map $e\colon \mathrm{Spec}k \rightarrow G$. Then it is well known~\cite[Proposition 1.8]{Milne17} that the following are equivalent. $G$ is smooth, $\mathcal{O}_{G,e}$ is a regular local ring and $G$ is geometrically reduced.

In characteristic zero all groups schemes are smooth as the following theorem of Cartier states.
\begin{theorem}[\cite{Cartier62},~\cite{Milne17},~\cite{Mu70}]
Any group scheme over a field $k$ of characteristic zero is smooth.
\end{theorem} 

However, this theorem fails in positive characteristic as the following examples show.

\subsection{Examples} Next we give some explicit examples of group schemes with emphasis on peculiarities of positive characteristic.

\begin{example}[The additive groups $\mathbb{G}_a$ and $\alpha_{p^n}$.]

Let $k$ be any field. The additive group scheme $\mathbb{G}_a$ over $k$ is defined as follows. As a scheme, $\mathbb{G}_a=\mathrm{Spec} k[x]$. The group scheme structure is given by the following $k$-algebra maps. Multiplication is defined by the map  $m^{\ast}\colon k[x]  \rightarrow k[x] \otimes k[x]$, where $m^{\ast}(x)=x\otimes 1+1\otimes x$. The identity element $e$ corresponds to the residue map $k[x]\rightarrow k$ given by $f(x)\mapsto f(0)$ and the inverse map $i$ corresponds to the map $i^{\ast}\colon k[x] \rightarrow k[x]$ defined by $i^{\ast}(x)=-x$.

As a functor of group valued points it is defined as follows. For any scheme $S$ over $k$, $\underline{\mathbb{G}}_a(S)=H^0(S,\mathcal{O}_S)$, as an abelian group.

Suppose that the field $k$ has characteristic $p>0$. Then the map $m^{\ast}$ defined above factors through $k[x]/(x^{p^n})$, for any $n>0$. Therefore the subscheme $\alpha_{p^n}$ of $\mathbb{G}_a$, where  $\alpha_{p^n}=\mathrm{Spec} k[x]/(x^{p^n})$ is a group scheme over $k$. It is a nonreduced  group scheme of length $p^n$ and is called the infinitesimal additive group scheme of order $p^n$.

As a functor of group valued points, $\alpha_{p^n}$   is defined as follows. Let $S$ be any scheme over $k$. Then
\[
\underline{\alpha}_{p^n}(S)=\{a\in H^0(S,\mathcal{O}_S),\; a^{p^n}=0\}.
\]
Since $k$ has characteristic $p>0$, this is an abelian group with addition.

\end{example}

\begin{example}[The multiplicative groups $\mathbb{G}_m$ and $\mu_{p^n}$.]   

Let $k$ be any field. The multiplicative group scheme $\mathbb{G}_m$ over $k$ is defined as follows. As a scheme, $\mathbb{G}_m=\mathrm{Spec} k[x,1/x]$. The group scheme structure is given by the following $k$-algebra maps. Multiplication is defined by the map  $m^{\ast}\colon k[x,1/x] \rightarrow  k[,1/x] \otimes k[x,1/x]$, where $m^{\ast}(x)=x\otimes x$. The identity element $e$ corresponds to the residue map $k[x,1//x]\rightarrow k$ given by $f(x)\mapsto f(1)$ and the inverse map $i$ corresponds to the map $i^{\ast}\colon k[x,1/x] \rightarrow k[x,1/x]$ defined by $i^{\ast}(x)=1/x$.

As a functor of group valued points it is defined as follows. For any scheme $S$ over $k$, $\underline{\mathbb{G}}_m(S)=H^0(S,\mathcal{O}_S)^{\ast}$, the group of units of $H^0(S,\mathcal{O}_S)$ as a multiplicative group.

For any $n>0$, the map $m^{\ast}$ factors through $k[x,1/x]/(x^n-1)\cong k[t]/(t^n-1)$ and therefore $\mu_n=\mathrm{Spec}k[x]/(x^n-1)$ becomes a group scheme. This is called the multiplicative group scheme of $n$-units. 

Suppose that $k$ has characteristic $p>0$ and $n=p^m$. Then $k[x]/(x^n-1)=k[x]/((x-1)^{p^m})$ and hence it is nonreduced. It is called the infinitesimal multiplicative group scheme of order $p^m$.

As a functor of group valued points, $\mu_{p^n}$   is defined as follows. Let $S$ be any scheme over $k$. Then
\[
\underline{\mu}_{p^n}(S)=\{a\in H^0(S,\mathcal{O}_S),\; a^{p^n}=1\}.
\]
Since $k$ has characteristic $p>0$, this is an abelian group with multiplication.
\end{example}

\begin{example}[The group scheme $D(G)$]\label{ex DG}
Let $G$ be a finitely generated group written multiplicatively. Let $k[G]$ be the group ring $k$-algebra. Its structure is the following. It is a vector space over $k$ with basis the elements $g\in G$. Multiplication in $k[G]$ is defined by
\[
\left(\sum_{i}a_i g_i\right) \cdot \left( \sum_j a^{\prime}_j g^{\prime}_j\right) =\sum_{i,j}(a_ia^{\prime}_j)(g_ig^{\prime}_j),
\]
Where $a_i, a^{\prime}_j \in k$ and $g_i, g^{\prime}_j\in G$.

Let $m^{\ast} \colon k[G] \rightarrow k[G]\otimes_k k[G]$, $i^{\ast} \colon k[G] \rightarrow k[G]$ and $e^{\ast} \colon k[G] \rightarrow k$ be the $k$-algebra homomorphisms defined by setting $m^{\ast}(g)=g\otimes g$, $i^{\ast}(g)=g^{-1}$ and $e^{\ast} (g)=1$, for all $g\in G$.
Then it is easy to verify that $k[G]$ with the maps $m^{\ast}$, $i^{\ast}$, $e^{\ast}$ is a Hopf algebra. The corresponding group scheme $G=\mathrm{Spec} k[G]$ is called the group scheme associated to $G$ and is denoted by $D(G)$. 

As a functor of group valued points it is defined as follows. Let $S$ be any scheme over $k$. Then 
\[
\underline{D}(G)(S)=\mathrm{Hom}_{group}(G,H^0(S,\mathcal{O}_S)^{\ast}).
\]

Suppose for example that $G$ is cyclic of order $n$. Let $g$ be a generator of $G$. Then the map $\phi\colon k[x]/(x^n-1) \rightarrow k[G]$ defined by $f(x)=g$ is an isomorphism of $k$ algebras. Then $m^{\ast}\colon k[G] \rightarrow k[G]\otimes k[G]$ defined above gives a $k$-algebra map 
\[
\Delta \colon \frac{k[x]}{(x^n-1)} \rightarrow \frac{k[x]}{(x^n-1)} \otimes \frac{k[x]}{(x^n-1)}
\]
Given by $\Delta(x)=x\otimes x$. Therefore $k[G]\cong \mu_n$. 

Let $p$ be the characteristic of $k$. If $p$ does not divide $n$, then $\mu_n$, and hence $D(G)$, is a smooth group scheme and the abstract group of its closed $k$-points is isomorphic to $\mathbb{Z}/n\mathbb{Z}$, the cyclic  closed points of order $n$.
If on the other hand $p$ divides $n$ then $\mu_n$ is not smooth. In particular, if $n=p^m$, then $D(G)\cong \mu_{p^m}$, which is an infinitesimal non reduced group scheme. 

Let $D(G)^{\ast}$ be the Cartier dual of $D(G)$. Then it is not difficult to check that this is a smooth group scheme and its group of $k$-points is isomorphic to $G$.

\end{example}

\begin{example}[The constant group scheme]\label{ex constant gs}
Let $G$ be a finite  group and $k$ any field. For any $g \in G$,  let $S_g=\mathrm{Spec} k$. Let also $G_k=\coprod_{g\in G} S_g$. As a scheme this is simply the disjoint union of $|G|$ closed points. Let $s_g \in G_k$ be the closed point of $S_g$. Then $G_k$ becomes a group scheme by defining multioplication $m \colon G_k \times G_k \rightarrow G_k$, inverse $i\colon G_k \rightarrow G_k$ and identity $e \colon \mathrm{Spec} k \rightarrow G_k$ as follows. For any $(s_g,s_{g^{\prime}}) \in G_k \times G_k$,  $m(s_g,s_{g^{\prime}})=s_{gg^{\prime}}$. For any $s_g \in G_k$, $i(s_g)=s_{g^{-1}}$ and $e$ is the inclusion $S_1 \rightarrow G_k$, where $1$ is the unit element of $G$.

Next we will describe the ring of regular functions $A$ of $G$ and the Hopf algebra structure of $A$ inducing the group scheme structure. This is particularly useful in order to relate $G_k$ with the group scheme $D(G)$ constructed in Example~\ref{ex DG}.

Let  $A=\{f\colon G \rightarrow k,\,\; \text{$f$ is a map of sets}\}$. Define addition and multiplication $A$ be setting $(f+h)(g)=f(g)+h(g)$ and $(f\cdot h)(g)=f(g)h(g)$. Then $A$ is a commutative $k$-algebra and it is simply isomorphic to the product $\times_{i=1}^{|G|} k$. Therefore, as a scheme $G_k =\mathrm{Spec} A$. Let $g\in G$. The maps $f_g \colon G \rightarrow k$ defined by setting $f_g(g^{\prime})$ to be $1$ if $g=g^{\prime}$ and $0$ if $g\not= g^{\prime}$ form a $k$-basis of $A$. The group scheme structure of $G_k$, as defined above, correspond to maps of $k$-algebras $m^{\ast}\colon A \rightarrow A\otimes_k A$, $i^{\ast} \colon A\rightarrow A$ and $e^{\ast} \colon A \rightarrow k$ defined as follows.
\begin{enumerate}
\item For any $g\in G$,
\[
m^{\ast} (f_g)=\sum_{g^{\prime}g^{\prime\prime}=g}f_{g^{\prime}}\otimes f_{g^{\prime\prime}}.
\]
\item For any $g\in G$, $i^{\ast}(f_g)=f_{g^{-1}}$.
\item For any $f\in A$, $e^{\ast}(f)=f(1)$.
\end{enumerate}
The group scheme $G_k$ is called the constant group scheme corresponding to $G$. Its set of closed points is a group isomorphic to $G$.

Suppose now that $G$ is finite and commutative.  Let $G_k^{\ast}$ be the Cartier dual of $G_k$. We will next show that $G_k^{\ast}$ is isomorphic as a group scheme to the group scheme $D(G)$ defined in the previous example.

By its definition, $G_k^{\ast} =\mathrm{Spec} A^{\ast}$, where $A^{\ast}=\mathrm{Hom}_k(A,k)$, is a commutative $k$-algebra with multiplication the dual $m^{\ast\ast}$ of $m^{\ast}$. The group scheme structure of $G_k^{\ast}$ is given by the dual map $\mu^{\ast}$, $i^{\ast}$ and $\varepsilon^{\ast}$, where $\mu \colon A\otimes A \rightarrow A$ is the multiplication map of $A$.

Let us now look more carefully to the multiplication in $A^{\ast}$. A $k$-basis for $A^{\ast}$ is the dual basis $f^{\ast}_g$ of the $k$-basis $f_g$ of $A$, $g\in G$, as defined above. Then 
$f_g^{\ast} \cdot f_{g^{\prime}}^{\ast} =m^{\ast\ast}(f_g^{\ast} \otimes f_{g^{\prime}}^{\ast})$. Let $g^{\prime\prime} \in G$. By checking the definitions of all the maps involved, it easily follows that $m^{\ast\ast}(f_g^{\ast} \otimes f_{g^{\prime}}^{\ast})(f^{\ast}_{g^{\prime\prime}})$ is $1$ if $g^{\prime\prime}=gg^{\prime}$ and zero otherwise. Therefore $f_g^{\ast} \cdot f_{g^{\prime}}^{\ast} =f^{\ast}_{gg^{\prime}}$. Similarly one can easily see that $i^{\ast\ast}(f^{\ast}_g)=f^{\ast}_{g^{-1}}$.

Let $k[G]$ be the group $k$-algebra of $G$. Define now the map 
\[
\Phi \colon A^{\ast} \rightarrow k[G]
\]
 as follows. Let $\phi=\sum_{g\in G}\lambda_gf_g^{\ast}\in A^{\ast}$, where $\lambda_g \in k$.  Then $\Phi(\phi)=\sum_{g\in G} \lambda_g g \in k[G]$. The description of multiplication in $A^{\ast}$ shows that this is an isomorphism of group algebras. One may also check that this is in fact a group scheme map. Therefore 
 \[
 G_k^{\ast}\cong D(G),
 \]
as group schemes.

\end{example}

\begin{example}[Group schemes of automorphisms]
Let $X$ be a scheme of finite type defined over a field $k$. The automorphism functor $\underline{\mathrm{Aut}}_{X/k}$ of $X$ is the  functor of groups which is defined as follows. For any $k$-scheme $S$, 
\[
\underline{\mathrm{Aut}}_k(X)(S)=\mathrm{Aut}_S(X\times S),
\]
the group of $S$-automorphisms of $X\times S$. If $X$ is proper over $k$ then the group valued functor $\underline{\mathrm{Aut}}_{X/k}$ is represented by a scheme $\mathrm{Aut}_{X/k}$ locally of finite type over $k$~\cite{Kollar96}. This is called the group scheme of automorphisms of $X$ and is a fundamental object of study in algebraic geometry.

Let $p$ be the characteristic of $k$. If $p=0$ then $\mathrm{Aut}_k(X)$ is smooth since group schemes in characteristic zero are always smooth. However, if $p>0$, then it may not be smooth~\cite{Tz17},~\cite{Tz17b},~\cite{Tz22},~\cite{ScTz23},~\cite{Martin22}.

A very simple example which exhibits this situation is the following. Let  $X=\mathrm{Spec} A$, where $A=k[x]/(x^2)$, where $k$ is a field of characteristic 2. We will show that 
\[
\mathrm{Aut}_k(X)\cong \alpha_2 \rtimes \mathbb{G}_m. 
\]
In particular it is not smooth. 

Let $S=\mathrm{Spec} R$ be an affine $k$-scheme. Then
\[
\mathrm{Hom}(S,\mathrm{Aut}_k(X))=\mathrm{Aut}_R \left(R[x]/(x^2)\right).
\]
Let now $\phi \colon R[x]/(x^2)\rightarrow R[x]/(x^2)$ be an $R$-isomorphism. This is determined by $\phi(x)$. Suppose that $\phi(x)=a+bx$, $a,b\in R$. Then 
\[
0=\phi(x^2)=(\phi(x))^2=a^2+b^2x^2=a^2.
\]
Hence $a^2=0$. Moreover, since $\phi$ is an $R$-isomorphism, $b\in R^{\ast}$.  Define now 
\begin{equation}\label{sec1-eq3}
\Psi \colon  \mathrm{Aut}_R \left(R[x]/(x^2)\right) \rightarrow   R^{[2]} \times R^{\ast}
\end{equation}
be setting $\Psi(\phi)=(a,b)$, where $\phi $ is an $R$-automorphism of $R[x]/(x^2)$, $\phi(x)=a+bx$ and $R^{[2]}=\{r\in R,\;\; r^2=0\}$. This is a bijection of sets.

Let now $\phi_1$, $\phi_2$ be $R$ automorphisms of $R[x]/(x^2)$. Let $\phi_i(x)=a_i+b_ix$, $i=1,2$, with $a_i^2=0$ and $b_i\in R^{\ast}$, $i=1,2$. Then $\phi_1\circ \phi_2(x)=a_1+b_1a_2+b_1b_2x$. Define now multiplication in $R^{[2]}\times R^{\ast}$ be setting 
\[
(a_1,b_1)\cdot (a_2,b_2)=(a_1+b_1a_2,b_1,b_2).
\]
Then $\Psi$ becomes a group isomorphism. Moreover, this multiplication exhibits $R^{[2]}\times R^{\ast}$ as the semidirect product $R^{[2]}\rtimes R^{\ast}$ via the map $f \colon R^{\ast}\rightarrow \mathrm{Aut}(R^{[2]})$, where for any $r\in R^{\ast}$, $f(r)$ is the automorphism of $R^{[2]}$ given by $f(r)(a)=ra$, $a\in R^{[2]}$. From this it immediately follows that $\mathrm{Aut}_k(X)\cong \alpha_2 \rtimes \mathbb{G}_m$ as claimed.

\end{example}

%\begin{example}[The Tate Oort group scheme.]

%\end{example}

The previous examples show that in positive characteristic a group scheme $G$ may not be smooth. In this case one may ask if the reduced subscheme $G_{red}$ is a group subscheme of $G$ and therefore smooth since smoothness and reducedness of group schemes are equivalent notions~\cite[Proposition 1.28]{Milne17}. This is true if $k$ perfect~\cite[Proposition 1.38]{Milne17} but not always if $k$ is not perfect, as the following example shows.

\begin{example}[$G_{red}$ may not be a group scheme, Exercise 1.57~\cite{Milne17}]
Let $k$ be a non perfect field of characteristic $p>0$ and let $\mathbb{G}_a^2$ be the product group scheme  $\mathbb{G}_a\times \mathbb{G}_a$. As a scheme it is $\mathrm{Spec}k[x,y]$ with group scheme structure defined by the map $m^{\ast} \colon k[x,y] \rightarrow k[x,y]\otimes_k k[x,y]$ such that  $m^{\ast}(x)= x\otimes 1+1\otimes x$ and $m^{\ast}(y)=y\otimes 1+1\otimes y$.  Let $G$ be the subscheme of $\mathbb{G}_a^2$ given by the equation $x^{p^2}-tx^p=0$, where $t\in k-k^p$. Then the group structure map $m^{\ast}$ factors through $G$ and therefore $G$ is a sub group scheme of $\mathbb{G}_a\times \mathbb{G}_a$.

The reduced subscheme $G_{red}$ of $G$ is defined by the equation $x(x^{p^2-p}-t)=0$. We will show that not only $G_{red}$ is not a subgroup scheme of $G$, but the scheme $G_{red}$ does not admit any group scheme structure. Suppose that this was true, i.e., $G$ admits some group scheme structure. It is easy to see that $G_{red}$ is a regular scheme. Since $G_{red}$ admits some group scheme structure,  it must be smooth~\cite[Proposition 1.28]{Milne17} (note that over non perfect fields regularity and smoothness are not equivalent notions). In particular it is geometrically reduced. Let $\bar{k}$ be the algebraic closure of $k$. Then in $\bar{k}[x,y]$, $x(x^{p^2-p}-t)=x(x^{p-1}-s)^p$, where $s^p=t$.  Therefore $G\otimes_k \bar{k}$ is not reduced, a contradiction.

\end{example}
The last example shows another property unique of group schemes in positive characteristic which indicates how much richer the theory of group schemes is in this case.

\begin{example}[Group schemes of order $p^2$ may not be commutative]\label{non commutative of order 2p}
In characteristic zero many properties from the theory of finite abstract groups carry to the theory of group schemes. For example, let $G=\mathrm{Spec}A$ be a finite group scheme over a field $k$ of characteristic $p\geq 0$. Then $A$ is a finite dimensional $k$-algebra. The dimension of $A$ as a $k$-vector space is called the order of $G$ and is denoted by $|G|$. Suppose that $p=0$. Then $G$ is smooth and hence just a disjoint union of $|G|$ reduced points. Then the set of closed points of $G$ is an abstract group many of whose properties carry to $G$. In particular, if $|G|=p^2$, where $p$ is a prime number, $G$ is commutative.  However, this is not true if $p>0$ and $|G|=p^2$.

Let for example $G=\mathrm{Spec} A$, where $A=k[s,t]/(s^p,t^p-1)$, where $k$ is a field of characteristic $p>0$. Then $G$ becomes a group scheme by defining $m^{\ast}\colon A\rightarrow A\otimes_k A$, $i^{\ast}\colon A\rightarrow A$ and $e^{\ast} \colon A\rightarrow k$ by setting $m^{\ast}(t)=t\otimes  t$, $m^{\ast}(s)=t\otimes s +s\otimes 1$, $i^{\ast}(t)=t^{p-1}$, $i^{\ast}(s)=-s$ and $e^{\ast}(f(s,t))=f(0,1)$, for any $f\in k[s,t]$. $G$ is not commutative of order $p^2$. 

The group scheme $G$ constructed above is isomorphic to the semidirect product $\alpha_p \rtimes \mu_p$ where the semidirect product structure is given as follows. Let $R$ be any commutative finitely generated  $k$-algebra. Let $\phi_R \colon \mu_p(R) \rightarrow \mathrm{Aut}(\alpha_p(R))$ be the group homomorphism defined as follows. Let $r\in \mu_p(R)=\{r\in R,\;\ \; r^p=1\}$ and $r^{\prime}\in \alpha_p(R)=\{r\in R,\;\;r^p=0\}$. Then $\phi(r)(r^{\prime})=rr^{\prime}$. Then the group structure on $\alpha_p(R)\rtimes \mu_p(R)$ is defined by setting
\[
(r_0,r_1)\cdot (r_0^{\prime},r_1^{\prime})=(r_0+r_0^{\prime}r_1,r_1r_1^{\prime}).
\]
The definition is functorial in $R$ and the group functor defined this way is represented by $G$.

\end{example}

\subsection{The Lie algebra of a group scheme}

Let $X$ be any algebraic scheme defined over a field $k$. A vector field on $X$ is a map of sheaves $D\colon \mathcal{O}_X \rightarrow \mathcal{O}_X$ such that for any open set $U\subset X$, $D(U)\colon \mathcal{O}_X(U)\rightarrow \mathcal{O}_X(U)$ is a $k$-derivation. A rational vector field on an integral algebraic scheme $X$ is a vector field of some non empty dense open subset $U$ of $X$. 
 
Vector fields of $X$ are in one to one correspondence with elements of $\mathrm{Hom}_X(\Omega_{X/k},\mathcal{O}_X)$. In the case when $X$ is smooth this is simply $H^0(X,T_X)$.

Let $X$ and $Y$ be algebraic schemes over a field $k$, let $X\times Y$ be their product and $p\colon X\times Y\rightarrow X$, $q\colon X\times Y \rightarrow Y$ be the projections. Then there exist a natural isomorphism 
$\Omega_{X\times Y/k} \rightarrow  p^{\ast} \Omega_{X/k} \oplus q^{\ast}\Omega_{Y/k}$. Let $D$ be a vector field on $Y$ which corresponds to a map $\Omega_{Y/k}\rightarrow \mathcal{O}_Y$. Then from the previous isomorphism it follows that $D$ induces a vector field denoted by $1\otimes D$ on $X\times Y$. Locally it is defined as follows. Suppose that $X=\mathrm{Spec} A$ and $Y=\mathrm{Spec}B$. Then $D$ is a $k$-derivation of $B$. Then $1\otimes D $ is the $k$-derivation of $A \otimes B$ defined by $(1\otimes D)(a\otimes b)=a\otimes Db$.

Next we define the notion of left invariant vector field on a group scheme $G$.

\begin{definition}\label{left invariant}
Let $G$ be a group scheme over a field $k$ and $D\colon \mathcal{O}_X \rightarrow \mathcal{O}_X$ a vector field on $G$. The vector field $D$ is called left invariant if and only if the following diagram is commutative.
\begin{equation}\label{left invariant diagram}
\xymatrix{
\mathcal{O}_G\ar[rr]^D \ar[d]^{m^{\ast}} & & \mathcal{O}_G \ar[d]^{m^{\ast}} \\
m_{\ast}\mathcal{O}_{G\times G} \ar[rr]^{m_{\ast}(1\otimes D)} && m_{\ast} \mathcal{O}_{G\times G}
}
\end{equation}
Where $m^{\ast}$ is the map induced from the multiplication map $m$.
\end{definition}

Let now $e\colon \mathrm{Spec} k \rightarrow G$ be the unit map of a group scheme $G$ defined over a field $k$ and let $T_{G,e}$ be the tangent space of $G$ at the closed point $e$ of $G$. Then 
\[
T_{G,e}=\mathrm{Hom}_k(\mathbf{m}_{G,e}/\mathbf{m}^2_{G,e}, k(e))\cong \mathrm{Der}_k(\mathcal{O}_{G,e}, k),
\]
Where $\mathbf{m}_{G,e}$ is the maximal ideal of $\mathcal{O}_{G,e}$. Hence $T_{G,e}$ is canonically isomorphic to the $k$-vector space of $k$-derivations $\delta \colon \mathcal{O}_{G,e}\rightarrow k$. In the following we will use this interpretation of the tangent space.

Let $D\colon \mathcal{O}_G \rightarrow \mathcal{O}_G$ be a vector field of $G$. Let $e^{\ast}\colon \mathcal{O}_{G,e}\rightarrow k$ be the $k$-algebra map inducing $e$. By composing $D$ with $e^{\ast}$ we get a tangent vector $\delta$ of $G$ at $e\in G$. 

The next proposition shows that  $\delta$ completely determines $D$ and that the tangent space of $G$ at $e$ is isomorphic to the space of left invariant vector fields of $G$.
\begin{proposition}\label{left invariant and tangent vectors}
Let $\delta\in T_{G,e}$ be a tangent vector of $G$ at $e$. Then there exists a unique left invariant vector field on $G$ which induces $\delta$ at $e$.
\end{proposition}
\begin{proof}
The proof will be in two steps. In the first step it will be shown that a left invariant vector field of a group scheme $G$ is completely determined by the tangent vector of $G$ at the identity point $e$ that it induces. In the second part of the proof, given a tangent vector $\delta$ of $G$ at $e$, a left invariant vector field $D$ on $G$ inducing $\delta$ will be explicitly constructed.

\textbf{Step 1.} Let $D \colon \mathcal{O}_G \rightarrow \mathcal{O}_G$ be a left invariant vector field of $G$ and let $\delta \colon \mathcal{O}_{G,e}\rightarrow k$ be the tangent vector (a $k$-derivation) of $G$ at $e$ that it induces.  In order to show that $D$ is uniquely determined by $\delta$, it suffices to show that $D=0$ if and only if $\delta =0$.  From the diagram (\ref{left invariant diagram}) in the definition of left invariant vector fields we get the commutative diagram
\begin{equation}\label{diagram 4}
\xymatrix{
\mathcal{O}_G \ar[rr]^D \ar[d]^{m^{\ast}} \ar@/^3pc/[drrr]^D& & \mathcal{O}_G \ar[d]^{m^{\ast}} \ar[dr]^{id}& \\
m_{\ast}\mathcal{O}_{G\times G} \ar[rr]^{m_{\ast}(1\otimes D)} \ar@/_2pc/[rrr]^{\sigma}& &m_{\ast} \mathcal{O}_{G\times G} \ar[r]^{\tau} & \mathcal{O}_G
&& \\
&&
}
\end{equation}
Where $\sigma=(1\otimes D)\circ \tau$ and $\tau$ is the sheaf of rings map induced from the map 
\[
G\cong G\times \mathrm{Spec}k \colon \overset{1_G \times e}{\rightarrow} G\times G.
\]
Let $g\in G$ be any closed $k$-point. Since $m(g,e)=g$, at the  point $g$ the map $\sigma$ is the composition
\begin{equation}\label{diagram 5}
\xymatrix{
\mathcal{O}_{G,g}\otimes \mathcal{O}_{G,e} \ar[r]^{1_{G,g} \otimes D_e} \ar[dr]^{\sigma}& \mathcal{O}_{G,g}\otimes \mathcal{O}_{G,e} \ar[d]^{\tau} \\
 & \mathcal{O}_{G,g}
 }
 \end{equation}
 Where $D_e$ is the $k$-derivation of $\mathcal{O}_{G,e}$ induced by $G$. Now from the definition of $\tau$ it follows that for any $a\in \mathcal{O}_{G,g}$ and $b\in \mathcal{O}_{G,e}$, $\tau(a\otimes b)=\bar{b}a$, where $\bar{b}\in k$ is the image of $b$ in $k=\mathcal{O}_{G,e}/\mathbf{m}_{G,e}$. Therefore $\sigma (a\otimes b)=\delta(b)a$. Hence if $\delta=0$ at $e$, then $\sigma =0$ at any closed $k$-point $g\in G$. From this, since $D=\sigma m^{\ast}$, it follows that $D$ is zero at any closed point $g\in G$ and hence $D=0$. 
 
 \textbf{Step 2.} Let $\delta \colon \mathcal{O}_{G,e}\rightarrow k$ be a tangent vector of $G$ at $e$. We will show that there exists a left invariant vector field $D$ on $G$ restricting to $\delta$ at $e$.

 The tangent vector $\delta$ of $G$ at $e$ defines a morphism $\phi \colon \mathrm{Spec} k[t]/(t^2) \rightarrow G$~\cite{Ha77}. Locally at an affine neighborhood $U=\mathrm{Spec} A$ of $e$ in $G$ the map $\phi$ is induced from the $k$-algebra map $\phi^{\ast} \colon A\rightarrow k[t]/(t^2)$ defined by $\phi^{\ast}(a)=\bar{a}+\delta(a)t$, where $\bar{a}$ is the image of $a$ in $A/\mathbf{m}_{G,e}=k$. 

The group scheme structure on $G$ defines the $G$-isomorphism
\[
\xymatrix{
G\times G\ar[rr]^{m \times p_2  }\ar[dr]_{p_2} & & G\times G \ar[dl]^{p_2}\\
 & G &
}
\]
Where $p_2$ is the second projection. Pulling back by $\phi$ we get the $ \mathrm{Spec} k[t]/(t^2) $ isomorphism $\Phi$
\[
\xymatrix{
G\times \mathrm{Spec} \frac{k[t]}{(t^2)}  \ar[rr]^{ \Phi} \ar[dr]_{p_2} & & G\times \mathrm{Spec} \frac{k[t]}{(t^2)} \ar[dl]^{p_2}\\
 &  \mathrm{Spec} \frac{k[t]}{(t^2)}  &
}
\]
which pulls back to the identity automorphism of $G$ by the residue map $k[t]/(t^2)\rightarrow k$. Therefore $\Phi$ is the identity map on the underlying topological spaces. Let $U=\mathrm{Spec} A$ be an affine open subset of $G$. Then on $U$, $\Phi$ is defined by a $k[t]/(t^2)$-algebra homomorphism  $\Phi^{\ast} \colon A[t]/(t^2) \rightarrow A[t]/(t^2)$ which becomes the identity modulo $(t)$. Hence it easily follows that $\Phi(a+bt)=a+(b+f(a))t$, for some map $f\colon A\rightarrow A$. The fact that $\Phi$ is a ring homomorphism shows that $f$ is a $k$-derivation of $A$. These local derivations patch to  a vector field $D$ on $G$. A careful look at the construction of $D$ shows that $D$ is left invariant and induces $\delta$ at $e$. The details are left to the reader.

Finally, note that if $G=\mathrm{Spec} A$ is an affine group scheme, then $D$ can be defined directly from the diagrams (\ref{diagram 4}) and (\ref{diagram 5}) by setting $D=\sigma m^{\ast}$, where $\sigma \colon A\otimes A\rightarrow k$ was defined by setting $\sigma(a\otimes b)=\delta b a$. One easily checks that this is a left invariant derivation of $A$.
\end{proof}

Next we define the notion of Lie algebra of a group scheme $G$ over a field $k$. It is a fundamental invariant of $G$ and it determines to some extent (completely in some cases that we discuss see in the next section

\begin{definition}[The Lie algebra of a group scheme~\cite{DG70},~\cite{Mu70},~\cite{Milne17}]
Let $G$ be a group scheme over a field $k$. The $k$-vector space of left invariant vector fields of $G$ with bracket defined by $[D_1,D_2]=D_1D_2-D_2D_1$, where $D_1,D_2$ are left invariant vector fields on $G$, is a Lie algebra which is called the Lie algebra of $G$ and is denoted by $\mathfrak{Lie}(G)$.

Suppose that $k$ has characteristic $p>0$. Then, as explained in Section~\ref{sec-1}, if $D$ is a vector field of $G$,  $D^p$ is also a vector field of $G$.  Then  $\mathfrak{Lie}(G)$ becomes a restricted Lie algebra with bracket as defined above and $p$-map $D^{[p]}=D^p$.
\end{definition}

Let $G$ be a group scheme over $k$. Define
\[
\mathrm{Hom}^e_k\left(\mathrm{Spec}\frac{k[t]}{(t^2)}, G\right)=\mathrm{Ker}\left(\underline{G}(k[t]/(t^2)\rightarrow \underline{G}(k)\right).
\]
This is the $k$-vector space of maps $\mathrm{Spec}k[t]/(t^2) \rightarrow G$ such that here is a factorization
\[
\xymatrix{
\mathrm{Spec}\frac{k[t]}{(t^2)}\ar[rr] & & G \\
& \mathrm{Spec}k \ar[ul] \ar[ur]^e& 
}
\]
Considering that such maps are in one to one correspondence with $k$-derivarions $\delta \colon \mathcal{O}_{G,e} \rightarrow k$, Proposition~\ref{left invariant and tangent vectors} implies the following.
\begin{corollary}\label{corollary 10}
There exists a natural isomorphism of $k$ vector spaces
\[
\mathfrak{Lie}(G) \cong \mathrm{Hom}_k(\mathbf{m}_{G,e}/\mathbf{m}^2_{G,e},k(e))\cong \mathrm{Hom}^e_k\left(\mathrm{Spec}\frac{k[t]}{(t^2)}, G\right),
\]
Where $k(e)=\mathcal{O}_{G,e}/\mathbf{m}_{G,e}=k$ and $\mathrm{Hom}^e_k\left(\mathrm{Spec}\frac{k[t]}{(t^2)}, G\right)$ is the $k$-vector space of maps 
\end{corollary}

The next proposition shows that the correspondence $G\mapsto \mathfrak{Lie}(G)$ defines a functor between the category of group schemes over $k$ and finite dimensional restricted Lie algebras over $k$.

\begin{proposition}\label{Lie to gs map}
Let $f\colon G_1\rightarrow G_2$ be a morphism of group schemes. Then $f$ induces a homomorphism of restricted Lie algebras  
$$
\mathfrak{Lie}(f) \colon \mathfrak{Lie}(G_1) \rightarrow \mathfrak{Lie}(G_2).
$$

The previous map is functorial in the following way. Let $f_1\colon G_1\rightarrow G_2$ and $f_2\colon G_2\rightarrow G_3$ are group scheme morphisms, then 
\[
\mathfrak{Lie}(f_2\circ f_1)=\mathfrak{Lie}(f_2)\circ \mathfrak{Lie}(f_1).
\]
\end{proposition}
\begin{proof}
From Corollary~\ref{corollary 10}, there are natural isomorphisms between $\mathfrak{Lie}(G_i)$ and the tangent space of $G_i$ at the identity point, $i=1,2$. Then it is clear that there exists a map $\mathfrak{Lie}(f)\colon \mathfrak{Lie}(G_1)\rightarrow \mathfrak{Lie}(G_2)$ which is  a homomorphism of $k$-vector spaces. We will show that this map is in fact a homomorphism or restricted Lie algebras, i.e., it preserves the bracket and the $p$-map. In order to show this it is necessary to explicitly describe the map $\mathfrak{Lie}(f)$. Let $D_1$ be a left invariant vector field of $G_1$ and let $\delta_1\colon \mathcal{O}_{G_1}\rightarrow k$ be the induced tangent vector at the identity of $G_1$. Composing with the map $f^{\ast}\colon \mathcal{O}_{G_2}\rightarrow \mathcal{O}_{G_1}$ we get  a tangent vector $\delta_2 \colon \mathcal{O}_{G_2}\rightarrow k$ to the identity of $G_2$. Then, from the proof of Proposition~\ref{left invariant and tangent vectors}, the left invariant vector field $D_2=\mathfrak{Lie}(f)(D_1)$ of $G_2$ fits, and in fact is constructed from  the following commutative diagram.

\begin{equation}\label{Lie Map}
\xymatrix{
\mathcal{O}_{G_2}\ar[rr]^{m_2^{\ast}}\ar[d]^{f^{\ast}} \ar@/^2pc/[rrrr]^{D_2}  &&  \mathcal{O}_{G_2}\otimes \mathcal{O}_{G_2}\ar[d]^{f^{\ast}\otimes f^{\ast}}  \ar[rr]^{1_{G_2}\otimes \delta_2}\ar[d]^{f^{\ast}} &&\mathcal{O}_{G_2}\ar[d]^{f^{\ast}} \\
\mathcal{O}_{G_1}\ar[rr]^{m_1^{\ast}} \ar@/_1pc/[rrrr]_{D_1} && \mathcal{O}_{G_1}\otimes \mathcal{O}_{G_1} \ar[rr]^{1_{G_1}\otimes \delta_1}&& \mathcal{O}_{G_1}
}
\end{equation}
By using the commutativity of this diagram one can easily check that $\mathfrak{Lie}(f)$ preserves the bracket and $p$-map and therefore is a homomorphism of restricted Lie algebras.

he proof of the second part of the proposition is by explicitly describing all the maps $\mathfrak{Lie}(f_1)$, $\mathfrak{Lie}(f_2)$ and $\mathfrak{Lie}(f_2\circ f_1)$ by using the diagram (\ref{Lie Map}).
\end{proof}

Finally we apply the previous theory to the case of the automorphism group scheme of an algebraic scheme. More precisely, let $X$ be a proper algebraic scheme $X$ defined over a field $k$ and let $\underline{\mathrm{Aut}}_{X/k}$ be the group valued functor defined by $\underline{\mathrm{Aut}}_{X/k}(S)=\mathrm{Aut}_S(X\times S)$, where $S$ is any scheme over $k$. Then it is well known~\cite{Kollar96} that $\underline{\mathrm{Aut}}_{X/k}$ is represented by a group scheme $\mathrm{Aut}_{X/k}$ over $k$.  The next proposition shows that the  Lie algebra of $\mathrm{Aut}_{X/k}$ is isomorphic to the Lie algebra of vector fields of $X$.

\begin{proposition}\label{Lie of Aut}
Let $X$ be an algebraic scheme defined over a field $k$ of characteristic $p\geq 0$ and $G=\mathrm{Aut}_{X/k}$ its automorphism group scheme. Then there exists a natural isomorphism of Lie algebras
\[
\Phi \colon \mathfrak{Lie}(\mathrm{Aut}_{X/k}) \rightarrow  \mathrm{Der}_k(X),
\]
Where $\mathrm{Der}_k(X)$ is the Lie algebra of vector fields of $X$. Moreover, if $p>0$, $\Phi$ is an isomorphism of restricted Lie algebras.
\end{proposition}
\begin{proof}
Let $A=k[t]/(t^2)$. Then there exist the following natural isomorphisms of $k$-vector spaces
\begin{gather*}
 \mathfrak{Lie}(\mathrm{Aut}_{X/k}) \xrightarrow{\tau} \mathrm{Hom}_k^e\left(\mathrm{Spec} A,\mathrm{Aut}_{X/k}\right)\xrightarrow{\sigma} \mathrm{Aut}_A(X\times \mathrm{Spec}A)\xrightarrow{\phi} \\
 \xrightarrow{\phi} \mathrm{Der}_k(X)\cong \mathrm{Hom}_X(\Omega_{X/k},\mathcal{O}_X)
\end{gather*}
Where $\tau$ is the isomorphism from Proposition~\ref{left invariant and tangent vectors} and Corollary~\ref{corollary 10}, $\sigma$ is from the definition of $\mathrm{Aut}_{X/k}$ as a group valued functor and its universal property and $\phi$ is the isomorphism explained explicitly in the proof of Proposition~\ref{left invariant and tangent vectors}.

The reader may now verify easily by following the definitions of each map that the resulting isomorphism of $k$-vector spaces between $ \mathfrak{Lie}(\mathrm{Aut}_{X/k}) $ and $\mathrm{Der}_k(X)$ is actually an isomorphism of Lie algebras (an restricted in positive characteristic).

\end{proof}

Next we give some examples of Lie algebras of certain group schemes with emphasis on positive characteristic cases when the Lie algebras are restricted.

\begin{example}[\cite{Milne17} Page 188]
Let $GL(n,k)$ and $SL(n,k)$ be the affine group schemes of invertible matrices and matrices with determinant $1$ over a field $k$ of characteristic $p>0$. As group valued functors they are defined as follows. Let $A$ be any finitely generated $k$-algebra. Then 
$\mathrm{Hom}_k(\mathrm{Spec} A, GL(n,k))$ is the group of $A$-automorphisms of $A^n$ while $\mathrm{Hom}_k(\mathrm{Spec} A, SL(n,k))$ is the group of $A$-automorphisms of $A^n$ with determinant $1$.

Then $\mathfrak{Lie}(GL(n,k))\cong \mathfrak{gl}_n(k)$ and $\mathfrak{Lie}(SL(n,k))\cong \mathfrak{sl}_n(k)$, where $\mathfrak{gl}_n(k)$ and $\mathfrak{sl}_n(k)$ are the restricted Lie algebras of matrices and matrices with trace zero, respectively. In both cases, the bracket is the bracket of matrices and the $p$-map is $p$-powers of matrices.

\end{example}

\begin{example}[The Lie algebra of $\mu_p$]\label{Lie algebra mp}
Let $k$ be a field of characteristic $p>0$. In this example it will be shown that $\mathfrak{Lie}(\mu_p)\cong\mathfrak{gl}_1(k)$, where $\mathfrak{gl}_1(k)=k$ with trivial bracket and $p$-map the Frobenius.

Let $D$ be a left invariant vector field of $\mu_p$ in characteristic $p>0$. Then $D$ is a $k$-derivation of $k[t]/(t^p-1)$ such that the following diagram commutes
\[
\xymatrix{
\frac{k[t]}{(t^p-1)} \ar[rr]^D\ar[d]^{m^{\ast}} & & \frac{k[t]}{(t^p-1)} \ar[d]^{m^{\ast}}\\
\frac{k[t]}{(t^p-1)} \otimes \frac{k[t]}{(t^p-1)} \ar[rr]^{1\otimes D} &&\frac{k[t]}{(t^p-1)} \otimes \frac{k[t]}{(t^p-1)}
}
\]
where $m^{\ast}(t)=t\otimes t$. Then $D=f(t)\frac{d}{dt}$, for some $f(t)\in k[t]/(t^p)$. Suppose that $f(t)=\sum_{i=0}^{p-1}a_it^i$, $a_i\in k$. Then from the above diagram it follows that 
\[
\sum_{i=1}^{p-1}a_i(t^i\otimes t^i) =\sum_{i=1}^{p-1}a_i (t\otimes t^i).
\]
This is possible only if $a_i=0$ for all $i\not=1$. Therefore $D=at\frac{d}{dt}$, $a\in k$. Let $\delta=t\frac{d}{dt}$. Then $\mathfrak{Lie}(\mu_p)=k\delta=\{a\delta,\;\; a\in k\}$. Moreover, $[a\delta,b\delta]=0$ and $\delta^p=\delta$. Hence $\mathfrak{Lie}(\mu_p)\cong\mathfrak{gl}_1(k)$.

\end{example}

\begin{example}[The Lie algebra of $\alpha_p$]\label{Lie algebra ap} 
Let $k$ be a field of characteristic $p>0$. Working identically as in the previous example, it follows that the $\mathfrak{Lie}(\alpha_p)=k\delta$, where $\delta=d/dt$. In this case $\delta^p=0$. Moreover, $[a\delta,b\delta]=0$, for any $a,b\in k$. Hence $\mathfrak{Lie}(\alpha_p)\cong k$, where $k$ is the trivial one dimensional restricted Lie algebra (trivial bracket and $p$-map).

\end{example}

\begin{example}[The Lie algebra of automorphisms of the Witt algebra~\cite{Salas00},~\cite{ScTz23}]
Let $k$ be a field of characteristic $p>0$ and $A=k[t]/(t^p-\omega)$, where $\omega \in k$. If $\omega\not\in k^p$ then $k\subset A$ is a purely inseparable field extension of height $1$. If on the other hand $\omega \in k^p$, then $A\cong k[t]/(t^p)$. To cover both cases we may assume that either $\omega \not\in k^p$ or $\omega=0$.

 Let $\mathrm{Aut}_{A/k}$ be the group scheme of automorphisms of $\mathrm{Spec} A$ over $k$.  The group scheme $\mathrm{Aut}_{A/k}$ can be described as follows. Let $R$ be a commutative $k$-algebra. Then $\underline{\mathrm{Aut}}_{A/k}(R)$ is the group of $R$-algebra automorphisms $\phi \colon R[t]/(t^p-\omega)\rightarrow R[t]/(t^p-\omega)$. 
Such an automorphism is determined by $\phi(t)$. Since $t^i$, $i=0,1,\ldots, p-1$ is an $R$-basis of $R[t]/(t^p-\omega)$, $\phi(t)=\sum_{i=1}^{p-1}\lambda_i t^i$, for some $\lambda_i\in R$. Since $t^p=\omega$ in $R[t]/(t^p-\omega)$, it follows that 
\begin{equation}\label{eq 10}
\lambda_0^p+\omega (\lambda_1-1)^p +\omega^2\lambda_2^p+\cdots +\omega^{p-1}\lambda_{p-1}^p=0.
\end{equation}
Moreover, since $\phi$ must be invertible it follows that the determinant of the map $\phi$ corresponding to the basis $t^i$, $i=0,1,\ldots, p-1$ must be non zero. From this it follows that $\mathrm{Aut}_{A/k}$ is isomorphic to the subscheme of $\mathbb{A}^p_k$ defined as the intersection of the equation (\ref{eq 10}) above and the open subset of $\mathbb{A}^p_k$ defined by the determinant. The reader can find a detailed description of these in~\cite{Salas00} and~\cite{ScTz23}.

Suppose that $\omega=0$. Then $\mathrm{Aut}_{A/k}$ is not reduced. Suppose on the contrary that $\omega\not\in k^p$. Then $\mathrm{Aut}_{A/k}$ is reduced. However it is  not geometrically reduced. In particular, in both cases $\mathrm{Aut}_{A/k}$ is not a smooth group scheme.

Let $\mathfrak{g}$ be the restricted Lie algebra of $\mathrm{Aut}_{A/k}$. This is called the Jacobson-Witt algebra. It is a restricted simple Lie algebra and of fundamental importance to the classification of restricted simple Lie algebras in positive characteristic~\cite{BlWi88}. By Proposition~\ref{Lie of Aut} it follows that $\mathfrak{g}\cong \mathrm{Der}_k(A)$. The $k$-derivations $t^i\delta$, where $\delta=d/dt$, $i=0,\dots, p-1$ form a  $k$-basis of $\mathfrak{g}$. The Lie bracket is given by the formulas
\[
[t^i\delta,t^j\delta]=(j-i)t^{i+j-1}\delta.
\]
Using this relation for $i=0$ and $j=0$ as well as the relation $[t^{p-1}\delta,t\delta]=(2-p)t^{p-1}\delta$ one easily sees that if $p\not= 2$, then $\mathfrak{g}$ is simple.

Consider in particular the cases when $p=2$ and $p=3$. 

Suppose that $p=2$. Then $\delta_1=\delta$ and $\delta_2=t\delta$ form a $k$-basis of $\mathfrak{g}$. Moreover, $\delta_1^2=0$, $\delta_2^2=\delta_2$ and $[\delta_1,\delta_2]=\delta$. From this description it easily follows that $\mathfrak{g}\cong k \rtimes \mathfrak{gl}_1(k)$. In particular it is not simple.

Suppose that $p=3$. Then $\delta_1=\delta$, $\delta_2=t\delta$ and $\delta_3=t^2\delta$ is a $k$-basis of $\mathfrak{g}$.  Moreover, $\delta_1^3=\delta_3^3=0$,  $\delta_2^3=\delta_3$, $[\delta_1,\delta_2]=\delta_1$, $[\delta_1,\delta_3]=2\delta_2$ and $[\delta_2,\delta_3]=\delta_3$. From this description it easily follows that $\mathfrak{g}\cong \mathfrak{sl}_2(k)$.
\end{example}

\subsection{Infinitesimal group schemes}
The purpose of this section is to present some results about finite affine group schemes over a field of  positive characteristic and to explain how their structure is related to the properties of its restricted Lie algebra. Standard references for the material that will be presented are~\cite{DG70} and~\cite{Milne17}.
\begin{definition}
Let $G=\mathrm{Spec} A$ be a finite group scheme over a field $k$ of characteristic $p>0$. The dimension of $A$ as a $k$-vector space is called the length of $G$ and is denoted by $|G|$.
\end{definition} 

Group schemes of length $p$ over an algebraically closed field are completely classified.

\begin{proposition}[\cite{OortTate70}]
Let $G$ be a finite group scheme of order $p$ over an algebraically closed field $k$ of characteristic $p>0$. Then $G$ is isomorphic to one of the following. The constant group scheme $\mathbb{Z}/p\mathbb{Z}$, the multiplicative group scheme $\mu_p$ or the additive group scheme $\alpha_p$.  In particular, $G$ is commutative.
\end{proposition} 
If $k$ is not algebraically closed this result is not true as there exist constant group schemes of length $p$ with less than $p$ points and therefore not isomorphic to $\mathbb{Z}/p\mathbb{Z}$. For example, $G=\mathrm{Spec}A$, where $A=k[x]/(x^p+\omega x)$, where $\omega \not\in k^p$.  

\begin{remark}
In contrast to the characteristic zero case where every group scheme of order $p^2$ is commutative, there exist non commutative group schemes of order $p^2$ over a field of characteristic $p>0$ as was shown in Example~\ref{non commutative of order 2p}. This exhibits how much richer the theory of group schemes in positive characteristic is compared to the one in characteristic zero.
\end{remark}

Next we define the notion of height of a finite group scheme.

\begin{definition}
Let $G=\mathrm{Spec} A$ be a finite group scheme with just one point over a field $k$ of characteristic $p>0$. The height of $G$ is the least positive integer $e$ such that $a^{p^e}=0$, for all $a\in \mathbf{m}_A$, where $\mathbf{m}_A$ is the unique maximal ideal of $A$. 
\end{definition}
The next proposition describes the scheme structure of a connected finite group scheme.
\begin{proposition}[\cite{Milne17}, Theorem 11.29]\label{height 1 structure}
Let $G=\mathrm{Spec} A$ be a connected finite group scheme over a perfect field of characteristic $p>0$. Then 
\[
A\cong k[t_1,\ldots, t_n]/(t_1^{p^{e_1}},\ldots, t_n^{p^{e_n}}),
\]
for some $e_1,\ldots, e_n \geq 1$. In particular, if $G$ has height $1$, then
\[
A\cong k[t_1,\ldots, t_n]/(t_1^p,\ldots, t_n^p).
\]
\end{proposition}
If the field $k$ is not perfect then the previous result is not correct as the following example shows.
\begin{example}[\cite{Wat79} Chapter 14, Exercise 1]]
Let $G=\mathrm{Spec} A$, where $A=k[x,y]/(x^{4},y^2-ax^2)$, where $a\not\in k^2$ and $k$ has characteristic $2$.  The $k$-algebra homomorphism $m^{\ast}\colon  A \rightarrow A\otimes A$ defined by $m^{\ast}(x)=x\otimes 1 + 1\otimes x$ and $m^{\ast}(y)=y\otimes 1+ 1\otimes y$ defines a group scheme structure on $G$. In fact it exhibits it as a closed group subscheme of $\mathbb{G}_a^2=\mathbb{G}_a\times \mathbb{G}_a$.

Let now $V=\{f\in A,\;\; f^2=0\}$. Suppose that $A$ has the form in Proposition~\ref{height 1 structure}. Then, since $\dim_k A=8$, $A$ would be isomorphic to $k[s,t]/(s^4,t^2)$. If this is the case then it is not hard to see that $\dim_k V=6$. However, we will show that this is not true. 

Let $f\in V$. Then we may write $f=f_0(x)+f_1(x)y$, where $f_i\in k[x]$ of degree at most $3$, $i=0,1$. Suppose that $f_0(x)=a_0+a_1x+a_2x^2+a_3x^3$ and $f_1(x)=b_0+b_1x+b_2x^2+b_3x^3$. Then
\[
f^2=f_0^2(x)+f_1^2(x)y^2=f_0^2(x)+f_1^2(x)(y^2-ax^2)+ax^2f^2_1(x).
\]
Since $f^2\in (x^4,y^2-ax^2)$, it follows that $f_0^2(x)+ax^2f_1^2(x)\in (x^4)$. Therefore, $a_0^2+a_1^2x^2+ab_0^2x^2\in (x^4)$ and hence $a_0^2+(a_1^2+ab_0^2)x^2=0$. This implies that $a_0^2=a_1^2+ab_0^2=0$. Now since $a\not\in k^2$, we get that $a_0=a_1=b_0=0$. Hence $\dim V=5$, a contradiction.
\end{example}

Let $G$ be a group scheme defined over a field $k$ of characteristic $p>0$. To $G$ there is an associated ascending filtration of infinitesimal group schemes which can be used in order to study the group scheme structure of  $G$. The filtration is derived from the notion of Frobenius Kernel which will be defined next.

\begin{definition}
Let $G$ be a group scheme defined over a field $k$ of characteristic $p>0$. Let $F^{(p^n)} \colon G\rightarrow G^{(p^n)}$ be the $n$-iterated geometric Frobenius morphism. $G^{(p^n)}$ is a group scheme and the morphism $F^{(p^n)}$ is a morphism of group schemes.  The kernel of $F^{(p^n)}$  is called the $n$-Frobenius Kernel of $G$ and is denoted by $G[nF]$. 
\end{definition}
The $n$-Frobenius kernel $G[nF]$ is an infinitesimal group scheme and from its definition it is not hard to see that it can also  be described as follows. Let $\mathcal{O}_{G,e}$ be the local ring of $G$ at the unit  point $e\in G$ and $\mathbf{m}_e$ its maximal ideal. Let $\mathbf{m}_e^{[n]}=\{a^{p^n}, \;\; a\in \mathbf{m}_e\}$. Then 
\[
G[nF]=\mathrm{Spec} \left(\mathcal{O}_{G,e}/\mathbf{m}_e^{[n]}\right).
\]
\begin{remarks}
\begin{enumerate}
\item From the definition of Frobenius kernel it is clear that for any group scheme $G$ over a field of characteristic $p>0$ , there exists a filtration of infinitesimal closed group subschemes of $G$,
\[
G[F]\subset G[2F] \subset \cdots \subset G[nF]\subset G[(n+1)F] \subset \cdots \subset G.
\]
\item The height of a finite group scheme $G$ is the smallest natural number $n\geq 1$ such that $G[nF]=G$, i.e., the Frobenius map $G\rightarrow G^{(p^n)}$ is trivial.
\item Since the tangent space at the identity of $G[nF]$ is the same as the tangent space at the identity of $G[F]$ it follows that $\mathfrak{Lie}(G[F])=\mathfrak{Lie}(G[nF])$, for any $n\geq 1$. This moreover shows that (in the case that $G$ is positive dimensional) that there may exist infinitely many non isomorphic group schemes with the same restricted Lie algebra.
\end{enumerate}
\end{remarks}

\begin{example}
Let $\mathbb{G}_a$, $\mathbb{G}_m$ be the additive and multiplicative  group schemes in characteristic $p>0$. 

The $n$-Frobenius Kernel $\mathbb{G}_a[nF]$ is the infinitesimal group scheme which is defined as a functor of groups as follows. Let $X$ be any $k$-scheme. Then from the definition of $n$-Frobenius kernel it follows that 
\[
\mathrm{Hom}_k(X,\mathbb{G}_a[nF])=\mathrm{Ker}[F^{(p^n)} \colon k[x]  \rightarrow H^0(X,\mathcal{O}_X)],
\]
where $F^{(p^n)}(f(x))=(f(x))^p$. Therefore,
\[
\mathrm{Hom}_k(X,\mathbb{G}_a[nF])=\{a\in H^0(X,\mathcal{O}_X),\;\; a^{p^n}=0\}.
\]
Hence $\mathbb{G}_a[nF]=\alpha_{p^n}$.

Similarly it can be shown that $\mathbb{G}_m[nF]=\mu_{p^n}$.
\end{example}

\subsection{The Demazure-Gabriel correspondence}
Let $k$ be a field of characteristic $p>0$. The purpose of this section is to show that there exists a one to one correspondence between infinitesimal group schemes of height $1$ over $k$ and finite dimensional restricted Lie algebras over $k$. In particular, a group scheme of height $1$ is completely determined by its restricted Lie algebra of derivations. This establishes a link between the representation theory of restricted Lie algebras and the algebraic theory of group schemes. The main references for this section is~\cite{DG70}, ~\cite{Milne17} and~\cite{Mu70}.

The first step in constructing the correspondence between restricted Lie algebras and height $1$ group schemes over $k$ is the construction of the universal algebra of a restricted Lie algebra $\mathfrak{g}$. 

Let $A$ be a $k$-algebra. Then $A$ is a restricted Lie algebra over $k$ by defining bracket by $[a,b]=ab-ba$, for all $a,b \in A$ and $p$-map by $a^{[p]}=a^p$, for all $a\in A$.
\begin{definition}
Let $\mathfrak{g}$ be a restricted Lie algebra over $k$. A universal restricted algebra of $\mathfrak{g}$ is a $k$-algebra $A$ together with a map of restricted Lie algebras $\rho \colon \mathfrak{g}\rightarrow A$ which satisfies the following universal property. Let $B$ be any other $k$-algebra and $\tau\colon \frak{g}\rightarrow B$ a map of restricted $k$-algebras. Then there exists a unique $k$-algebra homomorphism $f\colon A \rightarrow B$ such that the following diagram commutes.
\[
\xymatrix{ 
\mathfrak{g}\ar[dr]^{\tau}\ar[r]^{\rho} & A\ar[d]^f\\
& B
}
\]
\end{definition}
For any restricted Lie algebra $\mathfrak{g}$ its universal restricted algebra exists and can be constructed as follows. Let $T(\mathfrak{g})=\oplus_{s\geq 0} \mathfrak{g}^{[s]}$ be the tensor algebra of $\mathfrak{g}$ viewed as a $k$-vector space, where $\mathfrak{g}^{[s]}$ is the tensor product of $k$-vector spaces $\mathfrak{g}\otimes_k \cdots \otimes_k \mathfrak{g}$, $s$- times. Let $\rho \colon \mathfrak{g} \rightarrow T(\mathfrak{g})$ be the map defined by $\rho(x) =x$, for any $x\in \mathfrak{g}$, where by $x$ in the right hand side we mean the element $x\in \mathfrak{g}^{[1]}=\mathfrak{g}$ as an element of grade $1$ in $T(\mathfrak{g})$.  Let $I \subset  T(\mathfrak{g})$ de the two sided ideal generated by the elements $\rho([x,y])-\rho(x)\rho(y)+\rho(y)\rho(x)$ and $(\rho(x))^p-\rho(x^{[p]})$ of $T(\mathfrak{g})$, for all $x,y \in \mathfrak{g}$. Let now $U^{[p]}(\mathfrak{g})=T(\mathfrak{g})/I$. Then it is not hard to see that $\rho \colon \mathfrak{g}\rightarrow U^{[p]}(\mathfrak{g})$ is injective and that 

\begin{proposition}[\cite{Milne17},~\cite{DG70}]\label{Upg}
$U^{[p]}(\mathfrak{g})$ is the universal restricted Lie algebra of $\mathfrak{g}$. If in addition $\mathfrak{g}$ is finite dimensional, then $U^{[p]}(\mathfrak{g})$ is also finite dimensional as a $k$-vector space.
\end{proposition}
From now on and for the simplicity of notations,  we will identify $\mathfrak{g}$ with $\rho(\mathfrak{g})$ and write $x$ for $\rho(x)$ in $U^{[p]}(\mathfrak{g})$, for any $x\in \mathfrak{g}$. Any products will be products in $U^{[p]}(\mathfrak{g})$.

From the universality of $U^{[p]}(\mathfrak{g})$ there exist the following homomorphisms of $k$-algebras.
\begin{enumerate}
\item $m \colon U^{[p]}(\mathfrak{g}) \rightarrow U^{[p]}(\mathfrak{g})\otimes U^{[p]}(\mathfrak{g})$ such that $m(x)=x\otimes 1+1\otimes x$, for any $x\in \mathfrak{g}$.
\item $i \colon U^{[p]}(\mathfrak{g}) \rightarrow U^{[p]}(\mathfrak{g})$ such that $i(x)=-x$, for all $x\in \mathfrak{g}$.
\item $\varepsilon\colon U^{[p]}(\mathfrak{g}) \rightarrow k$ defined by $e(x)=0$, for all $x\in \mathfrak{g}$.
\end{enumerate}

Let  $(U^{[p]}(\mathfrak{g}))^{\ast}$ be the dual of $U^{[p]}(\mathfrak{g})$ as a $k$-vector space. Dualizing the  homomorphism $m$ above we obtain a homomorphism of $k$-vector spaces 

\[
m^{\ast}\colon (U^{[p]}(\mathfrak{g}))^{\ast}\otimes (U^{[p]}(\mathfrak{g}))^{\ast} \rightarrow (U^{[p]}(\mathfrak{g}))^{\ast}.
\]

From the definitions it easily follows that for any $f,g \in (U^{[p]}(\mathfrak{g}))^{\ast}$, 
\[
m^{\ast}(f\otimes g)(x)=f(1)g(x)+f(x)g(1)=m^{\ast}(g\otimes f)(x),
\]
for any $x\in \mathfrak{g}$. Therefore, considering that $\mathfrak{g}$ generates $U^{p]}(\mathfrak{g})$ as a $k$-algebra, it follows that $m^{\ast}(f\otimes g)=m^{\ast}(g\otimes f)$ and hence $m^{\ast}$ makes $(U^{[p]}(\mathfrak{g}))^{\ast}$ a commutative $k$-algebra.

Let  $\mu \colon U^{[p]}(\mathfrak{g}) \otimes U^{[p]}(\mathfrak{g}) \rightarrow U^{[p]}(\mathfrak{g})$ be the map defined by the multiplication of $U^{[p]}(\mathfrak{g})$, i.e., $\mu(x\otimes y)=xy$. Dualizing we get a $k$-algebra homomorphism 
$\mu^{\ast}\colon (U^{p]}(\mathfrak{g}))^{\ast}\rightarrow (U^{p]}(\mathfrak{g}))^{\ast}\otimes (U^{p]}(\mathfrak{g}))^{\ast}$. Let also $e^{\ast}\colon (U^{p]}(\mathfrak{g}))^{\ast}\rightarrow k$ be defined by $e^{\ast}(f)=f(1)$. Then after a careful verification by using the definitions it follows that
\begin{proposition}
The commutative $k$-algebra $(U^{[p]}(\mathfrak{g}))^{\ast}$ together with the homomorphisms $\mu^{\ast}$, $i^{\ast}$ and $e^{\ast}$ is a Hopf algebra. Therefore 
\[
G(\mathfrak{g})=\mathrm{Spec}\left(U^{[p]}(\mathfrak{g})\right)^{\ast}
\]
 is a group scheme with group scheme structure coming from the Hopf algebra structure of $(U^{[p]}(\mathfrak{g}))^{\ast}$.
\end{proposition}
(The reader may find the details in~\cite{Milne17} or~\cite{DG70}). 

Consider now the map 
\begin{equation}\label{phi map}
\phi \colon \mathfrak{g}\rightarrow \mathrm{Hom}_k((U^{p]}(\mathfrak{g}))^{\ast},k)
\end{equation}
defined as follows. Let $x\in \mathfrak{g}$. Then $\Phi(x)=\phi_x$, where  $\phi_x\colon (U^{p]}(\mathfrak{g}))^{\ast}\rightarrow k$ is defined by $\phi_x(f)=f(x)$, for all $f\in (U^{p]}(\mathfrak{g}))^{\ast}$. By the definition of product in $(U^{p]}(\mathfrak{g}))^{\ast}$, 
\[
\phi_x(f_1f_2)=f_1(1)f_2(x)+f_2(x)f_1(1)=e^{\ast}(f_1)\phi_x(f_2)+e^{\ast}(f_2)\phi_x(f_1).
\]
This shows that $\phi_x \colon (U^{p]}(\mathfrak{g}))^{\ast}\rightarrow k$ is a derivation.
 From Proposition~\ref{left invariant and tangent vectors} it follows that there exists a unique left invariant derivation $\Phi_x$ of $(U^{[p]}(\mathfrak{g}))^{\ast}$ restricting to $\phi_x$ at the identity point. Therefore there exists a map of restricted Lie algebras
\begin{equation}\label{eq DG 1}
\Phi \colon \mathfrak{g}\rightarrow \mathfrak{Lie}(G(\mathfrak{g}))
\end{equation}
defined by setting $\mathbf{\Phi}(x)=\Phi_x$.

\begin{proposition}[\cite{DG70},~\cite{Milne17}]
Let $\mathfrak{g}$ be a finite dimensional restricted Lie algebra over a field of characteristic $p>0$. Then 
\begin{enumerate}
\item The group scheme $G(\mathfrak{g})$ is a height $1$ group scheme.
\item The map $\Phi \colon \mathfrak{g}\rightarrow \mathfrak{Lie}(G(\mathfrak{g}))$ is an isomorphism of restricted Lie algebras.
\end{enumerate}
\end{proposition}
\begin{proof}
I will show the first part. The reader may find a detailed proof of the second in~\cite[Chapter II, 7]{DG70}.

Let $R=U^{[p]}(\mathfrak{g})$ and let $\mathbf{m}_{e^{\ast}}$ be the maximal ideal corresponding to the identity map $e^{\ast}\colon R^{\ast}\rightarrow k$. Let $x\in \mathbf{m}_{e^{\ast}}$. Hence $x(1)=e^{\ast}(x)=0$. To show that $G(\mathfrak{g})$ has height $1$ it suffices to show that $x^p=0$. 

Suppose on the contrary that $x^p\not= 0$. Let $\mathbf{m}_{\varepsilon}$ be the ideal of $R$ generated by $x\in \mathfrak{g}$. This is a maximal ideal and it is simply the kernel of the map $\varepsilon\colon R\rightarrow k$. Given then any $r\in R$, there exists a $y\in \mathfrak{g}$ and $\lambda\in k$ such that $r=y+\lambda$. 

Since $x^p\not=0$, there exists $r\in R$ such that $x^p(r)\not=0$. From the above discussion, there exist $y\in \mathfrak{g}$ and $\lambda \in k$ such that $r=y+\lambda$. Considering that $x(1)=0$, it follows that $x^p(y)\not=0$. But this means that, in the notation of equation (\ref{phi map}),  $\phi_y(x^p)=0$. But $\phi_y$ is a derivation. Therefore, $\phi_y(x^p)=px^{p-1}\phi(y)=0$, a contradiction. Hence $x^p=0$.

% and let $\Psi \colon R \rightarrow \mathrm{Hom}(R^{\ast},R^{\ast})$ be defined as follows. Let $r\in R$. Let $\psi_r \colon R^{\ast}\rightarrow k$ be defined by setting $\psi_r(f)=f(r)$. Then $\Psi_r$ is defined from the following diagram
%\[
%\xymatrix{
%R^{\ast} \ar[r]^{\mu^{\ast}} \ar[dr]^{\Psi_r} & R^{\ast}\otimes_k R^{\ast} \ar[d]^{\psi_r \otimes 1_{R^{\ast}}} \\
 %& R^{\ast}
%}
%\]
%Note that this is an extension of the map $\Phi$ in (\ref{eq DG 1}) in the whole $R$. If $r\in \mathfrak{g}$ then $\psi_r=\phi_r$ and $\Psi_r=\Phi_r$ is a $k$-derivation of $R^{\ast}$.

\end{proof}

The following proposition shows that morphisms between height $1$ group schemes are in one to one correspondence with homomorphisms between their restricted Lie algebras.
\begin{proposition}
Let $\mathfrak{g}_1$, $\mathfrak{g}_2$ be finite dimensional restricted Lie algebras over a field $k$ of characteristic $p>0$ and let $\phi \colon \mathfrak{g}_1 \rightarrow \mathfrak{g}_2$ be a homomorphism of restricted Lie algebras. Then there exists a unique morphism of group schemes $G(\phi) \colon G(\mathfrak{g}_1)\rightarrow G(\mathfrak{g}_2)$ such that the following diagram commutes.
\[
\xymatrix{
\mathfrak{Lie}(G(\mathfrak{g}_1)) \ar[rr]^{\mathfrak{Lie}(G(\phi))} & &\mathfrak{Lie}(G(\mathfrak{g}_1))\\
\mathfrak{g}_1\ar[u]^{\mathbf{\Phi}} \ar[rr]^{\phi} && \mathfrak{g}_2\ar[u]^{\mathbf{\Phi}}
}
\]
where $\Phi$ is the isomorphism (\ref{eq DG 1}).
\end{proposition}
\begin{proof}
By its definition, $G(\mathfrak{g}_i)=\mathrm{Spec} \left( U^{[p]}(\mathfrak{g}_i)\right)^{\ast}$, $i=1,2$. The map $\phi \colon \mathfrak{g}_1\rightarrow \mathfrak{g}_2$ induces a $k$-algebra 
homomorphism $U^{[p]}(\mathfrak{g}_1)\rightarrow U^{[p]}(\mathfrak{g}_2)$. Dualizing this map (as $k$-vector spaces) we get a $k$-linear map $(U^{[p]}(\mathfrak{g}_2))^{\ast} \rightarrow (U^{[p]}(\mathfrak{g}_1))^{\ast}$. By checking carefully the definition of the $k$-algebra structure of $(U^{[p]}(\mathfrak{g}_i))^{\ast}$ one sees that this map is actually a $k$-algebra homomorphism and the resulting morphism $G(\phi)\colon G(\mathfrak{g}_1)\rightarrow G(\mathfrak{g}_2)$ is a group scheme morphism.

The existence and commutativity of the diagram in the statement follows by a careful study of the definition of the $\mathfrak{Lie}(f)$ where $f\colon \mathfrak{g}_1\rightarrow \mathfrak{g}_2$ is a map of restricted Lie algebras as defined in Proposition~\ref{Lie to gs map}. The reader can find more details in~\cite[Chapter II, 7]{DG70}.

\end{proof}
Putting all the above together we get the Demazure-Gabrier Correspondence. A detailed exposition of this can be found in~\cite[Chapter II, 7]{DG70}. The case when $\mathfrak{g}$ is abelian is treated in detail in~\cite[Pages 124-134]{Mu70}.

\begin{theorem}[Demazure-Gabriel correspondence.~\cite{DG70},~\cite{Mu70}]\label{DG equivalence}
Let $k$ be a field of characteristic $p>0$. The functor $G\mapsto \mathfrak{Lie}(G)$ is an equivalence between the category of finite group schemes of height $1$ and the category of finite dimensional restricted Lie algebras. Its inverse is the functor 
$\mathfrak{g}\mapsto G(\mathfrak{g})$,where $\mathfrak{g}$ is a finite dimensional restricted Lie algebra over $k$.

In particular, let $\mathfrak{g}$ be a  finite dimensional restricted Lie algebra and $G$ a finite group scheme of height $1$, Then,

\begin{gather*}
\mathfrak{Lie}(G(\mathfrak{g}))\cong \mathfrak{g} \quadand  G(\mathfrak{Lie}(G))\cong G
\end{gather*}
\end{theorem}
\begin{remark}
According to the Demazure-Gabriel correspondence, the categories of finite group schemes of height $1$ and of finite dimensional restricted Lie algebras are equivalent. Therefore the classification of finite group schemes of height $1$ is equivalent to the classification of finite dimensional restricted Lie algebras. Restricted simple Lie algebras are classified~\cite{BlWi88}. However this classification does not provide a classification of simple height $1$ group schemes because the corresponding restricted Lie algebra of a simple group scheme is a simple restricted Lie algebra but may not be a simple Lie algebra. 
\end{remark}

We will finish the discussion about the Demazure-Gabriel correspondence with two explicit examples that will be useful in the next sections.
\begin{example}
Let $k$ be a field of characteristic $p>0$ and $\mathfrak{g}=k$ be the 1-dimensional restricted Lie algebra with trivial bracket and $p$-map. From the Example~\ref{Lie algebra ap} this is isomorphic to the Lie algebra of $\alpha_p$. We will explicitly construct the group scheme $G(\mathfrak{g})$ and show that it is isomorphic to $\alpha_p$. 

Since the bracket and $p$-map of $\mathfrak{g}$ are trivial, it easily follows that as $k$-algebras
\[
U^{[p]}(\mathfrak{g})= T(\mathfrak{g})/(xy-yx, x^p,\; x,y\in \mathfrak{g})=S(\mathfrak{g})/(x^p)\cong k[t]/(t^p).
\]
Then the map $m^{\ast}\colon A\rightarrow A\otimes_k A$, where $A=k[t]/(t^p)$ as defined after Proposition~\ref{Upg} is simply defined by $m^{\ast}(t)=t\otimes 1 +1\otimes t$. This makes $U^{[p]}(\mathfrak{g})$ into a Hopf algebra and hence 
$\mathrm{Spec}(U^{[p]}(\mathfrak{g}))$ is a group scheme. Moreover, the previous description of $m^{\ast}$ and $U^{[p]}(\mathfrak{g})$ show that
\[
\mathrm{Spec}(U^{[p]}(\mathfrak{g}))\cong \alpha_p.
\]
From this it follows that $G(\mathfrak{g})=\mathrm{Spec}\left(U^{[p]}(\mathfrak{g})\right)^{\ast}$ is the Cartier dual of $\alpha_p$ which is isomorphic to $\alpha_p$.
\end{example}
\begin{example}
Let $k$ be a field of characteristic $p>0$ and $\mathfrak{g}=\mathfrak{gl}_1(k)$ be the 1-dimensional restricted Lie algebra with trivial bracket and $p$-map the Frobenius map. Therefore $\mathfrak{gl}_1(k)$ is generated by an element $e$ such that $e^{[p]}=e$.
Then there are isomorphisms of $k$-algebras
\[
U^{[p]}(\mathfrak{g})= T(\mathfrak{g})/(xy-yx, x^p-x,\; x,y\in \mathfrak{g})=S(\mathfrak{g})/(x^p-x,\; x\in \mathfrak{g}) \cong [k[t]/(t^p-t).
\]
As in the previous example, the map $m^{\ast}\colon A\rightarrow A\otimes_k A$, where $A=k[t]/(t^p-t)$  is defined by $m^{\ast}(t)=t\otimes 1 +1\otimes t$ and it makes $U^{p]}(\mathfrak{g})$ into a Hopf algebra. Therefore $\mathrm{Spec}(U^{[p]}(\mathfrak{g}))$ is a group scheme. In fact from the description of $m^{\ast}$ this is simply the constant group scheme $(\mathbb{Z}/p\mathbb{Z})_k$. Hence $G(\mathfrak{g})$ is its Cartier which by Example~\ref{Cartier dual of mp} is isomorphic to $\mu_p$. 

\end{example}

\subsection{Actions of infinitesimal group schemes}

In this section we will discuss several useful properties, unique in positive characteristic, of an action of an infinitesimal group scheme $G$ on an algebraic variety $X$ over a field of characteristic $p>0$. Special emphasis will be given to the case when $G$ has height $1$ and in particular when $G=\alpha_p$ or $G=\mu_p$. These cases are of particular interest since $\alpha_p$ and $\mu_p$ may appear as subgroup schemes of the automorphism group scheme of a variety. The main references for this section are~\cite{Mu70},~\cite{MuFoKi93},~\cite{Milne17}. 

\begin{definition}
An action of a group scheme $G$ on a scheme $X$ over a field $k$ is a morphism $\mu \colon G\times X\rightarrow X$ such that
\begin{enumerate}
\item The composite 
\[
X\cong \mathrm{Spec} k\times X \overset{e\times 1_X}{\rightarrow} G\times X \overset{\mu}{\rightarrow} X
\]
is the identity.
\item The diagram
\[
\xymatrix{
G\times G \times G \ar[rr]^{m\times 1_X} \ar[d]^{1_G \times \mu} && G\times X\ar[d]^{\mu} \\
G\times X \ar[rr]^{\mu} && X
}
\]
is commutative, where $m\colon G\times G \rightarrow G$ and $e\colon \mathrm{Spec} k\rightarrow G$ is the multiplication and identity morphisms of $G$.
\end{enumerate}
\end{definition}
In the language of functor of points, the previous definition is equivalent to the following. A group scheme $G$ acts on a scheme $X$ if and only if for any scheme $S$ over $k$, $\mathrm{Hom}(S,G)$ acts on $\mathrm{Hom}(S,X)$ functorially, which means that for any morphism $f\colon S \rightarrow T$ of schemes over $k$, the following diagram commutes.
\[
\xymatrix{
\mathrm{Hom}(T,X) \ar[rr] && \mathrm{Hom}(S,X)\\
\mathrm{Hom}(T,G)\times \mathrm{Hom}(T,X) \ar[u]^{\mu_T}\ar[rr] & & \mathrm{Hom}(S,G)\times \mathrm{Hom}(S,X) \ar[u]^{\mu_S}
}
\]
where $\mu_T$ and $\mu_S$ define the action of $\mathrm{Hom}(T,G)$ and $\mathrm{Hom}(S,G)$ on $\mathrm{Hom}(T,X) $ and $\mathrm{Hom}(S,X)$, respectively. 

Suppose that $G$ acts on $X$. Consider then the morphism
\[
(\mu,p_2)\colon G\times X \rightarrow G\times X,
\]
where $p_2\colon G\times X \rightarrow X$ is the projection. We say that the action is proper or  free if and only if the morphism $(\mu,p_2)$ is either proper or a closed immersion, respectively. 

Next we will recall the notions of $G$-invariant maps and quotient of a scheme by a group~\cite{Mu70},~\cite{MuFoKi93},~\cite{Kollar97}.

\begin{definition}
Let $G$ be a group scheme acting on an algebraic scheme $X$ over a field $k$. A morphism $f\colon X\rightarrow Z$ is called $G$-invariant if and only if the following diagram commutes.
\[
\xymatrix{
G\times X\ar[r]^{\mu}\ar[d]_{p_2} & X\ar[d]^f\\
X\ar[r]^f & Z
}
\]
In particular, if $Z=\mathbb{A}^1_k$, $f$ is called a $G$-invariant function.
\end{definition}

\begin{definition}[\cite{Kollar97}]
Let $G$ be a group scheme acting on an algebraic scheme $X$ over a field $k$. Let also $\pi\colon X\rightarrow Z$ be a $G$-invariant morphism.
\begin{enumerate}
\item The morphism $\pi \colon X\rightarrow Z$ is called a categorical quotient of $X$ by $G$ if and only if for any other $G$-invariant morphism $g \colon X\rightarrow Y$, there exists a unique morphism $h\colon Z\rightarrow Y$ such that $h=g\circ \pi$.
\item The morphism $\pi\colon X\rightarrow Z$ is called a topological quotient of $X$ by $G$ if and only if the following conditions are satisfied.
\begin{enumerate}
\item For any algebraically closed field $E$ the map 
\[
f(E)\colon \underline{X}(E)/\underline{G}(E)\rightarrow \underline{Z}(E)
\]
 is an isomorphism of sets.
\item The morphism $\pi$ is universally submersive.
\end{enumerate}
\item The morphism $\pi$ is called a geometric quotient of $X$ by $G$ if and only if it is a topological quotient and in addition,
\[
\mathcal{O}_Z=(\pi_{\ast}\mathcal{O}_X)^G,
\]
the sheaf or rings of $G$-invariant functions of $X$. 
\item All three kinds of quotients are unique if they exist.
\end{enumerate}
\end{definition}

The next theorem establishes the existence of a geometric quotient under certain conditions in the case of finite group schemes. More general results can be found in~\cite{MuFoKi93},~\cite{Kollar97}.

\begin{theorem}[Theorem 1, page 104~\cite{Mu70},~\cite{Kollar97}]\label{existence of quot}
Let $G$ be a finite group scheme acting on a scheme $X$ defined over a field $k$. Then,
\begin{enumerate}
\item Suppose that the orbit of any point is contained in an affine open subset of $X$. Then the geometric quotient $\pi \colon X\rightarrow Z$ exists. Moreover, $\pi$ is a finite morphism of degree $n=|G|$.
\item Suppose that $k$ has characteristic $p>0$. Then the geometric quotient $\pi \colon X\rightarrow Z$ exists as an algebraic space~\cite{Kollar97}.
\end{enumerate}
Suppose that either of the above conditions is satisfies and therefore the geometric quotient $\pi \colon X \rightarrow Z$ exists. Suppose in addition that $G$ acts freely on $X$, then the morphism $\pi$ is flat of degree $n=|G|$.
\end{theorem}

\begin{remarks}
\begin{enumerate}
\item Suppose that $G$ is an infinitesimal finite group scheme, i.e., $G$ is a single point as a set. Then since every orbit is a point as a set, the condition (1) above  is satisfied automatically and hence in this case the geometric quotient exists. Moreover, every open subset $U$ of $X$ is $G$-invariant.
\item Suppose that $G=\mathrm{Spec} R$. Then according to the proof of part $(1)$ of the previous theorem found in~\cite[Theorem 1, Page 104]{Mu70}, locally the morphism $\pi \colon X\rightarrow Z$ is described as follows. Suppose $X=\mathrm{Spec} A$. Then the map $\mu \colon G \times X \rightarrow X$ defining the action of $G$ on $X$ corresponds to a homomorphism of $k$-algebras $\mu^{\ast}\colon A \rightarrow R\otimes_k A$. Then $Y=\mathrm{Spec} B$, where 
\[
B=\{a\in A,\;\; \mu^{\ast}a=1\otimes a\}.
\]
\end{enumerate}
\end{remarks}

We next consider the case of actions of group schemes of height $1$.
\begin{proposition}\label{height 1 actions}
Let $G$ be a finite group scheme of height $1$ over a field $k$ of characteristic $p>0$. To give an action of $G$ on a scheme $X$ over $k$  it is equivalent to give a map of restricted Lie algebras
\[
\mathfrak{Lie}(G)\rightarrow \mathrm{Der}_k(X),
\]
where $\mathrm{Der}_k(X)$ is the restricted Lie algebras of vector fields of $X$.
\end{proposition}
\begin{proof}
Let $\mathrm{Aut}_{X/k}$ be the automorphism group scheme of $X$. To give an action of $G$ on $X$ is equivalent to give a morphism of group schemes $\phi\colon G \rightarrow \mathrm{Aut}_{X/k}$. This map induces a map of Frobenius Kernels 
\[
\phi[F] \colon G[F] \rightarrow \mathrm{Aut}_{X/k}[F].
\]
Since $G$ has height $1$, $G[F]=G$. Then by the Demazure-Gabriel correspondence, Theorem~\ref{DG equivalence}, $\phi[F]$ is determined completely by the map between the  restricted Lie algebras of $G$ and $\mathrm{Aut}_{X/k}[F]$ that it induces. Taking into consideration that $\mathfrak{Lie}(\mathrm{Aut}_{X/k}[F])= \mathfrak{Lie}(\mathrm{Aut}_{X/k})$ and that by Proposition~\ref{Lie of Aut}, $\mathfrak{Lie}(\mathrm{Aut}_{X/k})\cong \mathrm{Der}_k(X)$, we get the statement of the proposition.
\end{proof}
The previous proposition makes it possible to obtain a description of the geometric quotient of a scheme by a height $1$ group scheme by using vector fields. Note that the geometric quotient exists by the first remark after Theorem~\ref{existence of quot}. 

We first define a special kind of purely inseparable map constructed by rational vector fields of a variety.
\begin{definition}\label{definition 500}
Let $X$ be an integral scheme of finite type over a field $k$ of characteristic $p>0$. Let $\mathfrak{h} \subset \mathrm{Der}_k(k(X))$ be a restricted sub Lie algebra of the restricted Lie algebra of derivations of the function field $k(X)$ of $X$. Let $D_1,\ldots, D_m$ be rational vector fields of $X$ forming a $k(X)$-basis of $\mathfrak{h}$. Let $Y^{\mathfrak{h}}$ be the  scheme whose underlying topological space is $X$ and its sheaf of rings is defined as follows. Let $U\subset X$ open. Then 
\[
\mathcal{O}_{Y^{\mathfrak{h}}}(U)=\{a\in \mathcal{O}_X(U), \;\; D_1(a)=\cdots =D_m(a)=0\} \subset \mathcal{O}_X(U).
\]
The resulting map $\pi \colon X\rightarrow X^{\mathfrak{h}}$ is called the quotient of $X$ by $\mathfrak{h}$. By Theorem~\ref{Jacobson Correspondence} $\pi$ is purely inseparable of degree $p^n$. Sometimes, $X^{\mathfrak{h}}$ is denoted by $X^{D_1,\dots, D_m}$.
\end{definition}

\begin{proposition}[Proposition 4.2~\cite{ScTz23}]\label{height 1 actions second}
Let $G$ be a height $1$ group scheme over a field $k$ of characteristic $p>0$ acting on an integral scheme $X$. Then by Proposition~\ref{height 1 actions}, the action of $G$ on $X$ is determined by  a map between of restricted Lie algebras
\[
\phi \colon \mathfrak{Lie}(G)\rightarrow \mathrm{Der}_k(X).
\]
Let $D_1,\ldots, D_n$ be vector fields of $X$ forming a $k$-basis for the image of $\phi$. Let also $\pi \colon X\rightarrow Z$ be the geometric quotient of $X$ by $G$, which exists by Theorem~\ref{existence of quot}. Then  
\[
Z=X^{D_1,\dots, D_n}
\]
Moreover, $\deg(\pi)=p^m$, where $m$ is the dimension, as a $k(X)$-vector space, of the restricted Lie subalgebra of $\mathrm{Der}_k(X)$ generated by $D_1,\ldots, D_n$. Suppose in addition that $X$ is normal. Then $Y$ is also normal.
\end{proposition}

\subsection{Actions of $\alpha_p$ and $\mu_p$.}
In this section we will study in more detail actions of $\mu_p$ and $\alpha_p$ on an algebraic variety $X$ over a field of characteristic $p>0$ and the quotient of $X$ by this action. A more detailed treatment of this topic can be found in~\cite{Tz17}.

\begin{proposition}\label{explicit des of ap and mp action}
Let $X$ be a scheme of finite type over a field $k$ of characteristic $p>0$. Then $X$ admits a nontrivial $\alpha_p$ or $\mu_p$ action if and only if $X$ has a nontrivial global vector field $D$ such that $D^p=0$ or $D^p=D$, respectively. Moreover,
\begin{enumerate}
\item Let $D$ be a vector field of $X$ such that $D^p=0$. Then the action of $\alpha_p$ on $X$ is defined by the map
\[
\Phi \colon \mathcal{O}_X \rightarrow \frac{\mathcal{O}_X[t]}{(t^p)}
\]
where 
\[\Phi(a)=\sum_{k=0}^{p-1}\frac{D^ka}{k!}t^k.
\]
 \item Suppose that $X$ admits a $\mu_p$ action induced by a vector field $D$ of $X$ such that $D^p=D$.Then the action is defined by the  map
\[
\Psi \colon \mathcal{O}_X \rightarrow \frac{\mathcal{O}_X[t]}{(t^p-1)}
\]
where 
\[
\Psi(a)=\sum_{k=0}^{p-1}\frac{D^ka}{k!}+\sum_{k=1}^{p-1}\frac{D^ka}{k!}t^k.
\]
\end{enumerate}
\end{proposition}
\begin{proof}

Let $G$ be either $\alpha_p$ or $\mu_p$. Since $G$ is a group scheme of height $1$, from Proposition~\ref{height 1 actions}, an action of $G$ on $X$ is completely determined by a map of restricted Lie algebras
\begin{equation}\label{eq 50}
\phi \colon \mathfrak{Lie}(G) \rightarrow \mathrm{Der}_k(X).
\end{equation}
By Examples~\ref{Lie algebra mp} and ~\ref{Lie algebra ap}, $\mathfrak{Lie}(\mu_p)$ and $\mathfrak{Lie}(\alpha_p)$ are one dimensional over $k$ generated by an element $e$ such that either $e^{[p]}=e$ or $e^{[p]}=0$, respectively. Let $D=\phi(e)$. Then either $D^p=D$ or $D^p=0$. Conversely, if such a vector field exists then there exists a map of restricted Lie algebras as in (\ref{eq 50}) above and hence an action of $\alpha_p$ or $\mu_p$ on $X$.

The above argument is based on the Demazure-Gabriel correspondence and even though short the form of the action is hiden. However, it is possible to write a simple self contained argument which describes the action explicitly and might be useful for calculations. We will only consider the case of an $\alpha_p$ action. The $\mu_p$ case is identical and is omitted.

Suppose that $X$ admits a nontrivial $\alpha_p$-action. Let \[
\mu \colon \alpha_p \times X \rightarrow X
\]
be the map that defines the action. Let 
\[
\mu^{\ast} \colon \mathcal{O}_X \rightarrow \mathcal{O}_X \otimes_k \frac{k[t]}{(t^p)}
\]
be the corresponding map on the sheaf of rings level. The additive group $\alpha_p$ is $\mathrm{Spec}\frac{k[t]}{(t^p)}$ as a scheme with group scheme structure given by
\[
m^{\ast}\colon \frac{k[t]}{(t^p)} \rightarrow \frac{k[t]}{(t^p)} \otimes \frac{k[t]}{(t^p)} \]
defined by $m^{\ast}(t)=1\otimes t + t \otimes 1$. Then by the definition of group scheme action, there is a commutative diagram.
\[
\xymatrix{
\mathcal{O}_X \ar[r]^{\mu^{\ast}} \ar[d]^{\mu^{\ast}} & \mathcal{O}_X \otimes_k \frac{k[t]}{(t^p)} \ar[d]^{\mu^{\ast}\otimes 1} \\
\mathcal{O}_X \otimes_k \frac{k[t]}{(t^p)} \ar[r]^{1\otimes m^{\ast}} & \mathcal{O}_X \otimes_k \frac{k[t]}{(t^p)} \otimes_k \frac{k[t]}{(t^p)}\\
}
\]
The map $\mu^{\ast}$ is given by \[
\mu^{\ast}(a)= a\otimes 1 +\Phi_1(a)\otimes t +\sum_{k=2}^{p-1}\Phi_k(a)\otimes t^k,
\] 
where $\Phi_k \colon \mathcal{O}_X \rightarrow \mathcal{O}_X$ are additive maps. The fact that $\mu^{\ast}$ is a sheaf of rings map shows that $\Phi_1$ is a derivation which we call $D$. We will show by induction that $\Phi_k(a)=D^ka/k!$. and hence the first part of the proposition.

From the commutativity of the previous diagram and the definition of $\mu^{\ast}$ and $m^{\ast}$ it follows that
\begin{gather*}
a \otimes 1 \otimes 1 + Da \otimes t \otimes 1 +\sum_{k=2}^{p-1} \Phi_k(a) \otimes t^k \otimes 1+ Da \otimes 1 \otimes t +D^2a\otimes t \otimes t + \\
\sum_{k=2}^{p-1}\Phi_k(Da)\otimes t^k 
\otimes t +\\
 \sum_{k=2}^{p-1}\left( \Phi_k(a)\otimes 1 \otimes t^k +D(\Phi_k(a))\otimes t \otimes t^k +\sum_{s=2}^{p-1}\Phi_s(\Phi_k(a))\otimes t^s \otimes t^k\right) =\\
a\otimes 1 \otimes 1 +Da \otimes t \otimes 1 + Da \otimes 1 \otimes t +\sum_{k=2}^{p-1}\sum_{s=0}^{k}\Phi_k(a) \binom{k}{s}t^s \otimes t^{k-s}.
\end{gather*}
Equating the coefficients of $t^{k-1}\otimes t $ on both sides of the equation we get that $k\Phi_k(a)t^{k-1}\otimes t = \Phi_{k-1}(Da)\otimes t^{k-1} \otimes t$. By induction, 
$\Phi_{k-1}=D^{k-1}/(k-1)!$. Therefore $\Phi_k(a)=D^ka/k!$ as claimed. Moreover, equating the coefficients of $t \otimes t^{p-1}$ it follows that $D^p=0$.

\end{proof}

The previous proposition motivates the following definition.

\begin{definition}
Let $D$ be a vector field on a scheme $X$ defined over a field of characteristic $p>0$. $D$ is called additive if $D^p=0$ and multiplicative if $D^p=D$.
\end{definition}

\begin{example}
In this example we will exhibit an action of $\alpha _p$ on $\mathbb{P}^3_k$ in two different ways in characteristic $p>0$. One by using the functorial definition of action and another by using the previous proposition, i.e., by exhibiting an explicit additive vector field on $\mathbb{P}^3_k$.

To write an action of $\alpha_p$ on $\mathbb{P}^3_k$ is equivalent to define a map of group schemes $\phi \colon \alpha_p \rightarrow \mathrm{Aut}_{\mathbb{P}^3/k}$. This can be done by defining a map 
\[
\phi_R \colon \underline{\alpha}_p(R) \rightarrow  \underline{\mathrm{Aut}}_{\mathbb{P}^3}(R)
\]
for any finitely generated $k$-algebra $R$. Now $\underline{\alpha}_p(R)=\{r\in R,\;\; r^p=0\}$ and $\underline{\mathrm{Aut}}_{\mathbb{P}^3}(R)$ is the group of  $R$-automorphisms of $\mathbb{P}^3_R$. Let now $r\in R$ such that $r^p=0$. Define $\phi_R(r)$ to be the $R$ automorphism of $\mathbb{P}^3_R$ given by the matrix $A(r)\in \mathrm{M}_3(R)$, where
\[
A(r)=
\begin{pmatrix}
1 & 0 & 0 & r \\
0 & 1 & 0 & 0\\
0 & 0 & 1 & 0\\
0 & 0 & 0 & 1
\end{pmatrix}
\]
It is not now difficult to see that $\phi_R$ is a group homomorphism and that it is functorial in $R$. Therefore it induces a group scheme morphism  $\phi \colon \alpha_p \rightarrow \mathrm{Aut}_{\mathbb{P}^3/k}$. Note that in this way one can define various actions of $\alpha_p$ on $\mathbb{P}^3_k$ or in $\mathbb{P}^n_k$ in general.

Let now $X\subset \mathbb{P}^3_k$ be a surface given by a homogeneous polynomial $f$ of degree $d$. The matrix $A(r)$ defined above corresponds to the automorphism $\phi_R(r)$ of $\mathbb{P}^3_R$ given by 
\[
\phi_R(r)([x,y,z,w])=[x+rw,y,z,w].
\]
If $X$ is invariant under this automorphism of $\mathbb{P}^3_R$ for any $R$,  then the action of $\alpha_p$ on $\mathbb{P}^3_k$ induces a non trivial action of $\alpha_p$ on $X$. 

For example, one can easily check that this happens for
\[
f(x,y,z,w)=w^py+wy^p+wyz^{p-1}+zx^p. 
\]
Moreover, with  some effort one can show that the surface $X$ defined by $f$ is singular with only rational double points. If in addition $p \geq 5$, $X$ is a surface of general type  with a nontrivial $\alpha_p$ action.

The $\alpha_p$ action can be defined also from the the following additive vector field on $X$. Let $\delta \colon k[x,y,z,w] \rightarrow k[x,y,z,w]$ be the homogeneous $k$-derivation given by $D=y\partial/\partial x$. Then $\delta^p=0$. Since $\delta(f)=0$, $\delta$ induces a homogeneous $k$-derivation of $k[x,y,z,w]/(f)$ and therefore gives a vector field of $D$ on $X$ with $D^p=0$. Then by Proposition~\ref{explicit des of ap and mp action}, $D$ induces an $\alpha_p$ action on $X$.
\end{example}

\begin{example}
In this example we will exhibit in a way similar to the previous example, nontrivial actions of $\mu_p$ on $\mathbb{P}^3_k$ in characteristic $p>0$. The methods can be generalized to produce actions of $\mu_p$ in $\mathbb{P}^n_k$ for any $n\geq 1$ and also to  hypersurfaces $X\subset \mathbb{P}^n_k$.

Let $a_1$, $a_2$ , $a_3$, $a_4$ be positive integers. Let $R$ be any finitely generated $k$-algebra. Then $\underline{\mu}_p(R)=\{r\in R,\;\; r^p=1\}$. Let also
\[
\phi_R \colon \underline{\mu}_p(R) \rightarrow  \underline{\mathrm{Aut}}_{\mathbb{P}^3}(R)
\]
be the map defined as follows. Let $r\in R$ such that $r^p=1$. Then $\phi_R(r)$ is the $R$-automorphism of $\mathbb{P}^3_R$ given by the matrix $A(r)\in \mathrm{M}_3(R)$, where
\[
A(r)=
\begin{pmatrix}
r^{a_1} & 0 & 0 & 0 \\
0 & r^{a_2} & 0 & 0\\
0 & 0 & r^{a_3} & 0\\
0 & 0 & 0 & r^{a_4}
\end{pmatrix}
\]
It is not hard to see that $\phi_R$ is a group homomorphism and that it is functorial in $R$. Therefore it induces a group scheme morphism  $\phi \colon \mu_p \rightarrow \mathrm{Aut}_{\mathbb{P}^3/k}$ and hence an action of $\mu_p$ on $\mathbb{P}^3_k$. 

This action can alternatively be defined by a vector field of $\mathbb{P}^3_k$ of multiplicative type in the following way. Let $\delta \colon k[x_1,x_2,x_3,x_4]\rightarrow k[x_1,x_2,x_3,x_4]$ be the homogeneous derivation defined by
\begin{equation}\label{eq 51}
\delta=a_1x_1\frac{\partial}{\partial x_1}+a_2x_2\frac{\partial}{\partial x_2}+a_3x_3\frac{\partial}{\partial x_3}+a_4x_4\frac{\partial}{\partial x_4}
\end{equation}
It is not hard to see that $\delta^p=\delta$. Hence $\delta$ induces a vector field $D$ of $\mathbb{P}^3_k$ such that $D^p=D$. Then by Proposition~\ref{explicit des of ap and mp action}, $D$ induces an $\mu_p$ action on $X$.

Let $X\subset \mathbb{P}^3_k$ be a surface given by a homogeneous polynomial $f$ of degree $d$. If $X$ is invariant under this automorphism of $\mathbb{P}^3_R$ for any $R$,  then the action of $\mu_p$ on $\mathbb{P}^3_k$ induces a non trivial action of $\mu_p$ on $X$. We will give two such examples of special interest.

Let $p=2$,  $a_1=a_2=1$, $a_3=a_4=0$ and 
\[
f(x_1,x_2,x_3,x_4)=x_1^4+x_2^4+x_3^4+x_1x_2x_3x_4. 
\]
Then $X$ is a quartic in $\mathbb{P}^3_k$ with only rational double points. Such a quartic is called a generalized $K3$ surface~\cite{Mats23a},~\cite{Mats23b}. $X$ is invariant under the $\mu_2$ action defined with the weights $a_1, a_2, a_3, a_4$ above and therefore the action of $\mu_2$ on $\mathbb{P}^3_k$ induces a non trivial $\mu_2$ action on $X$. Note that in this case the vector field $D$ on $X$ which induces the $\mu_2$ action on $X$ is obtained from the $k$-derivation
\[
\delta=x_1\frac{\partial}{\partial x_1}+x_2\frac{\partial}{\partial x_2}.
\]
of $k[x_1,x_2,x_3,x_4]$. This derivation has the property that $\delta^2=\delta$ and moreover $\delta(f)=0$. Therefore $\delta$ induces a vector field $D$ on $X$ such that $D^2=D$. 

Let $p=5$, and  $a_i=i$, $i=1,2,3,4$. Let also
\[
f(x_1,x_2,x_3,x_4)=\sum_{i_1+i_2+i_3+i_4=5}\lambda_{i_1,i_2,i_3,i_4}x_1^{i_1}x_2^{i_2}x_3^{i_3}x_4^{i_4}
\]
where $i_1+2i_2+3i_3+4i_4=0 \mod 5$. Then $X$ is a quintic surface in $\mathbb{P}^3_k$. $X$ is invariant under the $\mu_5$ action defined with the weights $a_1, a_2, a_3, a_4$ above and therefore the action of $\mu_5$ on $\mathbb{P}^3_k$ induces a non trivial $\mu_5$ action on $X$. From the point of view of vector fields again, this action is induced by the vector $D$ on $X$ which is defined from the derivation $\delta$ in (\ref{eq 51}) with $a_i=i$, $i=1,2,3,4$. This derivation has the property that $\delta(f)=0$ and hence induces a vector field $D$ on $X$ such that $D^5=D$.

The surfaces $X$ defined above are fundamental in the classification of Godeaux surfaces in characteristic $5$~\cite{La81}. By~\cite{La81},  the singularities of the general $X$ are exactly 11 rational double points and the quotient $Y=X/\mu_5$ is a smooth classical Godeaux surface. In fact by~\cite{La81} all classical Godeaux surfaces in characteristic $5$ are of this form. 

\end{example}

In the remaining part of this section we will present several results regarding quotients of a variety by $\alpha_p$ or $\mu_p$ actions. To begin we will define the divisorial and isolated part of a vector field, or more generally of a rational vector field of an algebraic scheme $X$.  In the literature~\cite{RS76}, the definitions in the case of a rational vector field are local for smooth varieties only. Here we will give a global definition for normal varieties.

\begin{definition}
Let $D$ be a vector field on an algebraic scheme $X$ defined over a field $k$. The fixed part of $D$ is the closed subscheme of $X$ defined by the ideal sheaf $(D(\mathcal{O}_X))$. Moreover,
\begin{enumerate}
\item The divisorial part of the fixed part of $D$ is the closed subscheme of $X$ defined by the height $1$ components in the primary decomposition of $(D(\mathcal{O}_X))$.
\item The isolated part of $D$ is the closed subscheme of $X$ defined by the height $\geq 2$ components in the primary decomposition of $(D(\mathcal{O}_X))$.
\end{enumerate}
\end{definition} 

\begin{definition}
Let $D$ be a rational vector field on a normal integral  algebraic scheme $X$ defined over a field $k$.  Then $D$ gives a $k$-derivation of the function field $k(X)$ of $X$. Let $Y\subset X$ be a codimension $1$ integral scheme. Let $A=\mathcal{O}_{X,Y}$ be the local ring of $\mathcal{O}_X$ at the generic point of $Y$. Since $X$ is normal, this is a DVR.  Let $\mathbf{m}_{X,Y}$ be its maximal ideal. Then $\mathbf{m}_{X,Y}=(t)$. Then, considering $A$ as a subring of $k(X)=k(A)$ and because $X$ is of finite type over $k$, there exists $a\in A$ such that $aD$ is a derivation of $A$, Then  $(aD(A))=(t^n)$, for some   $n\in \mathbb{Z}$. Since $A$ is a DVR, $a=ut^m$, where $u\in A^{\ast}$ and $m\in \mathbb{Z}$. Define $v_D(Y)=n-m$. The divisor
\[
\mathrm{div}(D)=\sum_Y v_D(Y)Y,
\]
where $Y$ runs over all codimension $1$ integral subschemes of $X$ is called the divisorial part of $D$. (As an exercise one can prove that this is a finite sum)
 
Let now $Z\subset X$ be an integral subscheme of $X$ of codimension $\geq 2$. Let $\mathcal{O}_{X,Z}$ be the local ring of $\mathcal{O}_X$ at the generic point of $Z$. Viewing this as a subring of $k(X)$  then, since $X$ is of finite type over $k$,  there exists $a\in \mathcal{O}_{X,Z}$ such that $aD$ is a derivation of $\mathcal{O}_{X,Z}$. Let $\mathbf{m}_{X,Z}$ be the maximal ideal of $\mathcal{O}_{X,Z}$. Then $Z$ is called an isolated part of $D$ if and only if $\mathbf{m}_{X,Z}$ is an embedded prime the ideal $(aD(\mathcal{O}_{X,Z}))$.
\end{definition}

These definitions can be made locally explicit in the case when $X$ is smooth of dimension $n\geq 1$ over $k$. Suppose that this is the case. Let $P\in X$ be a closed point and let $x_1, x_2, \ldots, x_n \in A=\mathcal{O}_{X,P}$ be a regular system of parameters. Let $\hat{A}$ be the completion of $A$ at $P$. Then $\hat{A}\cong k[[x_1,\ldots, x_n]]$. Then
\[
D=\phi \left( f_1\frac{\partial}{\partial x_1}+\cdots + f_n \frac{\partial}{\partial x_n}\right)
\]
such that $\phi \in k(X)$ and $f_i \in A\subset k[[x_1,\ldots, x_n]]$, $i=1,\dots, n$ and the $f_1,\dots, f_n$ do not have a common factor. The divisorial part of $X$ in a neighborhood of $P\in X$ is $\mathrm{div}(\phi)$ and the isolated part is defined by the ideal $(f_1,\ldots, f_n)$. 

\begin{example}
Let $k$ be a field of characteristic $2$. Then
\[
\delta=xy\frac{\partial}{\partial x}+y^2\frac{\partial}{\partial y}.
\]
is a $k$-derivation of $A=k[x,y]$ such that $\delta^2=0$. The fixed part of $\delta$ is defined by the ideal 
\[
(\delta(A))=(\delta(x),\delta(y))=(xy,y^2)=(x,y)\cap (y).
\]
Therefore the divisorial part of $\delta$ is given by the ideal $(y)$ and the isolated part by the ideal $(x,y)$.

\end{example}

Suppose  that $X$ is an algebraic scheme over a field $k$ of characteristic $p>0$ such that $X$ admits a nontrivial $\alpha_p$ or $\mu_p$ action. In the remaining part of this section we will present some results about the structure of the quotient of $X$ by this action.

By the remark after Theorem~\ref{existence of quot},  the geometric quotient $\pi \colon X\rightarrow Y$ exists. Moreover, by Proposition~\ref{explicit des of ap and mp action} the $\alpha_p$ or $\mu_p$ action on $X$ is determined by a vector field $D$ on $X$ such that $D^p=0$ or $D^p=D$, respectively. In addition, from Proposition~\ref{height 1 actions second}, locally the morphism $\pi$ is described as follows. Suppose $X=\mathrm{Spec}A$. Then $Y=\mathrm{Spec} B$, where $B=\{a\in A,\;\;D(a)=0\}$. 

The next proposition provides information about the map $\pi$ and the singularities of $Y$. The reader may find more information in~\cite{AA86},~\cite{Mats23a},~\cite{Mats23b}~\cite{Tz17}.

\begin{proposition}[\cite{Tz17}]\label{singularities of the quotient}
Let $X$ be an integral scheme of finite type over a field $k$ of characteristic $p>0$. Suppose $X$ has a non trivial  $\alpha_p$ or $\mu_p$ action induced by a vector field $D$ on $X$ such that $D^p=0$ or $D^p=D$, respectively. Let $\pi \colon X \rightarrow Y$ be the quotient. Then
\begin{enumerate}
\item The map $\pi$ is a purely inseparable morphism of degree $p$.
\item If $X$ is normal then $Y$ is normal.
\item If $X$ is $S_2$ then $Y$ is $S_2$ as well.
\item If $X$ is smooth then the singularities of $Y$ are exactly the image of the isolated part of the fixed part of $D$  ~\cite{AA86}.
\item If $X$ is normal and $\mathbb{Q}$-Gorenstein, then $Y$ is also $\mathbb{Q}$-Gorenstein. In particular, let $B$ be a divisor in $Y$ and $\tilde{B}$ be the divisorial part of $\pi^{-1}(B)$. Then if $n\tilde{B}$ is Cartier, $pnB$ is Cartier too.
\end{enumerate}
\end{proposition}

%\subsubsection{Singularities with $\mu_p$ actions}
%\begin{remark}
%The singularities of the quotient of a smooth variety with a nontrivial $\mu_p$ action are in a great extend similar to cyclic quotient singularities in characteristic zero~\cite{Hi99}. However, the singularities of even a smooth variety with a non trivial $\alpha_p$-action can be much worse than what can be expected of a quotient of a smooth variety by a finite group in characteristic zero. In particular, in characteristic zero, the quotient of a smooth variety by a finite group has rational singularities~\cite{Burns74}. This is not always true in positive characteristic for quotients by $\alpha_p$ actions, as the next example illustrates. This makes $\alpha_p$ quotients very hard to handle in general.
%\end{remark}

\noindent\textbf{Singularities of the quotient of a smooth surface by a non trivial $\mu_p$ action}.

The singularities of the quotient of a smooth surface with a nontrivial $\mu_p$ action are in a great extend similar to cyclic quotient singularities in characteristic zero, as we will see next.

\begin{definition}
Let $P \in X$ be a normal surface singularity defined over an algebraically closed field $k$. Let $f \colon Y \rightarrow X$ be its minimal resolution.  $P \in X$ is called a toric singularity of type $\frac{1}{p}(1,m)$, $1\leq m\leq p-1$, where $p$ is a prime number,  if the exceptional set of $f$ is a chain of smooth rational curves
\[
\underset{E_1}{\bullet}-\underset{E_2}{\bullet}-\cdots -\underset{E_m}{\bullet}
\]
such that the intersection numbers $E_i^2=-b_i$, are obtained from the continuous fraction decomposition of $p/m$. 
\end{definition}

Let  $P\in X$ be a toric singularity of type $\frac{1}{p}(1,m)$. If the characteristic of the base field is not $p$, then a toric singularity $P\in X$ of type $\frac{1}{p}(1,m)$  is locally analytically isomorphic to the cyclic quotient singularity $\mathbb{A}^2_k/\mu_p$, where the group $\mu_p$ of $p$-roots of unity act on $\mathbb{A}^2_k$ by $\zeta \cdot x =\zeta x$, $\zeta \cdot y =\zeta^m y$, where $\zeta$ is a primitive $p$-root of unity and $x$, $y$ the coordinates of $\mathbb{A}^2_k$. 

Suppose that the characteristic of the ground field $k$ is $p>0$. Let $X$ be a smooth surface with a nontrivial $\mu_p$ action and let $\pi \colon X \rightarrow Y$ be the quotient of $X$ by the $\mu_p$ action. According to Proposition~\ref{explicit des of ap and mp action}, the $\mu_p$ action is induced by a nontrivial vector field $D$ on $X$ of multiplicative type. Moreover, by Proposition~\ref{singularities of the quotient}, the singularities of $Y$ are the image of the isolated fixed points of $D$. 

Let $P\in X$ be an isolated fixed point of $D$ and $Q=\pi(P)\in Y$. Then locally analytically at $P\in X$, $D=x\partial/\partial x +m y \partial /\partial y$, for some $1\leq m\leq p-1$. Let $Q=\pi(P)\in Y$~\cite{RS76}. Then the following holds.

\begin{proposition}\cite{Hi99}\label{class-of-mult-quot}
 The singularity $Q\in Y$ is a toric singularity of type $1/p(1,m)$. In particular, it is a rational singularity.
\end{proposition}

Therefore, according to the previous result,  toric singularities are precisely the singularities of the quotient of a smooth surface by a non trivial $\mu_p$ action. Hence the singularities of the quotient of a smooth surface by a $\mu_p$ action in characteristic $p>0$ are the same as those of the quotient of a smooth surface by a nontrivial $\mathbb{Z}/p\mathbb{Z}$ action in characteristic zero. However, the situation is very different with $\alpha_p$ actions as we will see next.

\medskip

\noindent\textbf{Singularities of the quotient of a smooth surface by a non trivial $\alpha_p$ action.}

In contrast to the case of the singularities of the quotient of a smooth surface by a nontrivial $\mu_p$ action,  the singularities of a smooth surface with a non trivial $\alpha_p$-action can be much worse than what can be expected of a quotient of a smooth variety by a finite group in characteristic zero. In particular, in characteristic zero, the quotient of a smooth variety by a finite group has rational singularities~\cite{Burns74}. This is not always true in positive characteristic for quotients by $\alpha_p$ actions, as the next example illustrates. This makes $\alpha_p$ quotients very hard to handle in general.

\begin{example}[A non rational quotient singularity]\label{a non rational sing}

Let $k$ be a field of characteristic $2$ and let $\delta$ be the $k$-derivation of $k[x,y]$ given by
\[
\delta=x^4\frac{\partial}{\partial x}+y^4\frac{\partial}{\partial y}.
\]
It is easy to see that $\delta^2=0$, i.e., $\delta$ is a additive derivation. The derivation $\delta$ induces an additive vector field $D$ of $\mathbb{A}^2_k$ which by Proposition~\ref{explicit des of ap and mp action} induces a non trivial $\alpha_2$ action on $\mathbb{A}^2_k$.  Let $\pi \colon \mathbb{A}^2_k \rightarrow Y$ be the quotient of $\mathbb{A}^2_k$ by this $\alpha_2$ action. We will show that $Y$ has an elliptic singularity and therefore it is not rational.

The description of the singularities of $Y$ will be done by constructing its minimal resolution $g \colon Z \rightarrow Y$. In fact we will show that there exists a diagram
\[
\xymatrix{
X\ar[r]^f\ar[d]^{\nu} & \mathbb{A}^2_k \ar[d]^{\pi} \\
Z  \ar[r]^g & Y
}
\]
with the following properties.
\begin{enumerate}
\item $X$ is smooth and $f$ is birational. 
\item The morphism $g$ is the minimal resolution of $Y$.
\item The vector field $D$ lifts to an additive vector field $\tilde{D}$ on $X$ without isolated part.
\item $Y$ is the quotient of $X$ with the $\alpha_2$ action on $X$ induced by $\tilde{D}$. 
\end{enumerate}

The first part is to construct $X$. This  will be constructed by repeatedly blowing up the reduced isolated points of the fixed part of $D$ and by showing that by doing so one arrives at a smooth surface $X$ and a vector field $\tilde{D}$ on $X$ lifting $D$ which has no isolated part.

From the definition of $D$, its fixed part is given by the ideal $(x^4,y^4)$ of $k[x,y]$. Therefore the fixed part of $D$ has only isolated part and no divisorial part. 

Let $f_1 \colon X_1 \rightarrow X$ be the blow up of $\mathbb{A}^2_k$ at $(0,0)$ and let $E_1\cong \mathbb{P}^1_k$ be its exceptional curve. Since the point $(0,0)$ is in the fixed part of $D$, $D$ lifts to a vector field $D_1$ on $X_1$ which can be explicitly described as follows. From the definition of blow up, $X_1=U\cup V$, where $U, V\cong \mathbb{A}^2_k$. The restriction of $f_1$ on $U$ is induced from the map $f _1\colon k[x,y]\rightarrow k[x,y]$ given by $f_1(x)=xy$ and $f_1(y)=y$, while the restriction of $f_1$ on $V$ is given by $f_1\colon k[x,y]\rightarrow k[x,y]$ given by $f_1(x)=x$ and $f_1(y)=xy$. 

Restricting $f_1$ on $U$ we get a commutative diagram
\begin{equation}\label{diagram 61}
\xymatrix{
k[x,y] \ar[r]^{f_1} \ar[d]^{\delta} & k[x,y] \ar[d]^{\delta_1}\\
k[x,y] \ar[r]^{f_1 }& k[x,y]
}
\end{equation}
where $\delta_1$ is a $k$-derivation of $k[x,y]$ inducing $D_1$ on $X$. Then $\delta_1=\phi(x,y)\partial/\partial x +\psi(x,y)\partial/ \partial y$. From the commutativity of the previous diagram it follows that $\phi(x,y)=x^4y^3-xy^3$ and $\psi(x,y)=y^4$. Therefore,
\[
\delta_1=y^3\left((x^4-x)\frac{\partial}{\partial x}+y\frac{\partial}{\partial y}\right).
\]
From this description we see that the divisorial part of $\delta_1$ in $U$ is $3E$ given by $y^3=0$ and the isolated part consists of four reduced closed points given by $x^4-x=y=0$. Locally at every isolated fixed point of $\delta_1$, and for appropriate choice of local coordinates $x$ and $y$,
\begin{equation}\label{eq 60}
\delta_1 =y^3\left(x\frac{\partial}{\partial x} +y\frac{\partial}{\partial y}\right).
\end{equation}
By doing an identical calculation on $V$ we see that there is one more reduced isolated fixed point of $D_1$ there of the same form. Therefore, the fixed part of $D_1$ consists of $3E$ and five isolated reduced points such all of them lie on $E$ and  locally at each one of them, $\delta_1$ has has the  form in (\ref{eq 60}) above.

Let now $f_2\colon X \rightarrow X_1$ be the blow up of $X_1$ at the five fixed points of $D_1$. Let also $E_i\cong \mathbb{P}^1_k$, $i=2,3,4,5,6$ be the $f_2$-exceptional curves and call also, for the sake of uniformity of notation, $E_1$ the birational transform of $E_1$ on $X$. Since $f_2$ blows up fixed points of $D_1$, $D_1$ lifts to a vector field $\tilde{D}$ on $X$. A straightforward calculation, similar to the one done for $X_1$ at every isolated fixed point of $D_1$ and using the form (\ref{eq 60}), shows that 
$\tilde{D}$ has no isolated fixed part and its divisorial part is the divisor $3E_1+\sum_{i=2}^6E_i$. Moreover, $E_i^2=-1$, $i=2,3,4,5,6$ and $E_1^2=-6$. 

Let now $\nu \colon X \rightarrow Z$ be the quotient of $X$ by the $\alpha_2$ action on it induced by $\tilde{D}$. Then since $\tilde{D}$ has no isolated part, from Proposition~\ref{singularities of the quotient}, $Z$ is smooth. Moreover there exists a morphism $g \colon Z\rightarrow Y$ producing the diagram (\ref{diagram 61}). Since $f$ is birational and both $\pi$ and $\nu$ are purely inseparable of degree $p$, $g$ is also birational. 

Next we will describe the exceptional set of $g$. The $g$-exceptional curves are the curves $F_i=\nu(E_i)$, $1\leq i \leq 6$. Local straight forward calculations show that $\nu^{\ast}(F_1)=E_1$ and that $\nu^{\ast}(F_i)=2E_i$, for $2\leq i \leq 6$. Therefore, $F_1^2=-3$ and $F_i^2=-2$, for $2\leq i \leq 6$. Moreover by looking at the construction of $X$, we see that $F_1 \cdot F_i=1$, $2\leq i \leq 6$ and $F_i\cdot F_j=0$, $2\leq i,j\leq 6$. Therefore,  the configuration of the $g$-exceptional curves is the following.

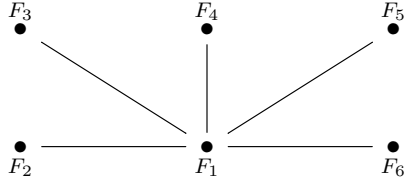
\[
\xymatrix{
    \overset{F_3}{\bullet}  \ar@{-}[drr]        &               &        \overset{F_4}{\bullet}\ar@{-}[d]      &   &\overset{F_5}{\bullet} \ar@{-}[dll]  \\
    \underset{F_2}{\bullet} \ar@{-}[rr]           &              & \underset{F_1}{\bullet}       &        &\underset{F_6}{\bullet} \ar@{-}[ll] \\                     
}
\]
From this it follows that the fundamental cycle $B$ of the singularities of $Y$ is 
\[
B=2E_1+E_2+E_3+E_4+E_5+E_6.
\]
Let $p_a(B)$ be the arithmetic genus of $B$~\cite{Art66}. From this we now get that 
\[
2p_a(B)-2=K_Z\cdot B+B^2=0.
\] 
Therefore, $ p_A(B)=1$ and hence $Y$ has one elliptic singularity. In particular the singularity of $Y$ is not rational.
\end{example}

We will finish this section by presenting an adjunction formula for quotients of surfaces by an $\alpha_p$ or $\mu_p$ action or more generally by a rational vector field of either additive or multiplicative type.

\begin{proposition}[\cite{RS76}]\label{adjunction2}
Let $X$ be a normal $\mathbb{Q}$-Gorenstein surface over a field of characteristic $p>0$. Let $D$ be a nontrivial rational vector field of $X$ such that $D^p=0$ or $D^p=D$. Let $\pi \colon X \rightarrow X^D$ be the quotient of $X$ by $D$ as defined in Definition~\ref{definition 500}. Then
\[
K_X\equiv \pi^{\ast}K_{X^D} +(p-1)\Delta,
\]
where $\Delta$ is the divisorial part of the fixed part of $D$ and $\equiv$ denotes numerical equivalence of divisors.
\end{proposition}
\begin{proof}
Since $X$ is normal and $\mathbb{Q}$-Gorenstein, $K_X$ is $\mathbb{Q}$-Cartier. Then by Proposition~\ref{singularities of the quotient}, $Y$ is normal and $K_Y$ is also $\mathbb{Q}$-Cartier. Therefore $\pi^{\ast}K_Y$ is well defined. Let $U=X-(\mathrm{Sing}(X)\cup Z)$, where $\mathrm{Sing}(X)$ is the singular part of $X$ and $Z$ the isolated part of the fixed part of $D$. Let $V=\pi(U)\subset Y$. Then $U$ and $V$ are open subsets of $X$ and $Y$ of codimension $2$ and moreover, since $D$ has no isolated part in $U$, $V$ is smooth.  Since $X$ is normal, to prove the formula of the proposition it suffices to prove that on $U$ the formula is actually linear equivalence of Cartier divisors. In this case, since $U$ is smooth it follows from~\cite{RS76}.

\end{proof}

By Proposition~\ref{explicit des of ap and mp action} the previous proposition applies directly to the case of quotients by $\alpha_p$ or $\mu_p$ actions.
\begin{corollary}
Let $\pi \colon X\rightarrow Y$ be the quotient of a normal $\mathbb{Q}$-Gorenstein surface $X$ over a field of characteristic $p>0$ by a nontrivial $\alpha_p$ or $\mu_p$ action. Then
\[
K_X\equiv \pi^{\ast}K_Y +(p-1)\Delta,
\]
where $\Delta$ is the divisorial part of the vector field $D$ of $X$ defining the $\alpha_p$ or $\mu_p$-action. 
\end{corollary}

\section{Automorphism group schemes in positive characteristic}\label{automorphisms}

The purpose of this section is to present certain results about the structure of the automorphism group scheme of a variety of general type in positive characteristic and to highlight through explicit examples some important differences between the positive characteristic case and the zero characteristic.

The automorphism group scheme of an algebraic variety is a fundamental invariant of it and an object of great interest and in a great extent determines its deformation and moduli theory. 

A case of particular interest is the automorphism group scheme of smooth varieties of general type and more generally of canonically polarized varieties which are defined next.

\begin{definition}
Let $X$ be a projective variety over an algebraically closed field $k$. Then
\begin{enumerate}
\item $X$ is called of general type if and only if $X$ is smooth and the pluricanonical linear system $|mK_X|$ defines a birational map to its image for $m$ sufficiently large.
\item $X$ is called canonically polarized if and only if it has canonical singularities and $K_X$ is an ample $\mathbb{Q}$-Cartier divisor. These are the canonical models of varieties of general type. In the case of surfaces these are projective surfaces with rational double points and ample canonical class.
\end{enumerate}
\end{definition}
Varieties of general type and their canonical models play a fundamental role in the classification of smooth projective varieties. In characteristic zero canonical models exist in all dimensions but in positive characteristic only in dimension at most $2$~\cite{Bombieri73}. Moreover, the automorphism group scheme of a smooth variety of general type or more generally of a canonically polarized variety is zero dimensional and of finite type over the ground field $k$~\cite{LeMe78}.

Let $X$ be a smooth projective variety of general type over an algebraically closed field of characteristic $p\geq 0$. 

Suppose that $p=0$. Then the following are known.
\begin{theorem}\label{results in char zero}
Let $X$ be a smooth variety of general type over an algebraically closed field of characteristic zero. Then.
\begin{enumerate}
\item $\mathrm{Aut}_{X/k}$ is smooth over $k$ (because every group scheme is smooth in characteristic zero).
\item Suppose that $\dim X=1$. Then  $|\mathrm{Aut}_{X/k}|\leq 42(2p_g(X)-2)$~\cite{Ha77}.
\item Suppose that $\dim X=2$. Then  $|\mathrm{Aut}_{X/k}|\leq 42^2K_X^2$~\cite{Xiao94},~\cite{Xiao95}.
\item Let $n\geq 3$. Then there exists a number $c\in \mathbb{R}$ which depends only on $n$ such that for any smooth variety $X$ of general type and dimension $n$, $|\mathrm{Aut}_{X/k}|\leq cK_X^n$~\cite{HaMcKXu13} .
\end{enumerate}
\end{theorem}
Suppose on the other hand that $p>0$. Then we will show by examples that in dimension $\geq 2$, all of the above fail. In fact one of the fundamental problems in the theory of automorphism groups schemes of varieties of general type is to find explicit bounds for its length and structure results about its group scheme structure and in particular its component $\mathrm{Aut}^0_{X/k}$ containing the identity.

\medskip

\noindent\textbf{Failure of the characteristic zero bounds in positive characteristic}

The next example shows that the bounds that hold true in characteristic zero are no longer true in positive characteristic.

\begin{example}[~\cite{HaMcKXu13}]\label{ex HMK}
Let $p>0$ be a prime number and $q=p^m$, for $m>0$. Let $K=\mathbb{F}_{q^2}$ and let $\lambda \colon K \rightarrow K$ be defined by $\lambda(a)=a^q$. Then $\lambda^2=\lambda$. This map is the positive characteristic analog of complex conjugation. Denote for any $a\in K$, $\bar{a}=\lambda(a)$.

Let $A=(a_{ij}) \in \mathrm{M}_n(K)$ be a square matrix of order $n$ with entries from $K$. Let $\bar{A}=(\bar{a}_{ij})$. In complete analogy with the complex case, a matrix $A\in \mathrm{GL}_n(K)$ is called unitary if and only if $A\bar{A}=I_n$. Let $U_n(K)$ be the set of all unitary matrices of order $n$. $U_n(K)$ is a subgroup of $GL_n(K)$ and therefore acts on $\mathbb{P}^n_K$ by linear automorphisms. 

Let $X \subset \mathbb{P}_K^{n+1}$ be the hypersurface of degree $q+1$ given by the equation
\[
x_0^{q+1}+x_1^{q+1} +\cdots +x_n^{q+1} +x_{n+1}^{q+1}=0.
\]
$X$ is smooth and $K_X=\mathcal{O}_X(q-n-1)$. This is ample if $q-n-1>0$ and therefore it is canonically polarized and in particular of general type. Note that $\mathrm{Aut}_{X/k}$ is smooth and finite since the automorphism group scheme of any smooth hypersurface in $\mathbb{P}^n$ is smooth~\cite{MO67}.

It is not difficult to ckeck that $X$ is invariant under  the action of $U_{n+2}(K)$ on $\mathbb{P}^{n+1}_K$ and therefore $U_{n+2}(K)$ acts faithfully on $X$ and hence it is a subgroup of $\mathrm{Aut}_{X/k}$. Therefore, $|\mathrm{Aut}_{X/k}|\geq |U_{n+2}(K)|$.  Now it is known~\cite{GLS94} that
\begin{equation}\label{order of un}
|U_{n+2}(K)|=\frac{1}{(n+2,q+1)}q^{\binom{n+2}{2}}\prod_{i-2}^{n+2}(q^i-(-1)^i).
\end{equation}
Moreover, $X$ has dimension $n$ and $K_X^n=(q-n-1)^n(q+1)$. In particular, $K_X^n$ is a polynomial of degree $n+1$ in $q$. On the other hand, from the  equation (\ref{order of un}) it follows that the order of $U_{n+2}(K)$ is a polynomial in $q$ of degree
\[
\binom{n+2}{2}+\binom{n+3}{2}-1.
\]
This shows that given an integer $n\geq 2$ there does not exist a number $c\in \mathbb{R}$ which depends only on $n$ such that if $X$ is smooth of general type of dimension $n$, then $|\mathrm{Aut}_{X/k}|\leq cK_X^n$ and therefore  the bound in Theorem~\ref{results in char zero}.2, are false in positive characteristic. 

\end{example}
\begin{problem}
Find an explicit upper bound for the length of the automorphism group scheme of a smooth variety of general type of dimension $n\geq 2$.
\end{problem}

In dimension $1$ things are much simpler. $\mathrm{Aut}_{X/k}$ is smooth and there is the following explicit bound for its length.

\begin{theorem}[\cite{Singh74}]
Let $X$ be a smooth curve of genus $g\geq 2$ defined over an algebraically closed field $k$ of characteristic $p>0$. Then
\[
|\mathrm{Aut}_{X/k}|\leq \frac{4p}{p-1}g^2\left(\frac{2g}{p-1}+1\right).
\]
\end{theorem}
Other than this result I am not aware of any explicit bound for the size of the automorphism group scheme in higher dimensions.

To get an idea what might be expected as a bound, consider  the case $n=2$ in Example~\ref{ex HMK}. Then $X$ is surface of degree $q+1$ in $\mathbb{P}^3_K$ and $K_X^2=(q-3)^2(q+1)$, which is a polynomial of degree $3$ in $q$. On the other hand,
\[
|\mathrm{Aut}_{X/k}|\geq |U_4(K)|=\frac{1}{(4,q+1)}q^6(q^2-1)(q^3+1)(q^4-1),
\]
which is a polynomial of degree $15$ in $q$ and much greater than the bound $42^2K_X^2$ that holds in characteristic zero. In fact this shows that in the case of surfaces the best possible analog of the Xiao or  Hacon-McKerna-Xu bound in Theorem~\ref{results in char zero} would be something like, for any $n>0$ there exists a $c\in \mathbb{R}$ such that $|\mathrm{Aut}_{X/k}|\leq c(K^2_X)^5$.

\medskip

\noindent\textbf{Non reduced automorphism group schemes in positive characteristic.}

Let $X$ be a projective variety of dimension at least $2$ over an algebraically closed field $k$ of characteristic $p>0$. Unlike the case of characteristic zero, when all group schemes are smooth, the automorphism group scheme $\mathrm{Aut}_{X/k}$ of $X$ may not be smooth over $k$ anymore. This fact complicates the moduli theory of $X$ compared to the characteristic zero case. 

If $X$ is a smooth curve of genus $g\geq 0$, then $\mathrm{Aut}_{X/k}$ is smooth over $k$. Indeed. If $g=0$ then $\mathrm{Aut}_{X/k}=PGL(2,k)$, which is smooth. If $g=1$, then $X$ is an elliptic curve and also $\mathrm{Aut}_{X/k}$ is smooth~\cite{Silv86}. Suppose that $g\geq 0$. Then $\mathrm{Aut}_{X/k}$ is finite and its tangent space at the identity is $\mathrm{Hom}_X(\Omega_{X/k},\mathcal{O}_X)=H^0(\omega_X^{-1})=0$, because $\omega_X$ is ample. Therefore $\mathrm{Aut}_{X/k}$ is smooth again.

Examples of smooth surfaces of general type with non smooth automorphism group scheme have been found in~\cite{La83},~\cite{Li08},~\cite{SB96},~\cite{Russell84}. Examples of del Pezzo surfaces with rational double points with non smooth automorphism group scheme have been found in~\cite{MS20},~\cite{M21}.

Next we will present two explicit examples, one of a canonically polarized surface and one of a smooth surface of general type with non smooth automorphism group scheme. The two examples will be obtained by two different methods and they are indicative how such examples can be obtained. Another example will be obtained in Section~\ref{failure of KV}.

A key observation is the following. 

\begin{proposition}
Let  $X$ be a canonically polarized variety over a field $k$ of characteristic $p>0$. Then $\mathrm{Aut}_{X/k}$ is not smooth if and only if $X$ has a non trivial global vector field $D$ such that either $D^p=0$ or $D^p=D$. Equivalently if and only if there exists a non trivial action of $\alpha_p$ or $\mu_p$ on $X$.
\end{proposition}
\begin{proof}
It is well known~\cite{LeMe78} that  $\mathrm{Aut}_{X/k}$ is zero dimensional and of finite type over $k$. Therefore $\mathrm{Aut}_{X/k}$ is not smooth if and only if its tangent space at the identity is not trivial. The tangent space at the identity is $\mathrm{Hom}_X(\Omega_{X/k},\mathcal{O}_X)=\mathrm{Der}_k(X)$. Hence $\mathrm{Aut}_{X/k}$ is not smooth if and only if $X$ admits a nontrivial global vector field $D$. Then by~\cite{RS76} there exists a non trivial vector field $D$ on $X$ such that $D^p=D$ or $D^p=0$. Therefore by Proposition~\ref{explicit des of ap and mp action} there exists a non trivial $\mu_p$ or $\alpha_p$ action on $X$.
\end{proof}

%Next we will give two examples of surfaces of general type with non reduced automorphism group scheme. The examples will exhibit two methods that can be used to obtain such surfaces. In the first example a surface will be constructed having an explicit non trivial vector field and in the second the surface will be the quotient of a rational surface by a non trivial rational vector field.  

\begin{example}
We will next construct a canonically polarized surface $X$  such that Frobenius Kernel of its automorphism scheme is $\mu_p$, for any prime $p\geq 2$. In particular, $\mathrm{Aut}_{X/k}$ is not smooth. This example appeared in the first version of~\cite{ScTz23} that can be found on the arxiv.

Let 
\[
D=x_1\frac{\partial}{\partial x_1}-x_2\frac{\partial}{\partial x_2}
\]
be a graded derivation of the polynomial ring $k[x_1,x_2,x_3,x_4]$. Then it is not hard to see that $D^p=D$. 

Let us next examine what the fixed part of $D$ is. This can be seen by working locally in each open set  of the standard affine cover of $\mathbb{P}^3_k$. Consider first $U_1\cong \mathbb{A}^3_k$ to be the affine open subset of $\mathbb{P}^3_k$ with coordinates $y_1=x_2/x_1$, $y_2=x_3/x_1$ and $y_3=x_4/x_1$. Then $D(y_1)=D(x_2/x_1)=-2x_2/x_1=-2y_1$. Similarly $D(y_2)=-y_2$ and $D(y_3)=-y_3$. Therefore  on $U_1$, $D$ is given by
\[
D=-2y_1\frac{\partial}{\partial y_1}-y_2\frac{\partial}{\partial y_2}-y_3\frac{\partial}{\partial y_3}.
\]
Hence the fixed part of $D$ on $U_1$ is given by the ideal $(2y_1,y_2,y_3)$. If $p\not=2$, this is the point $[1,0,0,0]$ of $\mathbb{P}^3_k$. If on the other hand $p=2$, then this is the line $y_2=y_3=0$ which is the restriction of the line $x_2=x_3=0$ on $U_1$. Similar calculations in the remaining affine open sets $U_i$ on $\mathbb{P}^3_k$ given by $x_i\not=0$, $i=2,3,4$ show that  if $p\not=2$, the fixed part of $D$ is the union of the line given by $x_1=x_2=0$ and the points $[1,0,0,0]$ and $[0,1,0,0]$. In particular it is one-dimensional. If on the other hand $p=2$, the fixed part of $D$ is the union of the lines $L_1$ given by  $x_1=x_2=0$ and $L_2$ given by $x_3=x_4=0$.

let $X\subset \mathbb{P}^3_k$ be the hypersurface given by the equation
\begin{gather}\label{examples-eq-6}
f(x_1,x_2,x_3,x_4)=x_1x_2(x_3^{2p-1}+x_4^{2p-1})+x_1^{2p}x_3+x_2^{2p}x_4+x_3^{2p+1}+x_4^{2p+1}=0
\end{gather}
Then $D(f)=0$ and therefore $D$ induces a nontrivial global vector field on $X$ which we still call $D$. If $p\not=2$, since the line $L_1$ is not contained in $X$, the fixed part of $D$ on $X$ is zero dimensional and hence has only isolated part. If on the other hand $p=2$, then since the line $L_2$ given by $x_3=x_4=0$ is contained in $X$, $L_1$ is the fixed part of the restriction of $D$ on $X$. In particular if $p=2$, the fixed part of $D$ has no isolated part.

Let $\lieg=\mathfrak{Lie}(\mathrm{Aut}_{X/k})=Der_k(X)$ be the Lie algebra of vector fields on $X$ and $\mathrm{Aut}_{X/k}[F]$ the Frobenius Kernel of the automorphism scheme $\mathrm{Aut}_{X/k}$ of $X$. Then

\begin{proposition}\label{example 2}
$X$ is a canonically polarized surface and
\[
\lieg \cong \lieg\mathfrak{l}_1(k) \quadand \mathrm{Aut}_{X/k}[F]\cong \mu_p.
\]
Moreover, $c_1^2=(2p-3)^2(2p+1)$ and $c_2=8p^3-4p^2-14p+3$. In particular, $\mathrm{Aut}_{X/k}$ is not smooth.
\end{proposition}

\begin{proof}

We will first show that  $X$ is normal and that its singularities are rational double points of type $A_n$, for some $n$.

In order to do so we will study the equation defining $X$ on each one of the open affine sets of the standard affine cover $U_i$  of $\mathbb{P}^3_k$, where $U_i=\{x_i\not= 0\}$, $i=1,2,3,4$.  We will only work in detail the case of $U_1$. The rest are identical and are omitted. 

Let $y_2=x_2/x_1$, $y_3=x_3/x_1$ and $y_4=x_4/x_1$ be the coordinates of $U_1\cong \mathbb{A}_k^3$. The singular points of $X$ in $U_1$ are the solutions of the system of equations

%Then a straightforward calculation shows that in $U_1$,
%\[
%D=-2u_2\frac{\partial}{\partial u_2}-u_3\frac{\partial}{\partial u_3}-u_4\frac{\partial}{\partial u_4}.
%\]
%In particular we see that in $U_1$, the fixed part of $D$ is exactly one point, if $p\not=2$, and the line $u_3=u_4=0$, if $p=2$ . The singular points of $X$ in $U_1$ are the solutions of the system of equations
\begin{gather}\label{ex2-1}
f=y_2(y_3^{2p-1}+y_4^{2p-1})+y_3+y_2^{2p}y_4+y_3^{2p+1}+y_4^{2p+1}=0\\\label{ex2-2}
\frac{\partial f}{\partial y_2}=y_3^{2p-1}+y_4^{2p-1}=0\\\label{ex2-3}
\frac{\partial f}{\partial y_y}=-y_2y_3^{2p-2}+1+y_3^{2p}=0\\\label{ex2-4}
\frac{\partial f}{\partial y_4}=-y_2y_4^{2p-2}+y_2^{2p}+y_4^{2p}
\end{gather}
It is not hard to see that the equations (\ref{ex2-1})-(\ref{ex2-4}) have finitely many solutions and therefore $X$ is normal. 

By~\cite[Lemma 14.3 and Proposition 14.4]{ScTz23},  the singularities of $X$ are not of type $A_n$,  if the following equations are satisfied
\begin{gather}\label{ex2-5}
\left(\frac{\partial^2 f}{\partial y_2 \partial y_3}\right)^2-\frac{\partial^2 f}{\partial y_2^2}\frac{\partial^2 f}{ \partial y_3^2}=y_3^{4p-4}=0\\\label{ex2-6}
\left(\frac{\partial^2 f}{\partial y_2 \partial y_4}\right)^2 -\frac{\partial^2 f}{\partial y_2^2}\frac{\partial^2 f}{ \partial y_4^2}=y_4^{4p-4}=0\\\label{ex2-7}
\left(\frac{\partial^2 f}{\partial y_3 \partial y_4}\right)^2 -\frac{\partial^2 f}{\partial y_3^2}\frac{\partial^2 f}{ \partial y_4^2}=4y_2^2y_3^{2p-3}y_4^{2p-3}=0
\end{gather}
The equations (\ref{ex2-5}) and (\ref{ex2-3}) are incompatible. Hence every singularity of $X$ is of type $A_n$. Therefore $X$ is a normal surface whose singularities are  rational double points and therefore it is a canonically polarized surface.

By the adjunction formula in $\mathbb{P}^3$ it follows that $\omega_X=\mathcal{O}_X(2p-3)$, and hence $\omega_X$ is ample, for all $p>0$. Therefore $X$ is a canonically polarized surface. Moreover, $c_1^2=K_X^2=(2p-3)^2(2p+1)$. To find $c_2$ we go through the minimal resolution $f \colon Y\rightarrow X$. Since $X$ has rational double points, $K_Y=f^{\ast}K_X$ and $\chi(\mathcal{O}_X)=\chi(\mathcal{O}_Y)$. Then by Noether's formula on $Y$,
\[
c_2=c_2(X)=12\chi(\mathcal{O}_Y)-K_Y^2=12\chi(\mathcal{O}_X)-c_1^2.
\]
A straightforward calculation of $\chi(\mathcal{O}_X)$ and the above formula give the formula for $c_2$ stated in the proposition.

Let $\Delta$ be the divisorial part of the fixed part of $D$. We will next show that $H^0(\mathcal{O}_X(\Delta))=k$. This is obvious if $p\not=2$, since in this case $\Delta=0$. Suppose that $p=2$. Then $\Delta=L_2$ is a line in $\mathbb{P}^3$. Then by Serre duality, $H^0(\mathcal{O}_X(\Delta))=H^2(\omega_X(-\Delta))$. Then from the exact sequence
\[
0 \rightarrow \omega_X(-\Delta)\rightarrow \omega_X\rightarrow \omega_X\otimes \mathcal{O}_{\Delta}\rightarrow 0,
\]
we get by taking cohomology the exact sequence
\[
H^1(\omega_X\otimes \mathcal{O}_{\Delta})
\rightarrow H^2(\omega_X(-\Delta)) \rightarrow H^2(\omega_X)\rightarrow 0
\]
Now since $\omega_X$ is ample and $\Delta\cong \mathbb{P}^1$, it follows that $\omega_X\otimes \mathcal{O}_{\Delta}=\mathcal{O}_{\mathbb{P}^1}(d)$, $d\geq 0$. Therefore $H^1(\omega_X\otimes \mathcal{O}_{\Delta})=0$ and hence 
\[
 H^2(\omega_X(-\Delta)) = H^2(\omega_X)=H^0(\mathcal{O}_X)=k.
 \]
 Then by~\cite[Lemma 14.5]{ScTz23}, there exists an exact sequence
 \[
 0 \rightarrow \mathcal{O}_X(\Delta)\rightarrow \mathcal{H}om_X(\Omega_X,\mathcal{O}_X)  \rightarrow \omega_X^{-1}(-\Delta) \rightarrow N \rightarrow 0,
\]
where $N$ is a coherent sheaf with zero dimensional support. By taking cohomology in the previous sequence it follows, since $\omega_X$ is ample and $\Delta$ effective, that 
 \[
Der_k(X)= \mathrm{Hom}_X(\Omega_{X/k},\mathcal{O}_X)=H^0(\mathcal{H}om_X(\Omega_{X/k},\mathcal{O}_X))=k. 
 \]
Therefore the Lie algebra $\lieg$ of $\mathrm{Aut}_{X/k}$ is one dimensional generated by the vector field $D$ such that $D^p=D$. Hence $\lieg \cong \lieg\mathfrak{l}_1(k)$. Therefore by the Demazure Gabriel correspondence, since  $\mathrm{Aut}_{X/k}[F]$ is a height $1$ group scheme with restricted Lie algebra $\mathfrak{g}=\mathfrak{gl}_1(k)$, it  follows that $\mathrm{Aut}_{X/k}[F]\cong \mu_p$.

\end{proof}

\end{example}

The following  example is due to  C. Liedtke~\cite{Li08}. 

\begin{example}[\cite{Li08}]
Let $S=\mathbb{P}^1_k\times \mathbb{P}^1_k$, where $k$ is an algebraically closed field of characteristic $2$. Let $x_0,x_1$ and $y_0, y_1$ be the homogeneous coordinates in the first and second copy of $\mathbb{P}^1_k$ in $S$, respectively. Let $K=k(S)$ be the field of rational functions of $S$. Then $K$ may be naturally identified with $k(x,y)$, where $x=x_1/x_0$ and $y=y_1/y_0$. Consider now the $k$-derivation $D$ of $K$ defined as follows.
\[
D=\prod_{i=1}^n\frac{(x-a_i)^2}{(x-b_i)^2}\frac{\partial}{\partial x} +\prod_{i=1}^n\frac{(y-c_i)^2}{(y-d_i)^2}\frac{\partial}{\partial y} ,
\]
where $a_i,b_i.c_i,d_i \in k$, $1\leq i \leq n$ are all pairwise disjoint.  Then $D$ defines a rational vector field on $S$ such that $D^2=0$. 

Let $U_i$, $V_i$ be the affine open subsets of the two copies of $\mathbb{P}^1_k$ in $S$ given by $x_i\not= 0$ and 
$y_i\not= 0$, respectively, $i=1,2$. Then $U_i\times V_j$, $1\leq i,j \leq 2$ is an affine open cover of $S$ and $D$ in the form described above  is defined in the affine open set $U_0\times V_0$. In $U_0\times V_0$, the divisorial part of $D$ is the divisor of the rational function 
\[
\phi=\prod_{i=1}^n\frac{1}{(x-b_i)^2(y-d_i)^2}
\]
Moreover the isolated part of the fixed part of $D$ in $U_1\times V_1$ are the $2n^2$ points $(a_i,c_j)$ and $(b_i,d_i)$, $1\leq i,j\leq n$. Describing $D$ in the other affine open sets $U_i\times V_j$, we see that the points $([0,1],Q_j)$, $(P_i,[0,1])$ and $([0,1],[0,1])$ are also isolated fixed points of $D$, where $Q_j=[1,c_j]$ and $P_i=[1,a_i]$, $1\leq i,j\leq n$. Therefore, the total number of isolated points of the fixed part of $D$ is $2n^2+2n+1$ and no new divisorial part exists in the other open sets. Hence globally on $S$, the divisorial part of $D$ is the divisor 
\begin{equation}\label{div part}
\Delta=-2nl_1-2nl_2,
\end{equation}
where $l_1$, $l_2$ are two fibers of the projections $p_i\colon S\rightarrow \mathbb{P}^1_l$, $i=1,2$,  to the first and second factor respectively, and the isolated part consists of $2n^2+2n+1$ reduced points.

Let us now describe $D$ locally at an isolated point of its fixed part. From the description above it follows that locally analytically at an isolated point of the form $(a_i,c_j)$,
\begin{equation}\label{eq70}
D=\phi\left(u(x)x^2\frac{\partial}{\partial x} +w(y)y^2\frac{\partial}{\partial y}\right).
\end{equation}
On the other hand at an isolated point of the form $(b_i,d_j)$,
\begin{equation}\label{eq71}
D=\phi\left(u(x)y^2\frac{\partial}{\partial x} +w(y)x^2\frac{\partial}{\partial y}\right),
\end{equation}
where $u(x), w(y)\in k[[x,y]]$ are units in $x$ and $y$.

Let now $\pi \colon S \rightarrow Y$ be the quotient of $S$ by the rational vector field $D$ of $S$. This is a purely inseparable morphism of degree $p$. Also, according to Proposition~\ref{singularities of the quotient}, $Y$ is normal and its singular points are the image by $\pi$ of the isolated fixed points of $D$. Moreover, since locally at its isolated fixed points $D$ is given by (\ref{eq70}) or (\ref{eq71}), according to~\cite[Proposition 3.2]{Li08}, the singular locus of $Y$ consists of $2n^2+2n+1$  rational double points of type $D_4$. One could also see this directly by doing an explicit calculation as in the Example~\ref{a non rational sing}. Let also $g\colon X\rightarrow Y$ be the minimal resolution of $X$. Then.

\begin{proposition}
$Y$ is a canonically polarized surface and $X$ a smooth surface of general type such that 
\[
c_1(X)^2=c_1(Y)^2=4(n-1)^2 \quadand c_2(Y)=4  \quadand c_2(X)=8(n^2+n+1).
\]
Moreover,
\begin{enumerate}
\item The automorphism group schemes $\mathrm{Aut}_{X/k}$ and $\mathrm{Aut}_{Y/k}$ are not smooth.
\item There exists fibrations $Y\rightarrow \mathbb{P}^1$ and $X\rightarrow \mathbb{P}^1$ such that every fiber of both is a singular rational curve of arithmetic genus $n$.
\end{enumerate}
\end{proposition}

\begin{proof}
By the adjunction formula for purely inseparable morphisms in Proposition~\ref{adjunction2}, and the equation (\ref{div part}), it follows that
\[
K_S=\pi^{\ast}K_Y+\Delta=\pi^{\ast}K_Y-2nl_1-2nl_2.
\]
Considering that $K_S=-2l_1-2l_2$, it follows that 
\begin{equation}\label{eq 74}
\pi^{\ast}K_Y=2(n-1)l_1+2(n-1)l_2.
\end{equation}
In particular this is an ample divisor on $S$. Therefore, since $\pi$ is finite, $K_Y$ is an ample Cartier divisor on $Y$. Hence, as explained earlier, $Y$ has only rational double points and therefore is a canonically polarized surface. 

Let $g\colon X\rightarrow Y$ be the minimal resolution. Since $Y$ has rational double point singularities, $K_X=g^{\ast}K_Y$. Therefore $K_X^2=K_Y^2$. Then from this and  (\ref{eq 74}) we get the statement about the first Chern numbers.

Since $\pi$ is purely inseparable, $c_2(Y)=c_2(X)=4$. Moreover, since $Y$ has exactly $2n^2+2n+1$ singular points, all of which are of type $D_4$, $g$ has exactly $4(2n^2+2n+1)$ exceptional curves. Hence
\[
c_2(X)=c_2(Y)+8n^2+8n+4=8(n^2+n+1).
\]

We will next show the existence of the fibration $Y\rightarrow \mathbb{P}^1$ such that every fiber is a singular rational curve of arithmetic genus $n$. Since $\pi$ is of degree $\pi$, the geometric Frobenius morphism $F^{(p)}\colon S\rightarrow S^{(p)}\cong S$ factors through $\pi$. Hence there exists a commutative diagram
\begin{equation}\label{eq 80}
\xymatrix{
          &       Y  \ar[dr]^{\nu} &     \\
 S \ar[ur]^{\pi} \ar[rr]^{F^{(p)}} \ar[d]^{p_1}&             &        S \ar[d]^{p_1}     \\
     \mathbb{P}^1\ar[rr]^{F^{(p)}} && \mathbb{P}^1              
}
\end{equation}
Let $\phi \colon Y \rightarrow \mathbb{P}^1$ be the composition $p_1\circ \nu$. Let now $Y_t=\phi^{-1}(t)$, for any $t\in \mathbb{P}^1$. From the commutativity of the above diagram it follows that $\pi^{\ast}(Y_t)=2F_t$, where $F_t=p_1^{-1}(t)$. Then from the equation (\ref{eq 74}) it follows that $K_Y\cdot Y_t=2(n-1)$. Then for general $t\in \mathbb{P}^1$, $Y_t$ is in. the smooth part of $Y$  and hence from the adjunction formula, $p_a(Y_t)=n$, as claimed. Moreover, the restriction of $\pi$ on $F_t$ gives a map $ F_t\rightarrow Y_t$ which in fact must is birational. Considering that $F_t\cong \mathbb{P}^1$,  $Y_t$ is rational as claimed.

We will next show that $\mathrm{Aut}_{Y/k}$ is not smooth. Since $Y$ has rational double points and $K_Y$ is ample, $\mathrm{Aut}_{Y/k}$ is finite over $k$~\cite{LeMe78}. To prove that it is not smooth it suffices to show that its tangent space at the identity, which is $\mathrm{Hom}_Y(\Omega_Y,\mathcal{O}_Y)$,  is not trivial. 

Let $U\subset Y$ be the smooth part of $Y$. This is a codimension $2$ open set. Let $V=\nu(U) \subset S^{(p)}\cong S$, where $\nu$ is as in the equation (\ref{eq 80}) above,  Then, since both $U$ and $V$ are smooth,  dualizing the exact sequence in Proposition~\ref{general-structure-theory}.4,  we get the injection
\[
0 \rightarrow \mathcal{F} \rightarrow \mathcal{H}om_U(\Omega_{U/k},\mathcal{O}_U)
\]
Moreover, from Corollary~\ref{adjunction 2}, $\mathcal{F}=\omega_U\otimes \nu^{\ast}\omega_V$. Taking sections and since $U$ has codimension $2$ in $Y$, and $\omega_Y$, $\nu^{\ast}\omega_S$ and $\mathcal{H}om_Y(\Omega_{Y/k},\mathcal{O}_Y)$ are reflexive, we get the inclusion of sections
\begin{equation}\label{eq90}
0 \rightarrow H^0(\omega_Y\otimes \nu^{\ast}\omega_S^{-1})\rightarrow \mathrm{Hom}_Y(\Omega_{Y/k},\mathcal{O}_Y).
\end{equation}
We will show that $H^0(\omega_Y\otimes \nu^{\ast}\omega_S^{-1})\not=0$ and therefore $\mathrm{Hom}_Y(\Omega_{Y/k},\mathcal{O}_Y)\not= 0$. Hence $\mathrm{Aut}_{Y/k}$ is not smooth.

Note that $h^0(\nu^{\ast}\omega_S)=0$ because if this was not true then $h^0(\pi^{\ast}(\nu^{\ast}\omega_S))=h^0((F^{(p)})^{\ast}\omega_S)=h^0(\omega^{\otimes 2}_S)\not= 0$, which is not true. Therefore by Serre duality, $h^0(\omega_Y\otimes \nu^{\ast}\omega_S^{-1})=0$. Since $Y$ has rational double points and $\omega_Y$, $\nu^{\ast}\omega_S^{-1}$ are invertible, Rieman-Roch holds for $\omega_Y\otimes \nu^{\ast}\omega_S^{-1}$. Hence, since $h^0(\omega_Y\otimes \nu^{\ast}\omega_S^{-1})=0$, it follows that
\begin{gather}\label{eq91}
h^0(\omega_Y\otimes \nu^{\ast}\omega_S^{-1}) \geq \chi(\omega_Y\otimes \nu^{\ast}\omega_S^{-1}) =\\\nonumber
\chi(\mathcal{O}_Y) +\frac{1}{2}\left( (K_Y-\nu^{\ast}K_S)^2-K_Y\cdot (K_Y-\nu^{\ast}K_S)\right)=\\\nonumber 
\chi(\mathcal{O}_X) -\frac{1}{2}\nu^{\ast}K_S(K_Y-\nu^{\ast}K_S)=\chi(\mathcal{O}_X)-\frac{1}{4}(\pi^{\ast}\nu^{\ast}K_S)\pi^{\ast}(K_Y-\nu^{\ast}K_S)=\\\nonumber
\frac{1}{12}(c_1^2(X)+c_2(X))-\frac{1}{4}(2K_S)((2n-2)l_1+(2n-2)l_2+4l_1+4l_2)= \\\nonumber
\frac{1}{12}(4(n-1)^2+8n^2+8n+8)+4(n+1)=n^2+4n+5>0.\nonumber
\end{gather}
Therefore $\mathrm{Aut}_{Y/k}$ is not smooth. Morever, the above calculation shows that not only $\mathrm{Aut}_{Y/k}$ is not smooth but its length goes to infinity as $n$ goes to infinity, or equivalently as $K_Y^2$ goes to infinity.  

Next we will show that $\mathrm{Aut}_{X/k}$ is also not smooth. By~\cite{Hi99}, there exists a commutative diagram
\[
\xymatrix{
Z \ar[d]_{\tau} \ar[r]^f & S \ar[d]^{\pi} \\
X\ar[r]^g  & Y
}
\]
with the following properties.
\begin{enumerate}
\item $Z$ is smooth and $f$ birational.
\item The rational vector field $D$ of $X$ lifts to a rational vector field $\delta$ on $Z$ whose fixed part has no isolated part.
\item $X$ is the quotient of $Z$ by $\delta$.
\end{enumerate}
Note that one could obtain this diagram directly by repeating the process in Example~\ref{a non rational sing}. Locally at its isolated points, $D$ is given by either (\ref{eq70}) or (\ref{eq71}). Then blowing up repeatedly the isolated points of $D$ and making use of these two possible forms as in Example~\ref{a non rational sing} one gets $Z$ and the fact that $D$ lifts to a rational vector field $\delta$  on $Z$ without isolated part. Then the quotient $X$ of $Z$ by $\delta$ is smooth and one gets the claimed diagram. 

As in the case of $Y$, the geometric Frobenius $F^{(p)}\colon Z\rightarrow Z^{(p)}$ factors throught $\tau$. This means that there exists a map $\sigma \colon X\rightarrow Z^{(2)}$ such that $\sigma \tau=F^{(p)}$. As in the case of $Y$ there exists an inclusion
\[
0 \rightarrow H^0(\omega_X\otimes \nu^{\ast}\omega_{Z^{(p)}}^{-1})\rightarrow \mathrm{Hom}_X(\Omega_{X/k},\mathcal{O}_X).
\]
An identical but lengthier calculation as in $Y$ shows that $H^0(\omega_X\otimes \nu^{\ast}\omega_{Z^{(p)}}^{-1})\not= 0$ and hence $\mathrm{Aut}_{X/k}$ is not reduced. Note that it would have been possible to work directly in $X$ and not do the case of $Y$ and then descend to $Y$ through the map $g$ but the calculations on $Y$ are simpler and makes the exposition clearer. The reader may repeat the calculations identically for $X$.

\end{proof}
\end{example}

\begin{remark}
The Proposition does not only produce an example in characteristic $2$ of a surface of general type $X$ such that $\mathrm{Aut}_{X/k}$ is not reduced but in fact a series of examples  $X_n$ (depending on $n$) such that the  automorphism group schemes $\mathrm{Aut}_{X_n/k}$ is not smooth and moreover, the order $|\mathrm{Aut}_{X_n/k}|$ goes to infinity as $n$ does. This follows from the equations (\ref{eq90}) and (\ref{eq91}) in the proof. 

\end{remark}

\begin{exercise}
Show that $\mathrm{Aut}_{X/k}[F]\cong \alpha_p^{\oplus m}$, for some $m>0$. Moreover, find $m$ explicitly as a function of $n$.  (By Theorem 12.1~\cite{ScTz23}, it suffices to show that there does not exist a non trivial $\mu_p$ action on $X$, or $Y$. Then use the method of Proposition 14.2~\cite{ScTz23}).
\end{exercise}
Having established the existence of smooth varieties of general type, or more generally of canonically polarized varieties with finite and non reduced automorphism group scheme, a natural problem is the following.

\begin{problem}
What is the structure as a group scheme of $\mathrm{Aut}_{X/k}$, or its component containing the identity $\mathrm{Aut}^0_{X/k}$. Moreover, What is the geometry of varieties with non reduced automorphism group scheme?
\end{problem}
In dimension at least $3$ these are, to my knowledge, completely open problems. However, in dimension $2$ the following regarding the structure of the Frobenius Kernel of the automorphism group scheme is known.

\begin{theorem}[\cite{ScTz23}]
Let $X$ be a projective normal surface over a field $k$ of characteristic $p>0$ such that $h^0(\omega_X^{-1})=0$.  Then the Frobenius kernel $\mathrm{Aut}_{X/k}[F]$ of $\mathrm{Aut}_{X/k}$ is isomorphic to one of the following  group schemes:
$$
\mathrm{SL}_2[F] \quadand \alpha_p^{\oplus n} \quadand \alpha_p^{\oplus n}\rtimes\mu_p,
$$
for some integer $n\geq 0$. In particular, if $X$ is a canonically polarized surface, then 
\[
n \leq 
\begin{cases}
\frac{1}{144}(73c_1^2+c_2)^2-1	& \text{if $c_1^2\geq 2$;}\\
\frac{1}{144}(121c_1^2+c_2)^2-1	& \text{if $c_1^2=1$.}
\end{cases}
\]
\end{theorem}
The above result includes more cases other than surfaces of general type, for example properly elliptic surfaces.

Finally, regarding the geometry  of canonically polarized surfaces with non reduced automorphism group scheme the following is known. More detailed results in specific cases can be found in~\cite{Tz17b},~\cite{Tz22},~\cite{Tz23}.

\begin{theorem}[\cite{Tz22}]
Let $X$ be a canonically polarized surface over an algebraically closed field of characteristic $p>0$. Suppose that $\mathrm{Aut}_{X/k}$ is not reduced and  
\[
p > \mathrm{max} \{8(K_X^2)^3+12(K_X^2)^2+3, 4508K_X^2+3\}.
\]
Then $X$ is unirational and $\pi_1(X)=\{1\}$.
\end{theorem}

\section{Failure of the Kodaira vanishing theorem in positive characteristic.}\label{failure of KV}

One of the fundamental tools in the classification theory of smooth projective varieties defined over an algebraically closed field of characteristic zero is the Kodaira vanishing theorem and its variations.
\begin{theorem}[Kodaira Vanishing Theorem]
Let $X$ be a smooth projective variety defined over an algebraically closed field of characteristic zero. Let $L$ be an ample invertible sheaf on $X$. Then
\begin{enumerate}
\item $H^i(X,L^{-1})=0$,  for all $ i< \dim X$.
\item $H^i(X, L\otimes \omega_X)=0$, for all $ i\geq 1$.
\end{enumerate}
\end{theorem}
An excellent treaty of the Kodaira Vanishing Theorem and other vanishing theorems can be found in~\cite{EsVie92}.

In characteristic $p>0$ the Kodaira Vanishing Theorem as stated above does not hold anymore. This is a major obstacle to extending the classification theory of characteristic zero to positive characteristic.

Counterexamples to the Kodaira Vanishing Theorem in positive characteristic have been obtained initially by D. Mumford~\cite{Mu67} and M. Raynaud~\cite{Ray78} and later by many others like S. Mukai~\cite{Mukai11} and X. Zheng~\cite{Xudo17}.

The main purpose of this section is to discuss the failure of Kodaira Vanishing in positive characteristic. In particular, to present conditions under which it holds, or fails, and give an explicit counterexample. 

\subsection{Conditions for the Kodaira Vanishing to hold in characteristic $p>0$}
The reader may find a very detailed exposition of results concerning conditions under which the Kodaira vanishing holds in positive characteristic in~\cite{EsVie92}.

One of the most general results in this direction is a criterion obtained by P. Deligne and L. Illusie~\cite{DeIll87}. Before we state the theorem, we define what it means for a variety to lift to characteristic zero. 

\begin{definition}
Let $X$ be a scheme of finite type over a perfect field $k$ of characteristic $p>0$. Let $(R,\mathbf{m}_R)$ be a local integral domain such that the function field $k(R)$ of $R$ has characteristic zero and $R/\mathbf{m}_R=k$,  (For example $R=\mathbb{Z}_p$, $k=\mathbb{Z}_p/(p)=\mathbb{F}_p$.)

We say that $X$ lifts to characteristic zero over $R$ if there exists a flat morphism $f\colon \mathcal{X} \rightarrow \mathrm{Spec} R$, such that the special fiber is isomorphic to $X$.

\end{definition}

There are many ways that a scheme $X$  over a perfect field  of characteristic $p>0$ may lift to characteristic zero. These depend on the choice of the ring $R$ as in the previous definition. One particular such ring with especially nice properties is the Witt ring or the ring of Witt vectors~\cite{Serre79}.  

\begin{proposition}[\cite{Serre79},~\cite{Pink05}]
Let $k$ be a perfect field of characteristic $p>0$. Then there exists a complete discrete valuation ring $W(k)$ with the following properties.
\begin{enumerate}
\item The field of fractions of $W(k)$ has characteristic zero.
\item The maximal ideal of $W(k)$ is generated by $p$  and its residue field is isomorphic to $k$.
\item The Frobenius map $F\colon k\rightarrow k$ lifts to a ring homomorphism of $W(k)$.
\item Every complete discrete valuation ring with function field of characteristic zero and residue field $k$ contains $W(k)$ as a subring.
\end{enumerate}
\end{proposition}
Let $W(k)$ be the Witt ring where $k$ has characteristic $p>0$. The Artin ring ring $W_n(k)=W(k)/(p^n)$ is called the ring of $n$-Witt vectors. In particular, if $k=\mathbb{F}_p$, then $W(k)=\mathbb{Z}_p$, the ring of $p$-adic integers, and $W_n(k)=\mathbb{Z}/p^n\mathbb{Z}$.

The ring $W_2(k)$ is described as follows. Let $\mathbb{Z}[k]$ be the group ring $\mathbb{Z}$-algebra, where we consider $k$ as an abelian group. It is a free $\mathbb{Z}$-module with basis $k$. Its elements have the form $\sum_i m_i\lambda_i $, $\lambda_i\in k$ and $m_i \in \mathbb{Z}$. In $\mathbb{Z}[k]$, for any prime $p>0$,  denote by $k\cdot p=\{\lambda p,\;\; \lambda \in k\}$. Then, as an additive group, $W_2(k)=k\oplus k\cdot p$ and the multiplication is defined by 
\[
(x+y\cdot p)\cdot (x^{\prime}+y^{\prime}\cdot p)=xx^{\prime}+(xy^{\prime}+x^{\prime}y)\cdot p,
\]
for any $x,y \in k$.

If $X$ admits a lift to characteristic zero over the $W(k)$ then many results of characteristic zero remain valid  in positive characteristic~\cite[Page 279]{Liedtke13}. In particular, the following is true.

\begin{theorem}[\cite{DeIll87}]
Let $X$ be a smooth projective variety defined over an algebraically closed field $k$ of characteristic $p>0$. Let $L$ be an ample invertible sheaf on $X$. Assume that both $X$ and $L$ admit a lifting to $W_2(k)$ and that $\dim X \leq p$. Then 
\[
H^i(X,L^{-1})=0,\;\; i< \dim X.
\]
\end{theorem}

In the case of surfaces of general type, T. Ekedahl~\cite{Ekedahl88} has obtained the following.

\begin{theorem}[\cite{Ekedahl88}]
Let $X$ be a minimal surface of general type and $L$ an invertible sheaf that is numerically equivalent to $\omega_X^{\otimes i}$, for some $i\geq 1$. Then $H^1(X,L^{-1})=0$ except possibly for certain surfaces in characteristic $2$ with $\chi(\mathcal{O}_X)\leq 1$.
\end{theorem}

\subsection{Counterexample of Kodaira Vanishing in positive characteristic}

Many of the known counterexamples of Kodaira vanishing for surfaces and especially the first ones constructed by Raynaud~\cite{Ray78} are based in the following observation of D. Mumford~\cite{Mu67} and L. Szpiro~\cite{Szp78}.  A clear exposition of it with applications can also be found in~\cite{Tak91}

\begin{theorem}[\cite{Mu67},~\cite{Szp78}]\label{Mumford criterion}
Let $f\colon X \rightarrow C$ be a fibration from a smooth projective surface $X$ onto a smooth projective curve $C$ such that all fibers are reduced,  irreducible and have positive arithmetic genus. Suppose that there exists a section  $\Gamma$  of $f$ such that $\Gamma^2>0$. Then
\begin{enumerate}
\item The sheaf $L=\mathcal{O}_X(\Gamma)\otimes f^{\ast}N$ is ample on $X$, where $N$ is an ample invertible sheaf on $C$ isomorphic to $\mathcal{O}_X(\Delta)\otimes \mathcal{O}_{\Gamma}$.
\item $H^1(X,L^{-1})\not=0$. In particular this is a counterexample to Kodaira Vanishing.
\end{enumerate}

\end{theorem}
Note that the previous theorem does not mention the characteristic of the base field anywhere. Considering that the Kodaira Vanishing theorem holds in characteristic zero, the theorem  implies then that the condition $\Gamma^2>0$ can only happen in positive characteristic.

We next define the notion of a Tango curve~\cite{Tango72},~\cite{Tak92},~\cite{Mukai11}.

\begin{definition}\label{def of Tango curve}
Let $C$ be a smooth projective curve defined over an algebraically closed field of characteristic $p>0$. $C$ is called a Tango curve of type $(p,n,d)$ if and only if the following holds.
\begin{enumerate}
\item There exists an invertible sheaf $\mathcal{L}$ on $C$ of degree $d\geq 1$ such that $\omega_C\cong \mathcal{L}^{pn}$.
\item There exists an open affine cover $U_i$, $i\in I$, of $C$ and sections $\xi_i \in \Gamma(U_i,\mathcal{O}_{U_i})$ such that $\Omega_{U_i}$ is generated by $d\xi_i$ as an $\mathcal{O}_{U_i}$-module and moreover on the intersection $U_{ij}=U_i\cap U_j$, $d\xi_i=a_{ij}^{np}d\xi_j$, where $a_{ij}\in\Gamma (U_{ij},\mathcal{O}^{\ast}_{U_{ij}})$ are the transition functions defining $\mathcal{L}$.
\end{enumerate}
\end{definition}
Note that, in the notation of the previous definition, since on $U_{ij}$ $d(\xi_i-a_{ij}^{np}\xi_j)=0$, there exist $b_{ij}\in \Gamma(U_{ij},\mathcal{O}_{U_{ij}})$, such that 
\begin{equation}\label{eq	200}
\xi_i=a_{ij}^{np}+b_{ij}^p.
\end{equation}

In constructing Tango curves the following criterion obtained by H. Kurke~\cite{Kurke81} is useful.

\begin{proposition}[\cite{Kurke81},~\cite{Tak91}]\label{Tango criterion}
Let $C$ be a smooth projective curve defined over an algebraically closed field $k$ of characteristic $p>0$. Let $\omega$ be an exact rational differential form of $C$ such that $\mathrm{div}(\omega)=pnD$, where $\mathrm{div}(\omega)$ is the divisor of $\omega$, $D$ is a non zero effective divisor of degree $d>0$ and $n$ is not divisible by $p$. Then $C$ is a Tango curve of type $(p,n,d)$.
\end{proposition}

\begin{example}[\cite{Mukai11},~\cite{Tak92}]
Let $C\subset \mathbb{P}^2_k$ be the curve given by the equation
\[
y^{np}+yz^{np-1}-x^{np-1}z=0
\]
where $k$ has characteristic $p>0$, $n\geq 1$ and $np\geq 2$. It is not hard to see that $C$ is smooth and by adjunction $p_a(C)=(np-1)(np-2)/2$. Hence $\deg \omega_C=2p_a(C)-2=np(np-3)$.

Let $U_z\subset \mathbb{P}^2_k$ be the affine open set where $z=1$. In $U_z$,  $C$ is given by the equation $y^{np}+y-x^{np-1}=0$. In this open set, $dy=-x^{np-2}dx$ and therefore on $U_z$   $\Omega_C$ is generated by $dx$. Hence the restriction of the divisor $\mathrm{div}(dx)$ of the exact rational differential form $dx$ on $U_z$ is zero. Hence $\mathrm{div}(dx)$ is supported on the line $L$ given by $z=0$. But the intersection $L\cap C$ is supported on the point $P=[1,0,0]$. Considering then $\deg(\mathrm{div}(dx))=\deg \omega_C$, we get that 
\[
\mathrm{div}(dx)=np(np-3)P.
\]
Therefore by Proposition~\ref{Tango criterion}, $C$ is a tango curve of type  $(p,n, np-3)$.
\end{example}

We are now ready to give an example of a smooth projective surface of general type where the Kodaira Vanishing Theorem fails. It is a Raynaud surface~\cite{Ray78} and the reader may find the details of the construction and relevant calculations in~\cite{La83},~\cite{Tak91} and~\cite{Tak92}.

\begin{example}[Counterexample to Kodaira Vanishing.]
Let $C$ be a Tango curve of type $(p,n,d)$ over an algebraically closed field $k$ of characteristic $p>0$ such that $n=pm-1$, $m\geq 2$. One can get such curves from the previous example by setting $n=pm-1$.  

By the definition of a Tango curve, there exists an affine open cover $U_i$, $i\in I$ of $C$, $\xi_i \in \Gamma (U_i,\mathcal{O}_{U_i})$, $a_{ij}\in \Gamma(U_{ij}, \mathcal{O}_{U_{ij}}^{\ast})$ satisfying the cocycle condition, such that on $U_{ij}$, 
$ d\xi_i=a^{p(mp-1)}_{ij}d \xi_j$, and therefore (as mentioned after Definition~\ref{def of Tango curve}), there exists $b_{ij}\in \Gamma(U_{ij}, \mathcal{O}_{U_{ij}} )$ such that
\begin{equation}\label{eq201}
\xi_i=a_{ij}^{p(mp-1)}\xi_j +b_{ij}^{p}.
\end{equation}

Let now $\tilde{\mathbb{P}}_i$ be the two dimensional weighted projective space with coordinates $x_i$, $y_i$, $z_i$ and weights $w(x_i)=m$, $w(y_i)=w(z_i)=1$. Let $V_i \subset \tilde{\mathbb{P}}_i \times U_i$ be the closed subscheme given by the graded homogeneous equation
\begin{equation}\label{eq205}
y_i^{mp-1}z_i=x_i^p+\xi_iz_i^{mp}.
\end{equation}
Let also $f_i \colon V_i \rightarrow U_i$ be the restriction of the projection $\tilde{\mathbb{P}}_i\times U_i \rightarrow U_i$ on $V_i$. Let $U_{ij}=U_i\cap U_j\subset U_i$, for any $i, j\in I$, $i\not=j$. Then $V_{ij}=f_i^{-1}(U_{ij})\subset \tilde{\mathbb{P}}_i \times U_{ij}$ is given again by the equation (\ref{eq205}) but now $\xi_i \in \Gamma( \mathcal{O}_{U_{ij}})$. 

Let $ \phi_{ij} \colon \tilde{\mathbb{P}}_i \times U_{ij} \rightarrow \tilde{\mathbb{P}}_j \times U_{ji}$ be defined by setting $\phi_{ij}([x_i,y_i,z_i])=[x_j,y_j,z_j]$, where
\begin{eqnarray}\label{eqn350}
x_j=a_{ji}^{mp-1}x_i-b_{ji}z_i^m\\\nonumber
y_j=a_{ji}^py_i\\
z_j=z_i. \nonumber
\end{eqnarray}
It is now easy to see, by taking into consideration the equation (\ref{eq201}), that
\[
y_j^{mp-1}z_j-x_j^p-\xi_jz_j^{np}=a_{ji}^{p(pm-1)}(y_i^{mp-1}z_i-x_i^p-\xi_iz_j^{np})
\]
Hence $\phi_{ij}(V_{ij}) \subset V_{ji}$ and hence $\phi_{ij}$ is an isomorphism between $V_{ij}$ and $V_{ji}$ over $U_{ij}$. It is now easy to check, by using the fact that the $a_{ij}$ satisfy the cocycle condition and using the equation (\ref{eq201}), that on the triple intersection $U_{ijk}=U_i\cap U_j \cap U_k$, $\phi_{jk}\phi_{ij}=\phi_{ik}$. Therefore the maps $f_i \colon V_i \rightarrow U_i$ glue to a map $f \colon X \rightarrow C$.  This surface is called a generalized Raynaud surface. This surface has the following properties.

\begin{proposition}[\cite{La83},~\cite{Tak91},~\cite{Tak92}]\label{Raynaud}
Let $f\colon X\rightarrow C$ be the generalized Raynaud surface constructed above, where $C$ is a Tango curve of type $(p,n,d)$, where $p=mp-1$. Then the following hold.
\begin{enumerate}
\item $X$ is smooth. Moreover,
\[
K_X^2=dp(m^2p^3-2mp^2-m(m+4)p+2m+4) \quadand c_2(X)=2pd(1-mp).
\]
\item Every fiber of $f$ is reduced and irreducible and has arithmetic genus $\frac{1}{2}(m-1)p(p-1)+\frac{1}{2}(p-1)(p-2)$.
\item Every fiber of $f$ is singular and its singular locus consists of exactly one point locally isomorphic to the cusp given by the equation $y^{mp-1}-x^p=0$.
\item There exists a section $S$ of $f$ such that $S^2=d>0$.
\item Let $K_X$ be the canonical divisor of $X$ and $F$ be a section. Then 
\[
K_X\equiv (mp+p-1)F+(m(p^2-p)-2p)S.
\]
where "$\equiv$" denotes numerical equivalence.
\item $H^0(X,T_X)=H^0(C,L)$. 
\end{enumerate}
\end{proposition}
\begin{proof}
We will only prove the smoothness of $X$, the existence of the section $S$ of $f$ and the fact that all fibers of $f$ are singular.  Detailed proofs of the remaining parts of (1)-(5) can be found in~\cite{La83}. The proof of (6) can be found in~\cite{Tak92}.

Let $W_{x_i}$, $W_{y_i}$, $W_{z_i}$ be the affine open cover of the weighted projective space $\tilde{\mathbb{P}}_i$ given by $x_i\not=0$, $y_i\not=0$ and $z_i\not=0$, respectively. From the equation (\ref{eq205}) which defines $V_i\subset \tilde{\mathbb{P}}_i$ it is clear that it is not possible to have $z_i=y_i=0$ on $V_i$. Therefore, $V_i\subset W_{y_i}\cup W_{z_i}$. Hence in order to prove that $X$ is smooth, it suffices to show that $V_i\cap W_{y_i}$ and $V_i \cap W_{z_i}$ are smooth. In particular, $V_i$ lies in the smoothh part of $\tilde{\mathbb{P}}_i$. This makes calculations easy.

Let us describe $V_i \cap W_{y_i}$. $W_{y_i}$ is isomorphic to $\mathbb{A}^2_k$ with coordinates $\tilde{x}_i=x_i/y_i^n$ and $\tilde{z}_i=z_i/y_i$. Therefore,  from (\ref{eq205}), $V_i \cap W_{y_i}$ is the closed subscheme of $\mathbb{A}^2_k \times U_i$ given by the equation
\[
\tilde{z}_i=\tilde{x}_i^p+\xi_i \tilde{z}_i^{mp}.
\]
Since $U_i$ is smooth, this is also easily seen to be smooth.

Let us now describe $V_i\cap W_{z_i}$. $W_{z_i}$ is isomorphic to $\mathbb{A}^2_k$ with coordinates $\tilde{x}_i=x_i/z_i^n$ and $\tilde{y}_i=y_i/z_i$. Therefore,  from (\ref{eq205}), $V_i \cap W_{z_i}$ is the closed subscheme of $\mathbb{A}^2_k \times U_i$ given by the equation
\begin{equation}\label{eq300}
\tilde{y}_i^{mp-1}=\tilde{x}_i^p+\xi_i.
\end{equation}
We will show that this is smooth. By completing $U_i$ at any closed point $P\in U_i$, and since $\hat{\mathcal{O}}_{U_i,P}\cong k[[t]]$, it follows that over $\hat{\mathcal{O}}_{U_i,P}$, $V_i\cap W_{z_i}$ is given by the equation (\ref{eq300}) above, but in 
$k[[t]][x_i,y_i,z_i]$. Moreover, $\xi_i=\xi_i(t)$. If this is singular, then  $\xi_i^{\prime}(t)=0$. Therefore, $\xi_i=\xi_i(t^p)$. But then $d\xi_i=0$. This contradicts the fact that $d\xi_i$ generate locally $\Omega_C$ and in particular this is not zero. Therefore $V_i\cap W_{z_i}$ is also smooth and hence $X$ is smooth.

Next we will show that every fiber of $f$ is singular and will describe its singularities. Let $P\in C$ and let $X_P=f^{-1}(P)$ be the fiber of $f$ over $P$. Then $X_P=(X_P\cap W_{y_i})\cup (X_P\cap W_{z_i})$. $X_P\cap W_{z_i}$ is the affine plane curve given by the equation
\[
\tilde{y}_i^{mp-1}=\tilde{x}_i^p+\bar{\xi}_i,
\]
in $k[x_i,y_i,z_i]$, where $\bar{\xi}_i$ is the image of $\xi_i$ in $k=k(P)$, the residue field of $P\in C$. Since $k$ is algebraically closed, this defines a singularity locally isomorphic to the cusp $\tilde{y}_i^{mp-1}=\tilde{x}_i^p$ at the point  $(-\bar{\xi}_i^{1/p}, 0)$. Therefore every fiber of $f$ is singular and its singular locus is one cusp singularity as above.

Finally we will prove the existence of a section $S$ of $f$. Let $S_i\subset \tilde{\mathbb{P}}_i \times U_i$ be defined by $z_i=x_i=0$, i.e., $S_i=[0,1,0]\times U_i$. Then $S_i\subset V_i$ and a section of $f_i \colon V_i \rightarrow U_i$. Now it is not difficult to see that the $S_i$ glue, by the glueing data (\ref{eqn350}), to a curve $S \subset X$ which is a setion of $f$.

\end{proof}
Let $X$ be a surface of general type over an algebraically closed field of characteristic zero. Then $c_1^2(X)\leq 3c_2(X)$ and $c_2(X)>0$. The first inequality is the Bogomolov-Miyaoka-Yau inequality and the second is the Castelnuovo inequality. As a corollary of the previous proposition we get that the generalized Raynaud surface $X$ constructer above violates several properties that hold in characteristic zero. In particular the Kodaira Vanishing Theorem, the Bogomolov-Miyaoka-Yau and the Castelnuovo inequality.

\begin{corollary}\label{counterexamples}
With notation and assumptions as in the previous theorem.
\begin{enumerate}
\item $K_X$ is ample, and hence $X$ canonically polarized, unless $(p,m)=(3,1)$ or $(p,m)=(2,2)$.
\item Let $N=\mathcal{O}_X(S)\otimes f^{\ast} M$, where $M$ is an invertible sheaf on $C$ isomorphic to $\mathcal{O}_X(S)\otimes \mathcal{O}_S$. Then $N$ is ample and 
\[
H^1(X,N^{-1})\not= 0.
\]
Hence this is a counterexample to the Kodaira Vanishing Theorem.
\item Suppose that $(p,m)\not=(3,1)$ or $(p,m)\not=(2,2)$. Then 
\[
c_2(X)<0.
\]
 Hence this is a counterexample to both the Bogomolov-Miyaoka-Yau inequality and the Castelnuovo inequality.
\item The automorphism group scheme $\mathrm{Aut}_{X/k}$ is not smooth.
\end{enumerate}
\end{corollary}
\begin{proof}
The statement (\ref{counterexamples}.1) follows easily from Proposition~\ref{Raynaud}.1 and the Nakai-Moishezon criterion for ampleness. 

The statement (\ref{counterexamples}.2) follows from Proposition~\ref{Raynaud}.4 and Theorem~\ref{Mumford criterion}.

The statement (\ref{counterexamples}.3) is clear since from the first part, $c_1^2>0$ and $c_2<0$.

Finally, regarding the last statement.  Since $X$ is of general type, $\mathrm{Aut}_{X/k}$ is a finite group scheme~\cite{LeMe78} and its tangent space at the identity is  $H^0(X,T_X)$. Hence $\mathrm{Aut}_{X/k}$ is smooth if and only if $H^0(X,T_X)=0$. However, by Proposition~\ref{Raynaud}.6, $H^0(X,T_X)=H^0(C,L)\not=0$, Hence $\mathrm{Aut}_{X/k}$ is not smooth.

\end{proof}

\end{example}

% !TEX root = paper.tex

%\input{intro}
%\input{sec1}
%\input{bib}
\end{document}